\documentclass[11pt]{article}

\usepackage{amssymb}
\usepackage{amsmath}
\usepackage{theorem}
\usepackage{epsfig}
\usepackage{verbatim}
\usepackage{graphicx}
\usepackage{indentfirst}
\newtheorem{thm}{Theorem}[section]

\newtheorem{prop}[thm]{Proposition}

\newtheorem{rem}[thm]{Remark}

\newtheorem{lemma}[thm]{Lemma}

\newcommand{\R}{\Bbb{R}}

\newcommand{\Z}{\Bbb{Z}}

\newcommand{\T}{\mathbb{T}}
\newcommand{\D}{\displaystyle}

\newcommand{\sign}{{\rm sign}\thinspace}

\newcommand{\di}{{\rm div}\thinspace}

\newcommand{\dt}{\frac{d}{dt}}
\newcommand{\F}{\mathcal{F}}

\newcommand{\dxd}{\partial_{x_2}}

\newcommand{\dpt}{\partial_t}

\newcommand{\la}{\Lambda}

\newcommand{\al}{\alpha}
\newcommand{\ep}{\varepsilon}

\newcommand{\dpa}{\partial^{\bot}_{\alpha}}
\newcommand{\da}{\partial_{\alpha}}

\numberwithin{equation}{section}

\begin{document}

\author{Antonio C\'ordoba, Diego C\'ordoba and Francisco Gancedo}

\title{Interface evolution: the Hele-Shaw and Muskat problems.}
\date{23-5-08}

\maketitle

\begin{abstract}
We study the dynamics of the interface between two incompressible
2-D flows where the evolution equation is obtained from Darcy's
law. The free boundary is given by the discontinuity  among the
densities and viscosities of the fluids. This physical scenario is
known as the two dimensional Muskat problem or the two-phase
Hele-Shaw flow. We prove local-existence in Sobolev spaces when,
initially, the difference of the gradients of the pressure in the
normal direction has the proper sign, an assumption which is also
known as the Rayleigh-Taylor condition.
\end{abstract}

\maketitle

%%%%%%%%%%%%%%%%%%%%%%%%%%%%%%%%%%%%%%%%%%%%%%%%%%%%%%%%%%%%%%%%%%%%%%%%%%%
\section{Introduction}

We consider the following evolution problem for the active scalar
$\rho = \rho(x,t)$, $x\in\mathbb{R}^2$, and $t\geq 0$:
\begin{equation*}
\rho_t + v\cdot\nabla\rho = 0,
\end{equation*}
with a velocity $v = (v_1,v_2)$ satisfying the momentum equation
\begin{equation}\label{dlaw}
\frac{\mu}{\kappa} v=-\nabla p-(0,\mathrm{g}\,\rho),
\end{equation}
and the incompressibility condition $\nabla\cdot v = 0$.

In the following we achieve a rather complete local existence
analysis of the dynamics of the interface between two incompressible
2-D flows with different characteristics (i.e. distinct values of
$\mu$ and $\rho$) which are  evolving under \eqref{dlaw}, also known
as Darcy's law \cite{bear}. This system was studied by Muskat \cite{Muskat} in
order to model the interface between two fluids in a porous media,
where $p$ is the pressure, $\mu$ is the dynamic viscosity, $\kappa$
is the permeability of the medium, $\rho$ is the liquid density and
$\mathrm{g}$ is the acceleration due to gravity. Saffman and Taylor
\cite{S-T} made the observation that the one phase version (one of
the fluids has zero viscosity) was  also known as the Hele-Shaw cell
equation \cite{H-S}, which, in turn, is the zero-specific heat case of the
classical one-phase Stefan problem.

There is a vast literature about those problems (see \cite{Peter}
and \cite{Hou} for references). In order to frame our result let
us point out that in \cite{SCH} is treated the case where both
densities are equal, showing global existence for small data in
the stable case and ill-possednes in the unstable case. In
\cite{Ambrose} the well-possednes in the stable case was
considered under time dependent assumption of the arc-chord
condition. Finally, in the case where the viscosities are the same,
the character of the interphase as the graph of a function is
preserved and in \cite{DY} \cite{DY2} this fact has been used to
prove local existence and a maximum principle, in the stable case,
together with ill-possednes in the unstable situation.

Due to the direction of gravity, the horizontal and the vertical
coordinates play different rolls. Here we shall assume spatial
periodicity in the horizontal space variable, says  $\rho(x_1 +
2k\pi, x_2, t) = \rho(x_1,x_2,t)$. The free boundary is given by
the discontinuity on the densities and viscosities of the fluids,
where $(\mu,\rho)$ are defined by

\begin{equation}\label{density}
(\mu,\rho)(x_1,x_2,t)=\left\{\begin{array}{cl}
                    (\mu^1,\rho^1),& x\in\Omega^1(t)\\
                    (\mu^2,\rho^2),& x\in\Omega^2(t)=\mathbb{R}^2 - \Omega^1(t),
                 \end{array}\right.
\end{equation} and $\mu^1\neq\mu^2$, and $\rho^1\neq\rho^2$ are constants.

Let the free boundary be parameterized by
$$\partial\Omega^j (t) =
\{z(\alpha, t) = (z_1(\alpha,t), z_2(\alpha,t)) : \alpha \in \R\}
$$
such that $$(z_1(\alpha+2k\pi,t),
z_2(\alpha+2k\pi,t))=(z_1(\alpha,t)+2k\pi, z_2(\alpha,t)),$$ with
the initial data $z(\alpha, 0) = z_0(\alpha)$.

 Notice that each fluid is irrotational, i.e. $\omega = \nabla \times u
 = 0$, in the interior of each domain $\Omega^i$ ($i=1,2$). Therefore the
 vorticity $\omega$ has its support on the curve $z(\al,t)$ and it
 can be shown easily to be of the form
\begin{equation*}
\omega(x,t)=\varpi(\al,t)\delta(x-z(\al,t)).
\end{equation*}

Then $z(\alpha, t)$ evolves with a velocity field coming from
Biot-Savart law, which can be explicitly computed and is given by
the Birkhoff-Rott integral of the amplitude $\varpi$ along the
interface curve:

\begin{align}
\begin{split}\label{fibr}
BR(z,\varpi)(\al,t)=\big(\!-\!\frac{1}{4\pi}
PV\int_{\T}\varpi(\beta,t)\frac{\tanh(\frac{z_2(\al,t)-z_2(\beta,t)}{2})
(1+\tan^2(\frac{z_1(\al,t)-z_1(\beta,t)}{2}))}{\tan^2(\frac{z_1(\al,t)-
z_1(\beta,t)}{2})+\tanh^2(\frac{z_2(\al,t)-z_2(\beta,t)}{2})}d\beta,\\
\frac{1}{4\pi}
PV\int_{\T}\varpi(\beta,t)\frac{\tan(\frac{z_1(\al,t)-z_1(\beta,t)}{2})
(1-\tanh^2(\frac{z_2(\al,t)-z_2(\beta,t)}{2}))}{\tan^2(\frac{z_1(\al,t)-
z_1(\beta,t)}{2})+\tanh^2(\frac{z_2(\al,t)-z_2(\beta,t)}{2})}d\beta\big),
\end{split}
\end{align}
where $PV$ denotes principal value \cite{St3}. It gives us the velocity field at the interface to which we can
subtract any term in the tangential direction without modifying the
geometric evolution of the curve
\begin{align}
\begin{split}\label{fullpm}
z_t(\al,t)&=BR(z,\varpi)(\al,t)+c(\al,t)\da z(\al,t).
\end{split}
\end{align}
A wise choice of $c(\al,t)$ namely:
\begin{align}
\begin{split}\label{fc}
c(\al,t)&=\frac{\al+\pi}{2\pi}\int_\T\frac{\da z(\al,t)}{|\da
z(\al,t)|^2}\cdot \da BR(z,\varpi)(\al,t) d\al\\
&\quad-\int_{-\pi}^\al
\frac{\da z(\beta,t)}{|\da z(\beta,t)|^2}\cdot \partial_{\beta}
BR(z,\varpi)(\beta,t)d\beta,
\end{split}
\end{align}
allows us to accomplish the fact that the length of the tangent
vector to $z(\al,t)$ be just a function in the variable $t$ only:
$$
A(t)=|\da z(\al,t)|^2,
$$
as will be shown in section 2 (see also \cite{Hou} and \cite{Y}). Then we can close the system using
Darcy's law with the equation:
\begin{align}
\begin{split}\label{fw}
\varpi(\al,t)&=-2A_\mu BR(z,\varpi)(\al,t)\cdot
\da z(\al,t)-2\kappa \mathrm{g}\frac{\rho^2-\rho^1}{\mu^2+\mu^1}\da
z_2(\al,t),
\end{split}
\end{align}
where $$A_\mu=\frac{\mu_1-\mu_2}{\mu_1+\mu_2}$$ is the Atwood number.

Finally we give the function which measures the arc-chord condition in the periodic case

\begin{equation}\label{df}
\mathcal{F}(z)(\alpha,\beta,t)=\frac{\beta^2/4}{\tan^2(\frac{z_1(\alpha,t)-z_1(\alpha-\beta,t)}{2})+
\tanh^2(\frac{z_2(\alpha,t)-z_2(\alpha-\beta,t)}{2})}\qquad
\forall\,\alpha,\beta\in(-\pi,\pi),
\end{equation}
with $$\mathcal{F}(z)(\alpha,0,t)=\frac{1}{|\da z(\alpha,t)|^2},$$
(see \cite{Y} for a closed curve).

Our main result consists on the existence of a positive time $T$
(depending upon the initial condition) for which we have a
solution of the periodic Muskat problem (equations \eqref{fibr}-\eqref{fw}) during the  time interval
$[0,T]$ so long as the initial data satisfy $z_0(\alpha)\in
H^{k}(\T)$ for $k\geq 3$, $\F(z_0)(\al,\beta)<\infty$, and
$$\sigma_0(\al)= -(\nabla p^2(z_0(\alpha)) - \nabla
p^1(z_0(\alpha)))\cdot\partial^{\perp}_{\alpha}z_0(\alpha) > 0,$$
where $p^j$ denote the pressure in $\Omega^j$.

It is interesting to remark that the equality of pressure at each
side of the free boundary is obtained in section 2 directly from
Darcy's law without any other assumption.

\begin{thm}\label{theorem2}
Let $z_0(\alpha)\in H^{k}(\T)$ for $k\geq 3$,
$\F(z_0)(\al,\beta)<\infty$, and $$\sigma_0(\al)=-(\nabla
p^2(z_0(\alpha)) - \nabla
p^1(z_0(\alpha)))\cdot\partial^{\perp}_{\alpha}z_0(\alpha)
> 0.$$ Then there exists a time $T>0$ so that there is a solution to \eqref{fibr}-\eqref{fw} in $C^{1}([0,T];H^k(\T))$ with
$z(\al,0)=z_0(\al)$.
\end{thm}

We devote the rest of the paper to the proof of theorem \ref{theorem2}
 which is organized as follows. In section 2 we derive the system of equations \eqref{fibr}-\eqref{fw}
with the corresponding choice of $c(\al,t)$ and we also obtain the
properties of the pressure. In section 3 and 4 we present several
crucial estimates on the operator $T(u)(\al)=2BR(z,u)(\al)\cdot
\da z(\al)$ and on the inverse operator $(I-\xi T)^{-1}$,
$|\xi|\leq 1$. Our proofs rely upon the boundedness properties of
the Hilbert transforms associated to $C^{1,\alpha}$ curves, for
which we need precise estimates obtained with arguments involving
conformal mappings, Hopf maximum principle and Harnack
inequalities. We then provide upper bounds  for  the amplitude of
the vorticity, the Birkhoff-Rott integral, the parametrization of
the curve and the arc-chord condition, namely:
\begin{equation*}
\|\varpi\|_{H^k}\leq \exp C(\|\F(z)\|^2_{L^\infty}+\|z\|^2_{H^{k+1}})\quad\text{(section 5)},
\end{equation*}
\begin{eqnarray*}
\|BR(z,\varpi)\|_{H^k}\leq \exp(C(\|\F(z)\|^2_{L^\infty}+\|z\|^2_{H^{k+1}})\quad\text{(section 6)},
\end{eqnarray*}
\begin{align*}
\frac{d}{dt}\|z\|^2_{H^k}(t)&\leq \D-\frac{\kappa}{2\pi(\mu_1\!+\!\mu_2)}\,\int_\T \frac{\sigma(\al,t)}{|\da z(\al)|^2} \da^k z(\al,t)\cdot \la(\da^k z)(\al,t) d\al\quad\\
&\quad+\exp C(\|\F(z)\|^2_{L^\infty}(t)+\|z\|^2_{H^k})\qquad\qquad\qquad\qquad\qquad \text{(section 7)},
\end{align*}
and
\begin{equation*}
\D\dt\|\F(z)\|^2_{L^\infty}(t)\leq \exp C(\|\F(z)\|^2_{L^\infty}(t)+\|z\|^2_{H^3}(t))\quad\text{(section 8),}
\end{equation*}
where the operator $\la$ is defined by the Fourier transform $\widehat{\la f}(\xi)=|\xi|\widehat{f}(\xi)$ and $\sigma(\al,t)$
is the difference of the gradients of the pressure in the
normal direction. In section 9 we study the evolution of $m(t)=\D\min_{\al\in\T} \sigma(\al,t)$, which satisfies
the following lower bound
$$m(t)\geq m(0)-\int_0^t \exp C(\|\F(z)\|^2_{L^\infty}(s)+\|z\|^2_{H^3}(s)) ds.$$
Finally, in section 10, we introduce a regularized evolution
equation where we use the previous a priori estimates together
with a pointwise inequality satisfied by the non-local operator
$\Lambda$ \cite{CC} to show local existence.

By a similar approach, in \cite{DY3} we obtain local-existence in Sobolev spaces for the full water wave problem in 2-D 
when the Rayleigh-Taylor condition is initially satisfied.
%%%%%%%%%%%%%%%%%%%%%%%%%%%%%%%%%%%%%%%%%%%%%%%%%%%%%%%%%%%%%%%%%%%%%%%%%%%%%%%%%%%%%%%%%%%%%%%%%%%%%%%%%%%%

%%%%%%%%%%%%%%%%%%%%%%%%%%%%%%%%%%%%%%%%%%%%%%%%%%%%%%%%%%%%%%%%%%%%%%%%%%%%%%%%%%%%%%%%%%%%%%%%%%%%%%%%%%%%

\section{The evolution equation}

Here $(\mu,\rho)$ are defined by

\begin{equation*}%\label{density}
(\mu,\rho)(x_1,x_2,t)=\left\{\begin{array}{cl}
                    (\mu^1,\rho^1),& x\in\Omega^1(t)\\
                    (\mu^2,\rho^2),& x\in\Omega^2(t),
                 \end{array}\right.
\end{equation*} where $\mu^1\neq\mu^2$, and $\rho^1\neq\rho^2$. Then using the Biot-Savart law we get

$$
v(x,t)=\frac{1}{2\pi}PV\int_\R\frac{(x-z(\beta,t))^\bot}{|x-z(\beta,t)|^2}
\varpi(\beta,t)d\beta
$$
for $x\neq z(\al,t)$ where the principal value is taken at
infinity.

It is convenient to introduce the complex notation  $z=x_1+ix_2$,
then the complex conjugate $\bar{v}$ of the velocity field is
given by
$$
\bar{v}(z,t)=\frac{1}{2\pi
i}PV\int_\R\frac{\varpi(\beta,t)}{z-z(\beta,t)} d\beta.
$$
In our case of  periodic interface, $z(\al+2\pi k,t)=z(\al,t)+2\pi
k$, the following classical identity
$$
\frac{1}{\pi}\big(\frac1z+\sum_{k\geq 1}\frac{2z}{z^2-(2\pi
k)^2}\big)=\frac{1}{2\pi\tan(z/2)},
$$
yields
\begin{align*}
\begin{split}\label{BS}
v(x,t)=\big(-\frac{1}{4\pi}
\int_{\T}\varpi(\beta,t)\frac{\tanh(\frac{x_2-z_2(\beta,t)}{2})
(1+\tan^2(\frac{x_1-z_1(\beta,t)}{2}))}{\tan^2(\frac{x_1-
z_1(\beta,t)}{2})+\tanh^2(\frac{x_2-z_2(\beta,t)}{2})}d\beta,\\
\frac{1}{4\pi}
\int_{\T}\varpi(\beta,t)\frac{\tan(\frac{x_1-z_1(\beta,t)}{2})
(1-\tanh^2(\frac{x_2-z_2(\beta,t)}{2}))}{\tan^2(\frac{x_1-
z_1(\beta,t)}{2})+\tanh^2(\frac{x_2-z_2(\beta,t)}{2})}d\beta\big)
\end{split}
\end{align*}
for $x\neq z(\al,t)$.

We have

\begin{equation*}
v^2(z(\al,t),t)=BR(z,\varpi)(\al,t)+\frac12\frac{\varpi(\al,t)}{|\da
z(\al,t)|^2}\da z(\al,t),
\end{equation*}

\begin{equation*}
v^1(z(\al,t),t)=BR(z,\varpi)(\al,t)-\frac12\frac{\varpi(\al,t)}{|\da
z(\al,t)|^2}\da z(\al,t),
\end{equation*}
where $v^j(z(\alpha,t),t)$ denotes the limit velocity field obtained
approaching the boundary in the normal direction inside $\Omega^j$ and
$BR(z,\varpi)(\al,t)$ is given by \eqref{fibr}.

 Darcy's law implies
\begin{equation*}
\Delta p(x,t)=-\di (\frac{\mu(x,t)}{\kappa} v(x,t))-\mathrm{g}\,\dxd
\rho(x,t),
\end{equation*}
therefore
\begin{equation*}
\Delta p(x,t)=\Pi(\al,t)\delta (x-z(\al,t)),
\end{equation*}
where $\Pi(\al,t)$ is given by
$$\Pi(\al,t)=\frac{\mu^2-\mu^1}{\kappa}v(z(\al,t),t)\cdot\dpa
z(\al,t)+\mathrm{g}(\rho^2-\rho^1)\da z_1 (\al,t).$$
It follows that:
\begin{equation*}
 p(x,t)=-\frac{1}{2\pi}\int_\T \ln\big(\cosh(x_2-z_2(\al,t))-\cos(x_1-z_1(\al,t))\big)\Pi(\al,t) d\al,
\end{equation*}
for $x\neq z(\al,t)$,
 implying the important identity
\begin{equation*}
p^2(z(\al,t),t)=p^1(z(\al,t),t),
\end{equation*}
which is just a mathematical consequence of Darcy's law, making
unnecessary to impose it as a physical assumption.

Let us introduce the following notation:
$$
[\mu v](\al,t)=(\mu^2 v^2(z(\al,t),t)-\mu^1 v^1(z(\al,t),t))\cdot \da z(\al,t).
$$
Then taking the limit in Darcy's law we obtain
\begin{align*}
\begin{split}
\frac{[\mu v](\al,t)}{\kappa}&=-(\nabla p^2(z(\al,t),t)-\nabla
p^1(z^1(\al,t),t))\cdot \da z(\al,t)-\mathrm{g}(\rho^2-\rho^1)\,\da z_2(\al,t)\\
&=-\da (p^2(z(\al,t),t)-p^1(z(\al,t),t))-\mathrm{g}(\rho^2-\rho^1)\,\da z_2(\al,t)\\
&=-\mathrm{g}(\rho^2-\rho^1)\,\da z_2(\al,t),
\end{split}
\end{align*}
which gives us
\begin{align*}
\begin{split}
\frac{\mu^2+\mu^1}{2\kappa}\varpi(\al,t)+\frac{\mu^2-\mu^1}{\kappa}BR(z,\varpi)(\al,t)\cdot
\da z(\al,t)=-\mathrm{g}(\rho^2-\rho^1)\da z_2(\al,t),
\end{split}
\end{align*}
so that
\begin{align*}
\begin{split}
\varpi(\al,t)=-A_{\mu}2BR(z,\varpi)(\al,t)\cdot
\da z(\al,t)-2\kappa\mathrm{g}\frac{\rho^2-\rho^1}{\mu^2+\mu^1}\da
z_2(\al,t).
\end{split}
\end{align*}

Next we modify the velocity of the curve in the tangential
direction:

\begin{align}\label{eev}
\begin{split}
z_t(\al,t)=BR(z,\varpi)(\al,t)+c(\al,t)\da z(\al,t),
\end{split}
\end{align}
where the scalar $c(\al,t)$ is chosen in such a way that the
tangent vector only depends on the variable $t$ as follows:
\begin{equation}\label{cancelacionextra}
|\da z(\al,t)|^2=A(t).
\end{equation}

To find such a $c(\al,t)$ let us differentiate the identity
\eqref{cancelacionextra}
\begin{align*}
\begin{split}
A'(t)&=2\da z(\al,t)\cdot \da z_t(\al,t)=2\da z(\al,t)\cdot \da
BR(z,\varpi)(\al,t)+2\da c(\al,t) A(t),
\end{split}
\end{align*}
so that
\begin{align}
\begin{split}\label{dla}
\da c(\al,t)=\frac{A'(t)}{2A(t)}-\frac{1}{A(t)}\da z(\al,t)\cdot
\da BR(z,\varpi)(\al,t).
\end{split}
\end{align}
Because $c(\al,t)$ has to be periodic, we obtain
\begin{align}
\begin{split}\label{AppA}
\frac{A'(t)}{2A(t)}=\frac{1}{2\pi A(t)}\int_\T\da z(\al,t)\cdot
\da BR(z,\varpi)(\al,t) d\al.
\end{split}
\end{align}

Using \eqref{AppA} in \eqref{dla}, and integrating in $\al$, one
gets the following formula
\begin{align}
\begin{split}\label{fc}
c(\al,t)&=\frac{\al+\pi}{2\pi}\int_\T\frac{\partial_{\beta}
z(\beta,t)}{|\partial_{\beta}
z(\beta,t)|^2}\cdot \partial_{\beta} BR(z,\varpi)(\beta,t) d\beta\\
&\quad -\int_{-\pi}^\al \frac{\partial_{\beta}
z(\beta,t)}{|\partial_{\beta} z(\beta,t)|^2}\cdot
\partial_{\beta} BR(z,\varpi)(\beta,t) d\beta,
\end{split}
\end{align}
where we have chosen $c(-\pi,t)=c(\pi,t)=0$.

Let us consider now the solutions of  equation \eqref{eev} with
$c(\al,t)$ given by \eqref{fc}. It is easy to check that
$$
\dt |\da z(\al,t)|^2=c(\al,t)\da |\da z(\al,t)|^2+b(t)|\da
z(\al,t)|^2,
$$
where
$$b(t)=\frac{1}{\pi}\int_\T\frac{\partial_{\beta} z(\beta,t)}{|\partial_{\beta}
z(\beta,t)|^2}\cdot\partial_\beta  BR(z,\varpi)(\beta,t) d\beta.
$$ Next we solve this linear partial differential equation, assuming that \eqref{cancelacionextra} is satisfied
initially, to find that the unique solution is given by
\begin{align*}
\begin{split}
|\da z(\al,t)|^2=|\da z(\al,0)|^{2}+\frac{1}{\pi}\int_0^t\int_\T \da
z(\al,t)\cdot\partial_\beta BR(z,\varpi)(\al,t) d\al ds,
\end{split}
\end{align*}
which proves \eqref{cancelacionextra}.

Our next step is to find the formula for the difference of the
gradients of the pressure in the normal direction:
$$-(\nabla p^2(z(\alpha,t),t) - \nabla
p^1(z(\alpha,t),t))\cdot\partial^{\perp}_{\alpha}z(\alpha,t),$$
which we denote by $\sigma(\al,t)$. Approaching the boundary in Darcy's law, we get
$$
\sigma(\al,t)=\frac{\mu^2-\mu^1}{\kappa}BR(z,\varpi)(\al,t)\cdot\dpa
z(\al,t)+ \mathrm{g}(\rho^2-\rho^1)\da z_1(\al,t).
$$
It is easy to check that
\begin{align*}
\begin{split}
\frac{\mu^2-\mu^1}{\kappa}BR(z,\varpi)(\al,t)\cdot\dpa z(\al,t)=
\frac{1}{4\pi} \da\int_\T \varpi(\beta,t)\log G(\al,\beta,t) d\beta,
\end{split}
\end{align*}
with
\begin{align*}
G(\al,\beta,t)&=\sin^2(\frac{z_1(\al,t)-z_1(\beta,t)}{2})\cosh^2(\frac{z_2(\al,t)-z_2(\beta,t)}{2})\\
&\quad+\cos^2(\frac{z_1(\al,t)-z_1(\beta,t)}{2})\sinh^2(\frac{z_2(\al,t)-z_2(\beta,t)}{2}),
\end{align*}
and therefore
$$
\int_\T \frac{\mu^2-\mu^1}{\kappa}BR(z,\varpi)(\al,t)\cdot\dpa
z(\al,t)d\alpha =0.
$$
This shows that the condition $\rho^2\neq \rho^1$ is crucial in
order to have a constant  sign in the normal direction of the
difference of the gradient. Furthermore, since $z_1(\al,t)-\al$ is
periodic we have
$$
\int_{\T}\da z_1(\al,t) d\al =2\pi.
$$
\begin{rem}
If we consider a closed contour, then it is easy to check that
$$
\int_{\T} \sigma(\al,t) d\al= 0,
$$
which makes  impossible the task of prescribing a sign to $\sigma$
along a closed curve.
\end{rem}
%%%%%%%%%%%%%%%%%%%%%%%%%%%%%%%%%%%%%%%%%%%%%%%%%%%%%%%%%%%%%%%%%%%%%%%%%%%%%%%%%%%%%%%%%%%%%%%%%%%%%%%%%%

%%%%%%%%%%%%%%%%%%%%%%%%%%%%%%%%%%%%%%%%%%%%%%%%%%%%%%%%%%%%%%%%%%%%%%%%%%%%%%%%%%%%%%%%%%%%%%%%%%%%%%%%%%

\section{The basic operator}

Let us consider the operator $T$ defined by the formula

\begin{equation}
T(u)(\al)=2BR(z,u)(\al)\cdot \da z(\al). \label{operatorT}
\end{equation}

\begin{lemma}
Suppose that $\|\F(z)\|_{L^{\infty}}< \infty$ and $z\in
C^{2,\delta}$ with $0<\delta$. Then $T:L^2\rightarrow H^1$ and
$$\|T\|_{L^2\rightarrow H^1}\leq \|\F(z)\|^2_{L^{\infty}}\|z\|^4_{C^{2,\delta}}.$$
\end{lemma}

\begin{rem}
In section 5, lemma
5.2. there is a proof showing that T also maps $H^k$ into $H^{k +
1}$, $k\geq 1$.
\end{rem}

Proof: Since the formula \eqref{fibr} yields
$$T(u)(\al)=\frac1\pi\da \int_\T u(\beta)\arctan \Big(\frac{\tanh(\frac{z_2(\al)-z_2(\beta)}{2})}{\tan(\frac{z_1(\al)-z_1(\beta)}{2})}\Big)d\beta,$$
we have
$$\int_\T T(u)(\al)d\al=0, $$
which implies $\|T(u)\|_{L^2}\leq \|\da T(u)\|_{L^2}$.

Let us denote
\begin{align*}
V(\al,\beta)&=(V_1(\al,\beta), V_2(\al,\beta))=(\tan(\frac{z_1(\al)-z_1(\beta)}{2}),\tanh(\frac{z_2(\al)-z_2(\beta)}{2})).
\end{align*} In the following we shall refer to the Appendix for the
definition of $V_j, A_j$ and their properties.

We write first: $$\da T(u)=2BR(z,u)(\al)\cdot \da^2
z(\al)+2\da z(\al) \cdot \da BR(z,u)(\al)=I_1+I_2.$$ For $I_1$ we
have the expression
\begin{align*}
\begin{split}
I_1&=2(BR(z,u)(\al)-\frac{\dpa z(\al)}{|\da z(\al)|^2} H(u)(\al))\cdot \da^2 z(\al)+2H(u)(\al)\frac{\dpa z(\al)\cdot \da^2 z(\al)}{|\da z(\al)|^2}\\
&=J_1+J_2,
\end{split}
\end{align*}
where $H(u)$ is the (periodic) Hilbert transform of the function u.

Then
\begin{align*}
\begin{split}
J_1&=-\frac{1}{2\pi}\da^2 z_1(\al)\int_{\T}u(\beta)\frac{V_2(\al,\beta)V^2_1(\al,\beta)}{|V(\al,\beta)|^2} d\beta\\
&\quad-\frac{1}{2\pi}\da^2 z_2(\al)\int_{\T}u(\beta)\frac{V_1(\al,\beta)V^2_2(\al,\beta)}{|V(\al,\beta)|^2} d\beta\\
&\quad-\frac{1}{2\pi}\da^2
z_1(\al)\int_{\T}u(\al-\beta)[\frac{V_2(\al,\al-\beta)}{|V(\al,\al-\beta)|^2}
-\frac{1}{|\da z(\al)|^2}\frac{\da z_2(\al)}{\tan(\beta/2)}] d\beta\\
&\quad+\frac{1}{2\pi}\da^2
z_2(\al)\int_{\T}u(\al-\beta)[\frac{V_1(\al,\al-\beta)}{|V(\al,\al-\beta)|^2}
-\frac{1}{|\da z(\al)|^2}\frac{\da z_1(\al)}{\tan(\beta/2)}] d\beta\\
= &\quad  K_1+K_2+K_3+K_4.
\end{split}
\end{align*}
And we may use  that $|V_2(\al,\beta)|\leq 1$, to get
$|K_1|+|K_2|\leq C\|z\|_{C^2}\|u\|_{L^2}$.

To estimate  $K_3$ let us observe that the following term
$$
A_1(\al,\al-\beta)=\frac{V_2(\al,\al-\beta)}{|V(\al,\al-\beta)|^2}
-\frac{1}{|\da z(\al)|^2}\frac{\da
z_2(\al)}{\tan(\frac{\beta}{2})}
$$ satisfies $\|A_1\|_{L^{\infty}}\leq \|\F(z)\|_{L^{\infty}}\|z\|^2_{C^2}$ (see appendix lemma \ref{A12}).

In $K_4$ we have the term
$$
A_2(\al,\al-\beta)=\frac{V_1(\al,\al-\beta)}{|V(\al,\beta)|^2}
-\frac{1}{|\da z(\al)|^2}\frac{\da
z_1(\al)}{\tan(\frac{\beta}{2})},
$$
which satisfies $\|A_2\|_{L^\infty}\leq
C\|\F(z)\|_{L^{\infty}}\|z\|^2_{C^2}$.

Then we obtain $|K_3|+|K_4|\leq
C\|\F(z)\|_{L^{\infty}}\|z\|^3_{C^2}\|u\|_{L^2}$, and therefore
$J_1\leq C\|\F(z)\|_{L^{\infty}}\|z\|^3_{C^2}\|u\|_{L^2}$. Since the
estimate  $J_2\leq
C\|\F(z)\|^{1/2}_{L^{\infty}}\|z\|_{C^2}|H(u)(\al)|$ is immediate,
we finally have
\begin{align}\label{oti1}
\begin{split}
|I_1|\leq
C\|\F(z)\|_{L^{\infty}}\|z\|^3_{C^2}(\|u\|_{L^2}+|H(u)(\al)|).
\end{split}
\end{align}

Next we write $2BR(z,u)(\al)$ as follows:
\begin{align*}
\begin{split}
2BR(z,u)(\al)&=\frac{1}{2\pi}\int_{\T}u(\beta)(1-V_2^2(\al,\beta))\frac{V^{\bot}(\al,\beta)}{|V(\al,\beta)|^2}d\beta-
\frac{1}{2\pi}\int_{\T}u(\beta)V_2(\al,\beta)(1,0)d\beta\\
&= J_3(\al)+J_4(\al).
\end{split}
\end{align*}
Easily we have  $|\da J_4(\al)\cdot\da z(\al)|\leq
C\|z\|^2_{C^1}\|u\|_{L^2}$.Taking one derivative in $J_3(\al)$, and
using the cancellation $\da z(\al)\cdot\dpa z(\al)=0$, we get
$$\da J_3(\al)\cdot\da z(\al)=K_5+K_6+K_7+K_8+K_9,$$ where
$$
K_5=-\frac{1}{2\pi}\int_{\T}u(\beta)(1-V_2^2(\al,\beta))V_2(\al,\beta)\da
z_2(\al)\frac{V^{\bot}(\al,\beta)\cdot\da
z(\al)}{|V(\al,\beta)|^2}d\beta,
$$

$$
K_6=\frac{1}{4\pi}\int_{\T}u(\beta)(1-V_2^2(\al,\beta))\da
z_1(\al)\da z_2(\al) d\beta,
$$
$$
K_7=-\frac{1}{2\pi}\!\int_{\T}u(\beta)(1\!-\!V_2^2(\al,\beta))(\da
z_1(\al)V_1^3(\al,\beta)\!-\!\da
z_2(\al)V_2^3(\al,\beta)))\frac{V^{\bot}(\al,\beta)\cdot\da
z(\al)}{|V(\al,\beta)|^4}d\beta,
$$
$$
K_8=\frac{1}{2\pi}\!\int_{\T}\!u(\beta)V_2^2(\al,\beta)(\da
z_1(\al)V_1(\al,\beta)\!+\!\da
z_2(\al)V_2(\al,\beta))\frac{V^{\bot}(\al,\beta)\cdot\da
z(\al)}{|V(\al,\beta)|^4}d\beta,
$$
and
$$
K_9=-\frac{1}{2\pi}\!\int_{\T}\!u(\beta)(\da
z_1(\al)V_1(\al,\beta)\!+\!\da
z_2(\al)V_2(\al,\beta))\frac{V^{\bot}(\al,\beta)\cdot\da
z(\al)}{|V(\al,\beta)|^4}d\beta.
$$
We have $|K_5|+|K_6|+|K_7|+|K_8|\leq C\|z\|^2_{C^1}\|u\|_{L^2}$.

Next we split $K_9=-L_1-L_2$, where
$$
L_1=\frac{1}{2\pi}\!\int_{\T}\!u(\beta)\da
z_1(\al)V_1(\al,\beta)\frac{V^{\bot}(\al,\beta)\cdot\da
z(\al)}{|V(\al,\beta)|^4}d\beta,
$$
can be rewritten as follows
$$
L_1=\frac{1}{2\pi}\!\int_{\T}\!u(\al-\beta)\da
z_1(\al)V_1(\al,\al-\beta)\frac{(V(\al,\al-\beta)-\da
z(\al)\beta)^{\bot}\cdot\da z(\al)}{|V(\al,\al-\beta)|^4}d\beta.
$$
We have $L_1=M_1+M_2$ where
\begin{align*}
\begin{split}
M_1=\frac{(\da^2z(\al))^{\bot}\cdot\da z(\al)}{|\da z(\al)|^4}|\da z_1(\al)|^2 H(u)(\al),
\end{split}
\end{align*}
and
$$
M_2=\frac{1}{2\pi}\!\int_{\T}\!u(\al-\beta)\da
z_1(\al)B(\al,\al-\beta)d\beta,
$$
for
$$B(\al,\al-\beta)=V_1(\al,\al-\beta)\frac{V(\al,\al-\beta)^{\bot}\cdot\da z(\al)}{|V(\al,\al-\beta)|^4}-\da z_1(\al)\frac{(\da^2 z(\al))^{\bot}\cdot\da z(\al)}{|\da z(\al)|^4\tan(\beta/2)}.$$
The term  $M_1$ satisfies $|M_1|\leq
C\|\F(z)\|^{1/2}_{L^{\infty}}\|z\|_{C^2}|H(u)(\al)|$. We claim that
$$
|M_2|\leq C\|\F(z)\|^2_{L^{\infty}}\|z\|^4_{C^{2,\delta}}\int_\T
|\beta|^{\delta-1}|u(\al-\beta)|d\beta
$$ (see the appendix lemma \ref{B} for the proof).

A similar estimate can be obtained for $L_2$. Finally we have

$$
|I_2|\leq
C\|\F(z)\|^2_{L^{\infty}}\|z\|^4_{C^{2,\delta}}(\|u\|_{L^2}+|H(u)(\al)|+\int_\T
|\beta|^{\delta-1}|u(\al-\beta)|d\beta).
$$
This inequality together with \eqref{oti1} yields
$$
|\da T(u)(\al)|\leq
C\|\F(z)\|^2_{L^{\infty}}\|z\|^4_{C^{2,\delta}}(\|u\|_{L^2}+|H(u)(\al)|+\int_\T
|\beta|^{\delta-1}|u(\al-\beta)|d\beta).
$$
To finish we use the  $L^2$ boundedness of $H$ and  Minkowski's
inequality to obtain the estimate
$$
\|\da T(u)\|_{L^2}\leq
C\|\F(z)\|^2_{L^{\infty}}\|z\|^4_{C^{2,\delta}}\|u\|_{L^2}.
$$
q.e.d.

%%%%%%%%%%%%%%%%%%%%%%%%%%%%%%%%%%%%%%%%%%%%%%%%%%%%%%%%%%%%%%%%%%%%%%%%%%%%%%%%%%%%%%%%%%%%%%%%%%%%%%%%%%%%%%%%%%%%%%%%%%

%%%%%%%%%%%%%%%%%%%%%%%%%%%%%%%%%%%%%%%%%%%%%%%%%%%%%%%%%%%%%%%%%%%%%%%%%%%%%%%%%%%%%%%%%%%%%%%%%%%%%%%%%%%%%%%%%%%%%%%%%%

\section{Estimates on the inverse operator $(I-\xi T)^{-1}$}

In \ref{operatorT} we have considered the operator $T:L^2\rightarrow
H^1$
\begin{equation*}
T(u)=2BR(z,u)(\al)\cdot \da z(\al),
\end{equation*}
for $\F(z)(\al,\beta)<\infty$. Then $T$ is a compact operator
from Sobolev space $L^2$ to itself whose adjoint is given by the
formula
\begin{align*}
\begin{split}
T^*(u)(\al)&=\frac{1}{2\pi}\int_{\T}u(\beta)\frac{\da
z_2(\beta)\tan(\frac{z_1(\al)-z_1(\beta)}{2}) - \da
z_1(\beta)\tanh(\frac{z_2(\al)-z_2(\beta)}{2})}{\tan^2(\frac{z_1(\al)-
z_1(\beta)}{2})+\tanh^2(\frac{z_2(\al)-z_2(\beta)}{2})}d\beta\\
&\quad-\frac{1}{2\pi}\int_{\T}u(\beta)\frac{\da
z_2(\beta)\tan(\frac{z_1(\al)-z_1(\beta)}{2})\tanh^2(\frac{z_2(\al)-z_2(\beta)}{2})}
{\tan^2(\frac{z_1(\al)-
z_1(\beta)}{2})+\tanh^2(\frac{z_2(\al)-z_2(\beta)}{2})}d\beta\\
&\quad-\frac{1}{2\pi}\int_{\T}u(\beta)\frac{ \da
z_1(\beta)\tanh(\frac{z_2(\al)-z_2(\beta)}{2})\tan^2(\frac{z_1(\al)-z_1(\beta)}{2})}
{\tan^2(\frac{z_1(\al)-
z_1(\beta)}{2})+\tanh^2(\frac{z_2(\al)-z_2(\beta)}{2})}d\beta.
\end{split}
\end{align*}

We will show that, in $H^{\frac12}$, $I - \xi T$ has a bounded
inverse $(I - \xi T)^{-1}$ for $|\xi|\leq 1$, whose norm grows at
most like $\exp(C|||z|||^2)$ with
$|||z|||=\|z\|_{H^3}+\|\F(z)\|_{L^\infty}$.

Let $z$ be outside the curve $z(\al)$, then  we define
\begin{align*}
\begin{split}
f(z)&=\frac{1}{2\pi}\int_{\T}u(\beta)\frac{\da
z_2(\beta)\tan(\frac{z_1 -z_1(\beta)}{2}) - \da
z_1(\beta)\tanh(\frac{z_2 -z_2(\beta)}{2})}{\tan^2(\frac{z_1 -
z_1(\beta)}{2})+\tanh^2(\frac{z_2 -z_2(\beta)}{2})}d\beta\\
&\quad-\frac{1}{2\pi}\int_{\T}u(\beta)\frac{\da
z_2(\beta)\tan(\frac{z_1 -z_1(\beta)}{2})\tanh^2(\frac{z_2
-z_2(\beta)}{2})} {\tan^2(\frac{z_1 -
z_1(\beta)}{2})+\tanh^2(\frac{z_2 -z_2(\beta)}{2})}d\beta.\\
&\quad-\frac{1}{2\pi}\int_{\T}u(\beta)\frac{ \da
z_1(\beta)\tanh(\frac{z_2 -z_2(\beta)}{2})\tan^2(\frac{z_1
-z_1(\beta)}{2})} {\tan^2(\frac{z_1 -
z_1(\beta)}{2})+\tanh^2(\frac{z_2 -z_2(\beta)}{2})}d\beta.\\
&=\frac{1}{2\pi}\Im \int_{\T}\frac{u(\beta)\da
z(\beta)}{\tan(\frac{z-z(\beta)}{2})}d\beta
\end{split}
\end{align*}
that is $f$ is the real part of the Cauchy integral
\begin{align*}
\begin{split}
F(z)= f(z) + i g(z) =\frac{1}{2\pi i} \int_{\T}\frac{u(\beta)\da
z(\beta)}{\tan(\frac{z - z(\beta)}{2})}d\beta
\end{split}
\end{align*}
which is defined in both periodic domains $\Omega_1, \Omega_2$,
placed above and below, respectively, of the curve $z(\al)$. In
the following we shall denote by $\widetilde{\Omega}_j$ a
corresponding fundamental domain i.e. $\Omega_j =
\bigcup\{\widetilde{\Omega}_j + 2\pi n\}$.

Taking $z=z(\al)+\ep\dpa z(\al)$ we obtain
\begin{align*}
\begin{split}
f(z(\al)+\ep\dpa z(\al))=\frac{1}{2\pi}\Im
\int_{\T}\frac{u(\beta)\da
z(\beta)}{\tan(\frac{z(\al)-z(\beta)+\ep\dpa
z(\al)}{2})}d\beta
\end{split}
\end{align*}
and letting  $\ep\rightarrow 0$ we get
\begin{equation}\label{fec}
f(z(\al))=T^*(u) - \sign(\ep)u(\al).
\end{equation}
On the other hand we have
$$
\lim_{\ep\rightarrow 0} g(z(\al)+\ep\dpa z(\al)) = \lim_{\ep\rightarrow 0} \Im
F(z(\al)+\ep\dpa z(\al)) = \mathcal{G}(u)(\al)
$$
where
\begin{align*}
\begin{split}
\mathcal{G}(u)(\al)&=-\frac{1}{2\pi}PV\int_{\T}u(\beta)\frac{\da
z_1(\beta)\tan(\frac{z_1(\al)-z_1(\beta)}{2})(1-\tanh^2(\frac{z_2(\al)-z_2(\beta)}{2}))}{\tan^2(\frac{z_1(\al)-
z_1(\beta)}{2})+\tanh^2(\frac{z_2(\al)-z_2(\beta)}{2})}d\beta\\
&\quad-\frac{1}{2\pi}PV\int_{\T}u(\beta)\frac{\da
z_2(\beta)\tanh(\frac{z_2(\al)-z_2(\beta)}{2})(1+\tan^2(\frac{z_1(\al)-z_1(\beta)}{2}))}
{\tan^2(\frac{z_1(\al)-
z_1(\beta)}{2})+\tanh^2(\frac{z_2(\al)-z_2(\beta)}{2})}d\beta.
\end{split}
\end{align*}
independent of the sign of $\ep\rightarrow 0$.

First we will show that $T^*u = \lambda u$ $\Rightarrow$
$|\lambda|<1$, and since $T^*$ is a compact operator (of
Hilbert-Schmidt type) we can then conclude the existence of $(I -
\xi T^*)^{-1}$ with $|\xi|\leq 1$ (see also \cite{BMO}). To do that let us compute the
value of $\nabla f(z(\al))$. Denoting $z=x_1+ix_2$, the identity
\begin{align*}
\begin{split}
f(z)&=\frac{1}{2\pi}\Im \int_{\T}\frac{u(\beta)\da
z(\beta)}{\tan(\frac{z-z(\beta)}{2})}d\beta=-\frac{1}{2\pi} \Im
\int_{\T}u(\beta)\partial_{\beta}\ln(\sin(\frac{z-z(\beta)}{2}))
d\beta\\
&=\frac{1}{2\pi} \Im
\int_{\T}\partial_{\beta}u(\beta)\ln(\sin(\frac{z-z(\beta)}{2}))
d\beta,
\end{split}
\end{align*}
yields
\begin{align*}
\begin{split}
\nabla f(z)&=\frac{1}{2\pi} \Im
\int_{\T}\partial_{\beta}u(\beta)\nabla\ln(\sin(\frac{z-z(\beta)}{2}))
d\beta.
\end{split}
\end{align*}
That is
\begin{align*}
\begin{split}
\nabla f(x)=\big(&-\frac{1}{4\pi}
\int_{\T}\partial_{\beta}u(\beta)(\beta,t)\frac{\tanh(\frac{x_2-z_2(\beta,t)}{2})
(1+\tan^2(\frac{x_1-z_1(\beta,t)}{2}))}{\tan^2(\frac{x_1-
z_1(\beta,t)}{2})+\tanh^2(\frac{x_2-z_2(\beta,t)}{2})}d\beta,\\
&\frac{1}{4\pi}
\int_{\T}\partial_{\beta}u(\beta)\frac{\tan(\frac{x_1-z_1(\beta,t)}{2})
(1-\tanh^2(\frac{x_2-z_2(\beta,t)}{2}))}{\tan^2(\frac{x_1-
z_1(\beta,t)}{2})+\tanh^2(\frac{x_2-z_2(\beta,t)}{2})}d\beta\big).
\end{split}
\end{align*}
Taking the limit as before we get
\begin{equation}\label{gf}
\nabla f(z(\al))=2BR(z,\da u)(\al)+\sign(\ep)\frac{\da u(\al)}{2|\da
z(\al)|^2}\da z(\al)
\end{equation}
Assuming now that $T^* u=\lambda u$, the analyticity of $F(z)$
allows us to obtain:
\begin{align}
\begin{split}\label{autovalores1}
0 < &\int_{\widetilde{\Omega}_1}|F'(z)|^2 dx =-\int_{\T}
f(z(\al))\nabla f(z(\al))\cdot \frac{\dpa z(\al)}{|\da z(\al)|}
d\al\\
=&\int_{\T}(1-\lambda )u(\al) 2BR(z,\da u)(\al)\cdot \frac{\dpa
z(\al)}{|\da z(\al)|} d\al
\end{split}
\end{align}
and
\begin{align}
\begin{split}\label{autovalores2}
0 < &\int_{\widetilde{\Omega}_2}|F'(z)|^2 dx =\int_{\T}
f(z(\al))\nabla f(z(\al))\cdot \frac{\dpa z(\al)}{|\da z(\al)|}
d\al\\
=&\int_{\T}(\lambda +1)u(\al) 2BR(z,\da u)(\al)\cdot \frac{\dpa
z(\al)}{|\da z(\al)|} d\al
\end{split}
\end{align}
where we have used \eqref{fec} and \eqref{gf}. Multiplying
together both inequalities we get
\begin{align*}
\begin{split}
0 < &(1-\lambda^2)\Big( \int_{\T}u(\al) 2BR(z,\da u)(\al)\cdot
\frac{\dpa z(\al)}{|\da z(\al)|} d\al \Big)^2,
\end{split}
\end{align*}
and therefore $|\lambda| < 1$.

\begin{prop}\label{not}
The norms $||(I \mp T^*)^{-1}||_{L_0^2}$  are  bounded from above by
$\exp(C|||z|||^2)$ for some universal constant $C$ where the space
$L_0^2$ is the usual $L^2$ with the extra condition of mean value
zero i.e. the  subspace orthogonal to the constants.
\end{prop}

Proof: With the notation introduced before we have
\begin{align*}
\begin{split}
F_1 = F/\Omega_1 = f_1 + i g_1\\
F_2 = F/\Omega_2 = f_2 + i g_2\\
f_1/\partial \Omega = T^*u - u\\
f_2/\partial \Omega = T^*u + u\\
g_1/\partial \Omega =g_2/\partial \Omega = \mathcal{G}(u).
\end{split}
\end{align*}
The  proof follows easily from the estimate
\begin{eqnarray}
e^{-C|||z|||^2}\leq \frac{||u - T^*u||_{L_0^2}}{||u +
T^*u||_{L_0^2}}\leq  e^{C|||z|||^2}\label{comparacion}
\end{eqnarray}
valid for every nonzero $u\in L_0^2(\partial\Omega)$.

This is because if we assume  $||u - T^*u||_{L_0^2}\leq
e^{-2C|||z|||^2}$ for some $||u||_{L_0^2}=1$ then we obtain  $||u
+ T^*u||_{L_0^2}\geq 2||u||_{L_0^2} - e^{-2C|||z|||^2} \geq 1$
which contradicts \eqref{comparacion}. Therefore we must have $||u
- T^*u||_{L_0^2}\geq e^{-2C|||z|||^2}$ for all $||u||_{L_0^2}=1$
i.e. $||(I - T^*)^{-1}||_{L_0^2}\leq e^{2C|||z|||^2}$. Similarly
we also have $||(I + T^*)^{-1}||_{L_0^2}\leq e^{2C|||z|||^2}$.

Since $u - T^*u=\mathcal{H}_1(\mathcal{G}(u))$ and $u +
T^*u=\mathcal{H}_2(\mathcal{G}(u))$ where $\mathcal{H}_j$ denotes
the Hilbert transforms corresponding to each domain $\Omega_j$, then
\eqref{comparacion} is a consequence of the estimate
\begin{eqnarray}
||\mathcal{H}_j||_{L^2(\partial\Omega_j)}\leq
e^{C|||z|||^2}\label{estimacionH},
\end{eqnarray}
where $C$ denotes a universal constant not necessarily the same at
each occurrence.

This is because the identity $\mathcal{H}_j^2 = -I$ implies
\begin{align*}
||u - T^*u||_{L_0^2}&=||\mathcal{H}_1(\mathcal{G}(u))||_{L_0^2}\leq
e^{C|||z|||^2}||\mathcal{G}(u)||_{L_0^2}\\
&\leq e^{2C|||z|||^2}||\mathcal{H}_2(\mathcal{G}(u))||_{L_0^2}
= e^{2C|||z|||^2}||u + T^*u||_{L_0^2}
\end{align*}
and similarly we have $||u + T^*u||_{L_0^2}\leq e^{2C|||z|||^2}||u
- T^*u||_{L_0^2}$.

 It is enough to prove \eqref{estimacionH} for
$\Omega_1$ (the case $\Omega_2$ will follows  by symmetry) and  we
will rely on the following geometric fact whose elementary proof
is left to the reader.
\begin{lemma} Let $\Omega$ be a domain in $R^2$ whose
boundary is a $C^{2,\delta}$  parameterized curve $z(\alpha)$
satisfying the arc-chord condition
$||\mathcal{F}(z)||_{L^{\infty}}<\infty.$ Then we have tangent balls
to the boundary contained in both $\Omega$ and $R^2/\Omega$.
Furthermore, we can estimate from below the radius of those balls by
$C |||z|||^{-1}$, for some universal constant $C>0$.
\end{lemma}

Let $\phi = u + i v$ be the conformal mapping from $\Omega_1$ to the
upper half-plane $\mathbb{R}^2_+$. Then $v$ is a non-negative
harmonic function vanishing only on $\partial\Omega_1$. Let
$\phi^{-1}$ be the inverse transformation.

\begin{lemma}\label{geometricl}
Since $\Omega_1$ is $2\pi$ periodic in the horizontal direction we
have $\phi (z + 2\pi) = \phi (z) + \alpha$ for a certain fixed
real number $\alpha$.
\end{lemma}
Proof: Let us define $\psi(w) = \phi(\phi^{-1}(w) + 2\pi)$. Then
$\psi$ is a conformal mapping from $\mathbb{R}^2_+$ to itself and,
therefore, given by a linear fractional transformation $\psi(w)
=\frac{aw + b}{cw + d}$ satisfying $ad - bc = 1$, where  $a,b,c$ and
$d$ are real numbers. Since $\psi$ can not have a fixed point in
$\mathbb{R}^2_+$  then it follows that $c=0$ and $a=d$. Therefore
taking $z = \phi^{-1}(w)$ we get the formula $\phi (z + 2\pi) = \phi
(z) + \alpha$ with $\alpha = \frac{b}{d}$, proving lemma 4.3.\\

Next we observe that $\phi'(z + 2k\pi) = \phi'(z)$ for every
$z\in\Omega_1$ and  since $\partial\Omega_1$ is smooth enough we
know from general theory that $\phi$ and $\phi'$ extend
continuously to $\partial\Omega_1$. Furthermore, in order to
estimate the size of $\phi'|_{\partial\Omega_1}$ it will be enough
to consider the compact part of that boundary corresponding to a
full period.

Composing with $\phi$, $\phi^{-1}$ one easily gets the formula
\begin{eqnarray*}
\mathcal{H}_1f= H(f\circ\phi^{-1})\circ\phi\\
 H(f\circ\phi^{-1})=\mathcal{H}_1(f)\circ\phi^{-1}
\end{eqnarray*}
therefore our problem is reduced to a weighted estimate for the
Hilbert transform with respect to the weight
$|(\phi'(x,0))^{-1}|=w(x)$ for which we have to show that $w$
belongs to the Muckenhoupt class $A_2$ (see \cite{FKP}). Now it turns out that for
general $C^1$ chord-arc curves that statement is false, but we will
take advantage of the fact $\partial\Omega_1$ is of class $C^2$ (in
fact $C^{1,\alpha}$ will suffice) to show that in our case $w(x)$
trivializes i.e. it is bounded above and below, more precisely:
\begin{lemma}\label{geometric2}
Let $w(x)=|(\phi'(x,0))^{-1}|$ then we have
\begin{eqnarray*}
w(x_0) e^{-C|||z|||^2}\leq w(x)\leq w(x_0)
e^{C|||z|||^2}\label{mapping}
\end{eqnarray*}
where $C$ is a universal constant, $|||z|||$ is our usual norm in
the curve $\partial\Omega_1$ and $x_0$ is any point. Normalizing our
conformal mapping  $\phi$ one may take $w(x_0)=1$.
\end{lemma}

Proof: From the geometric lemma \ref{geometricl} we know the
existence of tangent balls to $\partial \Omega_1$ contained inside
$\overline{\Omega}_1$ of radius  $r=O(1/|||z|||)$ and such that
each one of those balls touches the boundary $\partial \Omega_1$
at a single point and its centers describe a parallel curve
$\Gamma$ to $\partial \Omega_1$ which is also of class $C^2$ with
norms $O(|||z|||)$. In the following we shall concentrate our
attention to the band $B$ of those points in $\Omega_1$ whose
distance to $\partial \Omega_1$ is less that $r$. Then the
boundary of $B$ consists of two parts, namely  $\partial \Omega_1$
and its parallel curve $\Gamma$ at distance $r$ which can also be
parameterized throughout $z(\al)$ in an obvious manner.

The length of the part of $\Gamma$ corresponding to a full period
$0\leq\alpha\leq 2\pi$, is clearly $O(|||z|||)$. Then, after
several applications of Harnack inequality in steps of length
$O(r)$, we obtain
\begin{eqnarray*}
e^{-C|||z|||^2}\leq \frac{v(z_1)}{v(z_2)}\leq e^{C|||z|||^2}
\end{eqnarray*}
 for any $z_1,z_2 \in \Gamma$. Let us consider a point $P\in \partial
 \Omega_1$ and  $Q\in \Gamma$ to be the center of the circle of radius
 $r$ tangent to $\partial \Omega_1$ at $P$, furthermore let us  denote
 by $\nu$ the inner normal vector. Then the non-negative harmonic
 function $v$ takes its strict minimum at the point $P$ and by
 Hopf principle we get the estimate
 \begin{equation}
\frac{\partial v}{\partial \nu}(P)\geq \frac{C}{r}
v(Q)\label{estimate1}
 \end{equation}
for some absolute constant $C>0$. On the other hand we may consider
a domain $D$ contained in the band $B$ in  such a way that its
boundary consists of a piece of $\partial \Omega_1$  of length $2r$
containing $P$ at its medium point, then the corresponding portion
of $\Gamma$, says $L_2$, obtained by vertical translation of the
points of $L_1$ and finally two arcs of $C^2$ curves smoothly
connecting $L_1$ and $L_2$ in such a way that $\partial D$ becomes a
$C^2$ curve with norm $O(|||z|||)$.

Let $\psi$ be conformal mapping from the unit ball $B_r$ to $D$ with
standard normalization. By the Kellogg-Warschawski theorem it
follows that $\psi$ extends continuously to the boundary and its
derivative is  bounded from above and below by  universal constants.
We also have the Poisson's kernel $K$ in $D$ obtained by conformal
mapping of the kernel for the ball of radius r. Then we may
represent the harmonic function $v$ as the integral of its boundary
values against the Poisson kernel:
 \begin{equation*}
  v(x) = \int_{\partial D}K(x,y)v(y)d\sigma(y)
 \end{equation*}
and
 \begin{equation*}
\frac{\partial v}{\partial \nu}(x)= \int_{\partial D}\frac{\partial
K}{\partial \nu_x}(x,y)v(y)d\sigma(y)
 \end{equation*}
which is a legitimate integral. We can take the limit (when
$x\rightarrow P\in\partial\Omega_1$) because $v$ vanishes
identically in $L_1$ and the points $y\in \partial D - L_1$ are at
distance at least $r$ from P to obtain the estimate
\begin{equation*}
\frac{\partial v}{\partial \nu}(P)\leq \frac{C}{r}sup_{x\in D}v(x)
\end{equation*}

To finish we can invoke  Dahlberg-Harnack principle up to the
boundary for the positive harmonic function $v$  (see \cite{CaC}
and \cite{Dahlberg}), which gives us the inequality

\begin{equation}
\frac{\partial v}{\partial \nu}(P)\leq
\frac{C}{r}v(Q)\label{estimate2}
\end{equation} for some fixed constant C. Then the estimates \eqref{estimate1}
and \eqref{estimate2} yield
\begin{equation}
C^{-1}\frac{v(Q_1)}{v(Q_2)}\leq \frac{\frac{\partial v}{\partial
\nu}(P_1)}{\frac{\partial v}{\partial \nu}(P_2)} \leq
C\frac{v(Q_1)}{v(Q_2)},
\end{equation}
but we know from Harnack that
\begin{eqnarray*}
e^{-C|||z|||^2}\leq \frac{v(Q_1)}{v(Q_2)}\leq  e^{C|||z|||^2}
\end{eqnarray*}
for two arbitrary points $Q_1, Q_2$ in $\Gamma$, and that ends the
proofs of lemma \ref{geometric2} and proposition \ref{not} because
$|\phi '(z(\al))| = |\nabla v(z(\al))| = \frac{\partial
v}{\partial\nu}(z(\al))$ since $\partial \Omega_1$ is the level
set $v=0$ of the positive harmonic function $v$ ($\phi= u + iv$)
q.e.d.

 The identity $ I + \xi T^* = \xi( I + T^*) + (1-\xi) I$ allows us
to conclude that
\begin{eqnarray*}
||(u + \xi T^*u)^{-1}||_{L_0^2}\leq  e^{C|||z|||^2}
\end{eqnarray*}
for $1-e^{-C_1|||z|||^2}\leq |\xi|\leq 1$ with an appropriate
constant $C_1$, but for general $\xi$ ($|\xi|\leq 1$) we have
\begin{prop}\label{not2}
For $|\xi|\leq 1$ the following estimate holds
\begin{eqnarray*}
||(I + \xi T)^{-1}||_{H_0^{\frac12}}= ||(I + \xi
T^*)^{-1}||_{H_0^{\frac12}}\leq e^{C|||z|||^2}
\end{eqnarray*}
for a universal constant C.
\end{prop}

Proof: First let us consider the inequality
\begin{eqnarray}
\int_{\Omega_j}|\nabla f_j|^2\geq e^{-C_2|||z|||^2}
||u||^2_{H^{\frac12}} \label{pag15}
\end{eqnarray}
where $F_j=f_j + i g_j$ is the Cauchy integral of $u$ in $\Omega_j$
which follows easily from  estimate (\ref{mapping}) for the
derivative of the conformal mapping $\phi$:

\begin{eqnarray*}
\int_{\Omega_j}|\nabla f_j|^2&=& \int_{\Omega_j}\Delta f_j^2=
\int_{\mathbb{R}^2_+} \Delta f^2_j(\phi^{-1}) |(\phi^{-1})'|^2\\
&=& \int_{\mathbb{R}^2_+} \Delta (f_j \circ \phi^{-1})^2 =
\int_{\partial\mathbb{R}^2_+} f_j \circ \phi^{-1}
\frac{\partial}{\partial \nu}f_j \circ \phi^{-1}
\end{eqnarray*}
where $\frac{\partial}{\partial \nu}$ is the derivative in the
normal direction
\begin{eqnarray*}
\frac{\partial f_j}{\partial y} (x,0) = \lim_{y\rightarrow
0}\frac{1}{\pi}\int\frac{u(x-t) - u(x)}{t^2 + y^2} dt = \Lambda
u(x).
\end{eqnarray*}
Therefore we can conclude that
\begin{eqnarray*}
\int_{\Omega_j}|\nabla f_j|^2&=&\int_{-\infty}^{+\infty} f_j \circ
\phi^{-1} \Lambda (f_j \circ \phi^{-1}) =\int_{-\infty}^{+\infty}
| \Lambda^{\frac12} (f_j \circ \phi^{-1})|^{2}\geq
e^{-C_2|||z|||^2} ||u||_{H^{\frac12}}^2
\end{eqnarray*}
 for a certain positive constant $C_2$ as a consequence of the
 following lemma:
\begin{lemma}
Let $\psi$ be a diffeomorphism of the real line such that
$0<C^{-1}\leq |\psi'(x)|\leq C$ then we have the equivalence of
Sobolev norms
\begin{eqnarray*}
C^{-(3+2s)}||f||_{H^s}\leq ||f\circ\psi||_{H^s}\leq
C^{3+2s}||f||_{H^s}
\end{eqnarray*}
for $0\leq s\leq \frac12$.
\end{lemma}
Proof: Given $f$ in $H^s$ we have
\begin{align*}
||\Lambda^s(f\circ\psi)||_{L^2}^2&= \int (f\circ\psi)(x)
\Lambda^{2s}(f\circ\psi)(x)dx= \int f(\psi(x)) \int\frac{f(\psi(x)) -
f(\psi(y))}{|x - y|^{1 + 2s}}dy dx \\
&=\frac12\int\int\frac{(f(\psi(x)) - f(\psi(y)))^2}{|x -
y|^{1 + 2s}}dy dx\\
& = \frac12\int\int\frac{(f(\bar{x}) -
f(\bar{y}))^2}{|(\psi^{-1})'(\bar{\bar{x}})|^{1 + 2s}|\bar{x} -
\bar{y}|^{1 + 2s}}(\psi^{-1})'(\bar{x})(\psi^{-1})'(\bar{y})d\bar{y}
d\bar{x}
\end{align*}
where $\bar{\bar{x}}$ comes from the application of the mean value
theorem. From our hypothesis we have
\begin{eqnarray*}
C^{-(3 + 2s)} \leq
\frac{(\psi^{-1})'(\bar{x})(\psi^{-1})'(\bar{y})}{|(\psi^{-1})'(\bar{\bar{x}})|^{1
+ 2s}}\leq C^{3 + 2s}
\end{eqnarray*}
which together with the equality
\begin{eqnarray*}
||\Lambda^s f||_{L^2}^2=\frac12\int\int\frac{(f(x) - f(y))^2}{|x -
y|^{1 + 2s}}dy dx
\end{eqnarray*}
allows us to finish the proof of the lemma.

\begin{rem}
In our case the diffeomorphism is given by
$\Psi(\alpha)=\Phi(z(\al))$ and we will use the periodic version
of (\ref{pag15}) i.e.
\begin{eqnarray*}
\int_{\widetilde{\Omega}_j}|\nabla f_j|^2\geq e^{-C_2|||z|||^2}
||u||^2_{H^{\frac12}(\mathbb{T})}.
\end{eqnarray*}
\end{rem}

To continue, let us assume that proposition \ref{not2} is false,
then there exist $u\in H^{\frac12}_0$, $||u||_{H^{\frac12}}=1$ and
$|\eta|>1$ such that $||\eta u - T^*u||_{H^{\frac12}}\leq e^{-
C_3|||z|||^2}$, where $C_3$ will be fixed later to be big enough for
our purposes.

 Let us also assume that the following estimate holds
\begin{align*}
| \int_{\T}&(\eta u -T^*u) 2BR(z,\da u)(\al)\cdot \frac{\dpa
z(\al)}{|\da z(\al)|} d\al|\leq \\
&\|\eta u -T^*u\|_{H^{\frac12}}\|2BR(z,\da u)(\al)\cdot \frac{\dpa
z(\al)}{|\da z(\al)|}\|_{H^{-\frac12}}\leq e^{-50 C_2|||z|||^2}
\end{align*}

 Then from identities (\ref{autovalores1}) and
(\ref{autovalores2}) we get
\begin{eqnarray}
e^{-C_2|||z|||^2}\leq  \int_{\T}(1-\eta )u(\al) 2BR(z,\da
u)(\al)\cdot \frac{\dpa z(\al)}{|\da z(\al)|} d\al   + e^{-50
C_2|||z|||^2}\label{estimateC2}
\end{eqnarray}
\begin{eqnarray*}
e^{-C_2|||z|||^2}\leq  \int_{\T}(1+\eta )u(\al) 2BR(z,\da
u)(\al)\cdot \frac{\dpa z(\al)}{|\da z(\al)|} d\al   + e^{-50
C_2|||z|||^2}.
\end{eqnarray*}
Adding these two inequalities together we obtain the positivity of
\begin{eqnarray*}
  \int_{\T}u(\al) 2BR(z,\da u)(\al)\cdot \frac{\dpa
z(\al)}{|\da z(\al)|} d\al> 0
\end{eqnarray*}
and then we get a contradiction when we substitute the value of
$\eta$ in the first inequality of (\ref{estimateC2}), if $\eta
\geq 1$, or in the second one, if $\eta \leq -1$. Therefore the
hypothesis $\|\eta u - T^*u\|_{H^\frac12}\leq e^{-C_3|||z|||^2}$ is
false for every $u$ in $H^{\frac12}_0$ and $\|u\|_{H^{\frac12}}=1$
and that gives us the desired estimate.

To finish the proof we need to show the following inequality
\begin{eqnarray*}
\|2BR(z,\da u)(\al)\cdot \frac{\dpa z(\al)}{|\da
z(\al)|}\|_{H^{-\frac12}}\leq e^{C|||z|||^2}\|u\|_{H^{\frac12}}
\end{eqnarray*}
for a universal constant $C$.

In order to prove it let us first observe that
\begin{eqnarray*}
\partial_\alpha(BR(z, u)(\al)\cdot \frac{\dpa
z(\al)}{|\da z(\al)|})= BR(z,\da u)(\al)\cdot \frac{\dpa
z(\al)}{|\da z(\al)|} + O_p(z)u
\end{eqnarray*}
where $O_p(z)$ is a bounded operator in $L^2$ whose norm is
controlled by $e^{C|||z|||^2}$ for a convenient value of $C$.
Therefore our task is equivalent to show the estimate
\begin{eqnarray*}
\|2BR(z, u)(\al)\cdot \frac{\dpa z(\al)}{|\da
z(\al)|}\|_{H^{\frac12}}\leq e^{C|||z|||^2}\|u\|_{H^{\frac12}}.
\end{eqnarray*}
Decoding the notation we have to consider the operators $\dpa
z_k(\al)\cdot T_j u(\al)$ where
\begin{eqnarray*}
T_j u(\al)= PV\int \frac{z_j(\al) - z_j(\al-\beta)}{|z(\al) -
z(\al-\beta)|^2}u(\al-\beta)d\beta.
\end{eqnarray*}
Let $\phi$ be a $C^{\infty}$ cut-off function supported on
$|x|\leq r$ such that $\phi\equiv 1$ on $|x|\leq \frac{r}{2}$
where $r=\frac{|||z|||}{2} $. Then $T_j u(\al)=T^1_j u + T^2_j u$ for
\begin{align*}
T^1_j u&=PV\int \phi(\beta)\frac{z_j(\al) -
z_j(\al-\beta)}{|z(\al) - z(\al-\beta)|^2}u(\al-\beta)d\beta,\\
T^2_j u&=PV\int (1 - \phi(\beta))\frac{z_j(\al) - z_j(\al-\beta)}{|z(\al)
-z(\al-\beta)|^2}u(\al-\beta)d\beta.
\end{align*}
It is  straightforward to check that $T^2_j$ is a smoothing
operator for which the desired estimate trivializes. Furthermore a
convenient Taylor expansion allows us to write $T_j^1 u(\al)=
m(\al)Hu(\al) + R(u)$ where $R$ is a smoothing operator, $H$ is
the Hilbert transform and the bounded smooth function $m$ depends
upon the curve $z$ in such a way that $\|\frac{\partial
m}{\partial\al}\|_{L^{\infty}}\leq e^{C|||z|||^2}$. Finally we may
invoke the following commutator estimate
\begin{eqnarray*}
\|\Lambda^{\frac12}(bv) -
b\Lambda^{\frac12}v\|_{L^2(\mathbb{T})}\leq C \|\nabla
b\|_{L^{\infty}}\|v\|_{L^2(\mathbb{T})}
\end{eqnarray*}
to complete our task. q.e.d.
$$
$$

\begin{rem}
Although it will not be needed to establish  our main theorem we
will improve the estimate on the eigenvalues of $T^*, T$ by showing
the existence of a constant $C_0=C_0(z)$ whose inverse $C_0^{-1}$
grows at most as a polynomial in $|||z|||$ and such that the
eigenvalues of $T^*$ must satisfy the estimate $|\lambda|\leq 1 -
C_0$. To see that let us consider the identities
\begin{align*}
\begin{split}
\int_{\Omega_1}|\nabla f_1|^2 dx + \int_{\Omega_2}|\nabla f_2|^2 dx
= 2\int_{\T}u(\al) 2BR(z,\da u)(\al)\cdot \frac{\dpa z(\al)}{|\da
z(\al)|} d\al\\
|\int_{\Omega_1}|\nabla f_1|^2 dx - \int_{\Omega_2}|\nabla f_2|^2
dx| = 2|\lambda|\int_{\T}u(\al) 2BR(z,\da u)(\al)\cdot \frac{\dpa
z(\al)}{|\da z(\al)|} d\al.
\end{split}
\end{align*}
Then it will be enough to show that both integrals
$\int_{\Omega_j}|\nabla f|^2 dx$ are  comparable i.e. there exists
a constant $1\leq C= C(|||z|||) < \infty$ such that
\begin{align*}
\begin{split}
\frac{1}{C}\int_{\Omega_1}|\nabla f_1|^2 \leq \int_{\Omega_2}|\nabla
f_2|^2 \leq C \int_{\Omega_1}|\nabla f_1|^2.
\end{split}
\end{align*}
Observe that the Cauchy-Riemman equations imply that this is
equivalent to show the analogous estimate for $g$ in place of $f$.

 The existence of such $C$ depending continuously upon $|||z|||$ follows easily
by a standard compactness argument. Nevertheless it is convenient
to have a control of the growth of the constants. In the following
we present an argument to show that $C(|||z|||)$  grows
polynomially with $|||z|||$.

\begin{prop} We shall
consider periodic curves $z(\alpha)$ (period $2\pi$). Because of the
smoothness and arc-chord conditions such a curve divide the cylinder
$\R/2\pi\Z\times(-\infty, \infty)$ in two regions $\Omega_j$
($j=1,2$,  above and below the curve respectively) containing
tangent balls as in the previous lemma. Then there exist a constant
$C= P(||z||_{C^{k,\delta}},||\mathcal{F}(z)||_{L^{\infty}})$,
polynomial in $|||z|||$,  such that
$$
\frac{1}{C}\int_{\Omega_1}|F'_1|^2 \leq \int_{\Omega_2}|F'_2|^2
\leq C \int_{\Omega_1}|F'_1|^2
$$ for any pair of periodic (in
$x_1$) holomorphic functions $F_j = f_j + ig_j$ ($j=1,2$), with
$f_j, g_j$ in Sobolev space $H^1(\Omega_j)$ and such that the
imaginary parts $g_j$, $j=1,2$ (or respectively the real part $f_j$
$j=1,2$) take the same boundary values.
\end{prop}

\underline{Proof:} In the following we shall use the expression
$P(\gamma)$ for different constants, to denote that they grow at
most polynomially with $\gamma$.

For $\frac{1}{r}= P(|||z|||)$ there exists two tangent circles to
the curve $z$ of radius $r$ and contained respectively in
$\Omega_1$ and $\Omega_2$. Therefore we can foliate the plane near
$z$ by parallel curves $z_{\epsilon}^{j}$ ($z_{0}^{j}=Z$), these
curves are the locus of points in $\Omega_j$ whose distance to
$z$ is $\epsilon$, in such a way that $|||z_{\epsilon}^{j}||| \leq
C |||z|||$ uniformly on $0\leq\epsilon\leq \frac{1}{10} r$ for
some universal finite constant $C$.

The Cauchy-Riemann equations for the holomorphic functions $F_j =
f_j + i g_j$ yields
$$
\int_{\Omega_j}|F'_j|^2 =  \int_{\Omega_j}|\nabla f_j|^2 =
\int_{\Omega_j}|\nabla g_j|^2
$$
Let us assume (without loss of generality)  that
$$
\int_{\Omega_1}|\nabla g_1|^2  \geq \int_{\Omega_2}|\nabla g_2|^2
$$
then we want to show the estimate
$$
\int_{\Omega_1}|\nabla g_1|^2  \leq P(|||z|||)
\int_{\Omega_2}|\nabla g_2|^2
$$
 and that will  finish the proof.

Let $\phi$ be a $C^{\infty}$ cut-off function such that
$\phi(t)\equiv 1$ when $|t|\leq\frac{1}{20} r$ and $\phi\equiv 0$
when $|t|\geq \frac{1}{10} r$, then we reflect the values of $g_1$
near $z(\al)$  by the formula
$$
\tilde{g}_1(P) = g_2(Q) \phi(dist(P,z))
$$
where $Q\in \Omega_2$ is obtained reflecting $P\in \Omega_1$ with
respect to $z$, that is $dist(P,z) = dist(Q,z)$, and the line
segment connecting them is normal to $z$ at its medium point.

By Dirichlet principle
$$
\int_{\Omega_1}|\nabla g_1|^2 \leq \int_{\Omega_1}|\nabla
\tilde{g}_1|^2\leq P(|||z|||)(\int_{\Omega_2}|\nabla g_2|^2|\phi|^2
+ \int_{\Omega_2}| g_2|^2|\nabla\phi|^2)
$$

Since $F_2$ is holomorphic we have the equalities
$$
\int_{0}^{2\pi} F_2(x,y_1) dx = \int_{0}^{2\pi} F_2(x,y_2) dx
$$
for $|y_j|$ big enough so that the horizontal lines $(x,y_j)$ do not
meet the curve $z$. The hypothesis that $f_j\in L^2(\Omega_j)$
implies that
$$
\int_{0}^{2\pi} g_2(x,y) dx =0
$$
for those $y$ which can be taken at distance $P(|||z|||)$ of the
curve $z$. For such a $y$ we get the estimate
$$
|g_2(x,y)|\leq\int_{0}^{2\pi} |\nabla g_2(t,y)| dt
$$
which implies
$$
\int_{0}^{2\pi}\int_{-y-1}^{-y}|g_2(x,s)|^2dsdx\leq
(2\pi)^2\int_{0}^{2\pi}\int_{-y-1}^{-y}|\nabla g_2(t,s)|^2dsdt
$$
therefore
$$
m\{|g_2(x,s)|\geq 10\cdot 2\pi||\nabla g_2||_{L^2(\Omega_2)} | 0\leq
x \leq 2\pi, -y-1\leq s \leq -y \}\leq \frac{1}{100}
$$
where $m$ denotes the Lebesgue measure.

Let $(x_m,y_m)$ be in the curve $z$ such that $y_m$ has a minimum
value. Then for all points $Q$ in $\Omega_2$ inside the band
$1/(20 r)\leq dist(Q,z)\leq 1/(20 r)$ whose distance to $(x_m,
y_m)$ is less than $1/P(|||z|||)$ (we shall denote by $N$ the set
of such $Q$) the segments connecting its points to those of
$\{(x,t), -y\leq t\leq -y-1\}$ are completely contained in
$\Omega_2$. For each $(x_0,y_0)\in N$ let us consider the line
segment connecting $(x_0,y_0)$ with the set
$$
E=\{ (x,s) | |g_2(x,s)| < 10\cdot 2\pi||\nabla
g_2||_{L^2(\Omega_2)} | 0\leq x \leq 2\pi, -y-1\leq s \leq -y \},
$$
then given $(x,s)\in E$ we have the estimate
$$
|g_2(x_0,y_0)| \leq 10\cdot 2\pi||\nabla g_2||_{L^2(\Omega_2)} +
\int_{0}^{L} |\nabla g_2((x_0,y_0) + t\omega)|dt
$$
where $\omega = \frac{(x-x_0, s-y_0)}{((x-x_0)^2 + (s-
y_0)^2)^{\frac12}}$ and $0\leq L\leq P(|||z|||)$.

Since the measure of E is big enough ($\geq \pi)$ the measure of
the region described in the unit circle by those $\omega$'s is
also big enough ($\geq 1/P(|||z|||)$). Therefore
$$
|g_2(x_0,y_0)| \leq P(|||z|||)\big (||\nabla g_2||_{L^2(\Omega_2)}
+ \int\int \frac{|\nabla g_2((x_0,y_0)
-(x,y))|}{||(x,y)||}dxdy\big ).
$$

This inequality implies that
$$
\int_{N}|g_2(x_0,y_0)|^2dx_0dy_0\leq P(|||z|||)
\int_{\Omega_2}|\nabla g_2(x,y)|^2dxdy
$$
To conclude the argument we just observe that because  the parallel
curves have tangent vector whose lengths are uniformly bounded by
$P(|||z|||)$, the integral $\int_{\Omega_2}| g_2|^2|\nabla\phi|^2$
is bounded by $P(|||z|||)(\int_{N}|g_2(x_0,y_0)|^2dx_0dy_0 +
\int_{\Omega_2}|\nabla g_2(x,y)|^2dxdy)$, q.e.d.
\end{rem}

%%%%%%%%%%%%%%%%%%%%%%%%%%%%%%%%%%%%%%%%%%%%%%%%%%%%%%%%%%%%%%%%%%%%%%%%%%%%%
%%%%%%%%%%%%%%%%%%%%%%%%%%%%%%%%%%%%%%%%%%%%%%%%%%%%%%%%%%%%%%%%%%%%%%%%%%%%%

\section{Estimates on $\varpi$}

%%%%%%%%%%%%%%%%%%%%%%%%%%%%%%%%%%%%%%%%%%%%%%%%%%%%%%%%%%%%%%%%%%%%%%%%%%%%%
%%%%%%%%%%%%%%%%%%%%%%%%%%%%%%%%%%%%%%%%%%%%%%%%%%%%%%%%%%%%%%%%%%%%%%%%%%%%%

In this section we show that the amplitude of the vorticity $\varpi$ is at the same level than $\da z$.
We prove the following result:

\begin{lemma}\label{lev}
Let $\varpi$ be a function given by
\begin{align}\label{formulavarpi}
\begin{split}
\varpi(\al)=-\frac{\mu^2-\mu^1}{\mu^2+\mu^1}2BR(z,\varpi)(\al)\cdot
\da z(\al)-2\kappa\mathrm{g}\frac{\rho^2-\rho^1}{\mu^2+\mu^1}\da
z_2(\al).
\end{split}
\end{align}
Then
\begin{equation}
\|\varpi\|_{H^k}\leq \exp C(\|\F(z)\|^2_{L^\infty}+\|z\|^2_{H^{k+1}}).
\end{equation}
for $k\geq 2$.
\end{lemma}

Proof:  We have $|A_\mu|\leq 1$, then the formula
\eqref{formulavarpi} is equivalent to
\begin{equation}\label{aux}
\varpi(\al)+A_\mu
T(\varpi)(\al)=-2\kappa\mathrm{g}\frac{\rho^2-\rho^1}{\mu^2+\mu^1}\da
z_2(\al),
\end{equation}
or
$$
\varpi(\al)=-2\kappa\mathrm{g}\frac{\rho^2-\rho^1}{\mu^2+\mu^1}(I+A_\mu
T)^{-1}(\da z_2)(\al).
$$
It yields
$$
\|\varpi\|_{H^{\frac12}}\leq C\|(I+A_\mu
T)^{-1}\|_{H^{\frac12}}\|\da z_2\|_{H^{\frac12}},
$$
and  proposition \eqref{not} gives
\begin{equation}\label{nl2varpi}
\|\varpi\|_{H^{\frac12}}\leq \exp
C(\|\F(z)\|^2_{L^\infty}+\|z\|^2_{H^3}).
\end{equation}

Taking the $k$ derivative of (\ref{aux}) we get:
\begin{equation*}
\da^k\varpi(\al)+A_\mu T(\da^k \varpi)(\al)= \Omega_k(\varpi) + C
\da^{k+1} z_2(\al), \quad
C=-2\kappa\mathrm{g}\frac{\rho^2-\rho^1}{\mu^2+\mu^1}
\end{equation*}
and using  Leibniz's rule we can write
\begin{align*}
\Omega_k(\varpi)(\al)&=\sum_{j=1}^{k} C_j\int
\Phi(\beta)\partial^j_\alpha\large (
\frac{(z(\alpha) - z(\alpha-\beta))^{\perp}\cdot
\da z(\alpha)}{|z(\alpha) - z(\alpha-\beta)|^2}\large ) \partial^{k-j}_\alpha
\varpi(\alpha-\beta) d\beta  + S(\varpi)(\al),
\end{align*}
where $S$ is a smoothing operator, $C_j$  are suitable constants and $\Phi$ is a $C^{\infty}$
cut-off such that $\Phi\equiv 0$ outside the ball $B(0,r)$ of
radius $r=\frac{1}{2|||z|||}$ and $\Phi\equiv 1$ in $B(0,
\frac{r}{2})$.

Next let us consider
\begin{align*}
\Lambda^{\frac12}\da^k\varpi(\al)+A_\mu T(\Lambda^{\frac12}\da^k
\varpi)(\al)&= A_\mu T(\Lambda^{\frac12}\da^k \varpi)(\al) - A_\mu
\Lambda^{\frac12}T(\da^k \varpi)(\al)\\
&\quad+\Lambda^{\frac12}\Omega_k(\varpi) + C \Lambda^{\frac12}\da^{k+1}z_2(\al),
\end{align*}
using the estimate for the inverse $(I+A_\mu T)^{-1}$ in the space $H^{\frac12}$ we get
\begin{align*}
\|\varpi\|_{H^{k+1}}&=\|\Lambda^{\frac12}\da^k\varpi\|_{H^{\frac12}}\\
&\leq  \exp
(C|||z|||^2) (A_\mu \|T\Lambda^{\frac12}\da^k\varpi\|_{H^{\frac12}}
+ A_\mu \|T\da^k\varpi\|_{H^{1}} + \|\Omega_k\|_{H^{1}} +
\|z\|_{H^{k+2}}).
\end{align*}
Then we have $$\|T\da^k\varpi\|_{H^{1}}\leq C|||z|||^4
\|\varpi\|_{H^{k}},$$ by Lemma 3.1, and
$$||\Omega_k||_{H^1}\leq C\|\F(z)\|^2_{L^{\infty}} ||\varpi||_{H^k}
||z||^2_{H^{k+2}},$$
by Lemma 5.2. (see below). Finally
\begin{align*}
\|T\Lambda^{\frac12}\da^k\varpi\|_{H^{\frac12}}&\leq \|T\Lambda^{\frac12}\da^k\varpi\|_{H^{1}}\leq C|||z|||^4
\|\varpi\|_{H^{k+\frac12}}\\
&\leq e^{C|||z|||^2}(\|\Omega_k\|_{H^{\frac12}} + \|z\|_{H^{k+\frac32}})\\
&\leq e^{C|||z|||^2}( ||\F(z)||_{L^{\infty}} ||\varpi||_{H^k}
||z||_{H^{k+2}} + ||z||_{H^{k+\frac32}})
\end{align*}
where we have used $\da^k\varpi (\alpha)=(I -
A_{\mu}T)^{-1}(\Omega_k + C\da^{k+1}z_2)$ and the estimate of the
norm in $H^{\frac12}$ of the inverse operator $(I-A_{\mu}T)^{-1}$.

A straightforward induction on $k\geq 2$  allows us to
finish the proof. The estimates for $k= 1, \frac32$ are obtained
similarly, but in all of them the norm $||z||_{H^3}$ has to appear
i.e. we have
\begin{equation*}
\|\varpi\|_{H^{\frac32}}\leq \exp
C(\|\F(z)\|^2_{L^\infty}+\|z\|^2_{H^3}).
\end{equation*}

\begin{lemma}
The operator $T$ maps Sobolev space $H^k$, $k\geq 1$, into
$H^{k+1}$ (so long as $\|z\|_{H^{k+2}} < \infty $) and satisfies
the estimate
\begin{eqnarray*}
\|T\|_{H^k\rightarrow H^{k+1}}\leq C |||z|||^2\|z\|^2_{H^{k+2}}
\end{eqnarray*}
\end{lemma}

Proof: We have
\begin{equation*}
T(u)(\al)=2BR(z,u)(\al)\cdot \da z(\al)\\
=\frac{1}{\pi} PV\int^{\infty}_{-\infty}\frac{(z(\alpha) -
z(\alpha-\beta))^{\perp} \da z(\alpha)}{|z(\alpha) -
z(\alpha-\beta)|^2}u(\alpha-\beta) d\beta
\end{equation*}
where, as usual and to simplify notation, we have dropped the time
dependence of all functions.

Let $\psi$ be a  $C^{\infty}$ cut-off such that $\psi\equiv 0$
outside the ball $B(0,r)$ of radius $r=\frac{1}{2|||z|||}$ and
$\psi\equiv 1$ in $B(0, \frac{r}{2})$. Then
\begin{align*}
T(u)(\al) &=\frac{1}{\pi}PV
\int^{\infty}_{-\infty}\psi(\beta)\frac{(z(\alpha) -
z(\alpha-\beta))^{\perp} \da z(\alpha)}{|z(\alpha) -
z(\alpha-\beta)|^2}u(\alpha-\beta) d\beta \\
&\quad +\frac{1}{\pi}
PV\int^{\infty}_{-\infty}(1 - \psi(\beta))\frac{(z(\alpha) -
z(\alpha-\beta))^{\perp} \da z(\alpha)}{|z(\alpha) -
z(\alpha-\beta)|^2}u(\alpha-\beta)d\beta\\
&= T_1 u(\alpha) + T_2 u(\alpha).
\end{align*}

i) Estimate of $T_2 u(\alpha)$: Leibniz's rule gives us
\begin{align*}
\da^{k+1}T_2 u(\alpha)&=\frac{1}{\pi} \int^{\infty}_{-\infty}(1 -
\psi(\beta))\frac{(z(\alpha) - z(\alpha-\beta))^{\perp}
\da z(\alpha)}{|z(\alpha) -
z(\alpha-\beta)|^2}\da^{k+1}u(\alpha-\beta)d\beta\\
&\quad+\frac{1}{\pi} \int^{\infty}_{-\infty}(1 -
\psi(\beta))\frac{(z(\alpha) - z(\alpha-\beta))^{\perp}
\da^{k+2}z(\alpha)}{|z(\alpha) -
z(\alpha-\beta)|^2}u(\alpha-\beta)d\beta\\
&\quad+ \text{\"{}other terms\"{}}\\
&= I_1 + I_2 + \text{\"{} other terms \"{} }
\end{align*}
The estimate for $\text{\"{}other terms\"{}}$ is straightforward.
For $I_1$ we integrate by parts:
\begin{align*}
I_1&=\frac{1}{\pi} \int^{\infty}_{-\infty}
\psi'(\beta)\frac{(z(\alpha) - z(\alpha-\beta))^{\perp}
\da z(\alpha)}{|z(\alpha) -
z(\alpha-\beta)|^2}\da^{k}u(\alpha-\beta)d\beta\\
&\quad -\frac{1}{\pi}\int^{\infty}_{-\infty}(1 -
\psi(\beta))\partial_{\beta}(\frac{(z(\alpha)-
z(\alpha-\beta))^{\perp} \da z(\alpha)}{|z(\alpha) -
z(\alpha-\beta)|^2})\da^{k}u(\alpha-\beta)d\beta.
\end{align*}
Then clearly we have:
\begin{eqnarray*}
||I_1||_{L^2}\leq C |||z|||^2||z||_{H^{k+2}}^2||u||_{H^k}
\end{eqnarray*}
Regarding $I_2$ we have
\begin{eqnarray*}
I_2(\alpha)=\sum^{2}_{j=1}\da^{k+2}z_j(\alpha)\cdot L_j u(\alpha)
\end{eqnarray*}
and clearly $ ||L_j u||_{L^{\infty}}\leq C |||z|||^2||u||_{H^k}$.
Therefore
\begin{eqnarray*}
||I_2||_{L^2}\leq C |||z|||^2||z||_{H^{k+2}}||u||_{H^k}
\end{eqnarray*}

ii)Estimate of $T_1u(\alpha)$: We have
\begin{eqnarray*}
\da^{k+1}T_1 u(\alpha)&=&\frac{1}{\pi}
\int^{\infty}_{-\infty}\psi(\beta)\frac{(z(\alpha) -
z(\alpha-\beta))^{\perp} \da z(\alpha)}{|z(\alpha) -
z(\alpha-\beta)|^2}\da^{k+1}u(\alpha-\beta)d\beta\\ &+&
\frac{1}{\pi} \int^{\infty}_{-\infty}\psi(\beta)\frac{(z(\alpha) -
z(\alpha-\beta))^{\perp} \da^{k+2}z(\alpha)}{|z(\alpha) -
z(\alpha-\beta)|^2}u(\alpha-\beta)d\beta\\&+& \frac{1}{\pi}
\int^{\infty}_{-\infty}\psi(\beta)\frac{(\da^{k+1}z(\alpha) -
\da^{k+1}z(\alpha-\beta))^{\perp} \da z(\alpha)}{|z(\alpha) -
z(\alpha-\beta)|^2}u(\alpha-\beta)d\beta\\ &+& \text{\"{}other terms\"{}}\\
&=& J_1 + J_2 + J_3+ \text{\"{}other terms\"{}}.
\end{eqnarray*}
As in the previous case the $\text{\"{}other terms\"{}}$ are easy
to handle and we shall show how to estimate the remainder two
cases.

We can write $J_2$ in the form
\begin{eqnarray*}
J_2(\alpha) = \sum\da^{k+2}z_j(\alpha) \int\psi(\beta)K_j(\alpha,
\alpha - \beta)u(\alpha-\beta) d\beta =
\sum^{2}_{j=1}\da^{k+2}z_j(\alpha)\cdot L_j u(\alpha)
\end{eqnarray*}
and observe that
\begin{eqnarray*}
||L_ju||_{L^{\infty}}\leq ||L_ju||_{H^1}\leq C
|||z|||^2||z||_{H^2}||u||_{H^1}
\end{eqnarray*}
which yields
\begin{eqnarray*}
||J_2||_{L^2}\leq C |||z|||^2||z||_{H^{k+2}}^2||u||_{H^k}.
\end{eqnarray*}

To estimate $J_1$ we integrate by parts
\begin{eqnarray*}
J_1(\alpha) &=& \frac{1}{\pi}
\int^{\infty}_{-\infty}\psi'(\beta)\frac{(z(\alpha) -
z(\alpha-\beta))^{\perp} \da z(\alpha)}{|z(\alpha) -
z(\alpha-\beta)|^2} \da^k u(\alpha-\beta) d\beta\\ &+& \frac{1}{\pi}
\int^{\infty}_{-\infty}\psi(\beta)\partial_{\beta}(\frac{(z(\alpha)
- z(\alpha-\beta))^{\perp} \da z(\alpha)}{|z(\alpha) -
z(\alpha-\beta)|^2}) \da^k u(\alpha-\beta) d\beta\\ &=& J_1^1 +
J_1^2.
\end{eqnarray*}
For the first part $J_1^1$ we have
\begin{eqnarray*}
||J_1^1||_{L^2}\leq C |||z|||^2||z||_{H^{3}}^2||u||_{H^k}.
\end{eqnarray*}
And
\begin{align*}
J_1^2(\alpha)&=\frac{1}{\pi}
\int^{\infty}_{-\infty}\psi(\beta)\frac{(\da
z(\alpha-\beta))^{\perp} \da z(\alpha)}{|z(\alpha) -
z(\alpha-\beta)|^2} \da^k u(\alpha-\beta) d\beta\\
&\quad+
\frac{2}{\pi}
\int^{\infty}_{-\infty}\psi(\beta)\frac{[(z(\alpha)\!
-\! z(\alpha\!-\!\beta))^{\perp} \da z(\alpha)][
\da z(\alpha\!-\!\beta)(z(\alpha) \!-\! z(\alpha\!-\!\beta))]}{|z(\alpha)
\!-\! z(\alpha\!-\!\beta)|^4} \da^k u(\alpha\!-\!\beta) d\beta\\
&=J_1^{2,1} + J_1^{2,2}.
\end{align*}
For $J_1^{2,1}$ $$(\da z(\alpha-\beta))^{\perp} \da z(\alpha)=(\da
z(\alpha-\beta)-\da z(\al) )^{\perp} \da z(\alpha)=-\da^2z^{\perp}(\alpha)\da z(\alpha)\beta+O(\beta^2),$$
and
\begin{eqnarray*}
|z(\alpha) - z(\alpha-\beta)|^2 = |\da z|^2\beta^2 +
O(\beta^3)
\end{eqnarray*}
where the constants in the "O" terms (and in theirs first
derivatives) are properly bounded in terms of $||z||_{H^3}$.

That is
\begin{eqnarray*}
J_1^{2,1}(\alpha) = -\frac{ \da^2 z^{\perp}(\alpha)\da
z(\alpha)}{|z_{\alpha}(\alpha)|^2}H\da^k u(\alpha) + \text{\"{}
bounded terms \"{}}
\end{eqnarray*}
where $H$ denotes the Hilbert transform. Therefore for the first integral we get
\begin{eqnarray*}
||J_1^{2,1}||_{L^2}\leq C |||z|||^2||z||_{H^{3}}^2||u||_{H^k}.
\end{eqnarray*}
Finally for $J_1^{2,2}$ we have
\begin{eqnarray*}
[(z(\alpha)\!
-\! z(\alpha\!-\!\beta))^{\perp} \da z(\alpha)][
\da z(\alpha\!-\!\beta)(z(\alpha) \!-\! z(\alpha\!-\!\beta))]=\frac12
\da^2 z^{\perp}(\alpha)\da z(\alpha)|\da z(\alpha)|^2\beta^3\!
+\! O(\beta^4)
\end{eqnarray*}
and
\begin{eqnarray*}
|z(\alpha) - z(\alpha-\beta)|^4 = |\da z|^4\beta^4 +
O(\beta^4).
\end{eqnarray*}
By a similar approach we obtain
\begin{eqnarray*}
J_1^{2,1}(\alpha) = \frac{\da^2 z^{\perp}(\alpha)\da z(\alpha)}{|\da z(\alpha)|^2}H\da^k
u(\alpha) + \text{\"{} bounded terms \"{}},
\end{eqnarray*}
and it yields
\begin{eqnarray*}
||J_1^{2,2}||_{L^2}\leq C |||z|||^2||z||_{H^{3}}^2||u||_{H^k}.
\end{eqnarray*}

To estimate $J_3$ we observe first that the substitution of
$u(\alpha -\beta)$ by $u(\alpha) - \da u(\alpha)\beta$ produces an
error bounded by $||z||^2_{H^{k+2}}|||z|||^2 ||u||_{H^k}$.

Using the expansions
\begin{eqnarray*}
\frac{\psi(\beta)}{|z(\alpha) - z(\alpha-\beta)|^2}=
\frac{\psi(\beta)}{|\da z(\alpha)|^2}\frac{1}{|\beta|^2} +
O(1)\psi(\beta)\\ \da^{k+1}z(\alpha) - \da^{k+1}z(\alpha-\beta) =
\beta\int_0^1\da^{k+2}z(\alpha-t\beta)dt
\end{eqnarray*}
and since the term corresponding to $\da u(\alpha)\beta$ can be
handled very easily, it remains to estimate:
\begin{eqnarray*}
\frac{u(\alpha)}{\pi}
\int^{\infty}_{-\infty}\psi(\beta)\frac{(\da^{k+1}z(\alpha) -
\da^{k+1}z(\alpha-\beta))^{\perp} \da z(\alpha)}{|z(\alpha) -
z(\alpha-\beta)|^2}d\beta = \int_0^1 K_t(\alpha)dt,
\end{eqnarray*}
where
\begin{eqnarray*}
K_t(\alpha) =\frac{u(\alpha)}{\pi|\da z(\alpha)|^2}
\int^{\infty}_{-\infty}\psi(\frac{\beta}{t})\frac{\da^{k+2}z^{\perp}(\alpha-\beta)
\cdot\da z(\alpha)}{\beta}d\beta.
\end{eqnarray*}
Finally, the $L^2$ boundedness of the Hilbert transform yields
\begin{eqnarray*}
||K_t||_{L^2} \leq ||z||^2_{H^{k+2}} ||u||_{H^k}|||z|||
\end{eqnarray*}
uniformly on $t$, allowing us to finish the proof.
%%%%%%%%%%%%%%%%%%%%%%%%%%%%%%%%%%%%%%%%%%%%%%%%%%%%%%%%%%%%%%%%%%%%%%%%%%%%%
%%%%%%%%%%%%%%%%%%%%%%%%%%%%%%%%%%%%%%%%%%%%%%%%%%%%%%%%%%%%%%%%%%%%%%%%%%%%%

\section{Estimates on $BR(z,\varpi)$}

%%%%%%%%%%%%%%%%%%%%%%%%%%%%%%%%%%%%%%%%%%%%%%%%%%%%%%%%%%%%%%%%%%%%%%%%%%%%%
%%%%%%%%%%%%%%%%%%%%%%%%%%%%%%%%%%%%%%%%%%%%%%%%%%%%%%%%%%%%%%%%%%%%%%%%%%%%%

This section is devoted to show that the Birkhoff-Rott integral is as regular as $\da z$.

\begin{lemma}The following estimate holds
\begin{eqnarray}\label{nsibr}
\|BR(z,\varpi)\|_{H^k}\leq \exp(C(\|\F(z)\|^2_{L^\infty}+\|z\|^2_{H^{k+1}}),
\end{eqnarray}
for $k\geq 2$.
\end{lemma}

\begin{rem}
Using this estimate for $k=2$ we find easily that
\begin{eqnarray}\label{nliibr}
\|\da BR(z,\varpi)\|_{L^\infty}\leq \exp(C(\|\F(z)\|^2_{L^\infty}+\|z\|^2_{H^{3}}),
\end{eqnarray}
which shall be used through out the paper.
\end{rem}

Proof: We show the proof for $k=2$, being the rest of the cases analogous. We have
\begin{align*}
BR(z,\varpi)(\al)&=\frac{1}{4\pi}
PV\!\!\!\int_{\T}\!\!\varpi(\al\!-\!\beta)\big(-\frac{V_2(\al,\beta)
(1\!+\!V^2_1(\al,\beta))}{|V(\al,\beta)|^2},\frac{V_1(\al,\beta)
(1\!-\!V^2_2(\al,\beta))}{|V(\al,\beta)|^2}\big)d\beta
\end{align*}
which is decomposed as follows:

\begin{align}
\begin{split}\label{p1p2p3}
BR(z,\varpi)(\al)&=\frac{1}{4\pi}PV\int_{\T}\varpi(\al-\beta)\frac{V^{\bot}(\al,\al-\beta)}{|V(\al,\al-\beta)|^2}d\beta\\
&\quad-\frac{1}{4\pi}\int_{\T}\varpi(\beta)V_2^2(\al,\beta)\frac{V^{\bot}(\al,\beta)}{|V(\al,\beta)|^2}d\beta-\frac{1}{4\pi}\int_{\T}\varpi(\beta)V_2(\al,\beta)d\beta (1,0)\\
&=P_1(\al)+P_2(\al)+P_{3}(\al).
\end{split}
\end{align}
Using that $|V_2(\al,\beta)|\leq 1$, we get $|P_2(\al)|+|P_3(\al)|\leq
C\|\varpi\|_{L^2}$, and lemma \ref{lev} yields $\|P_2\|_{L^2}+\|P_3\|_{L^2}\leq
\exp(C(\|\F(z)\|^2_{L^\infty}+\|z\|^2_{H^{3}})$.

Let us write
\begin{align*}
\begin{split}
P_1(\al)&=\frac{1}{4\pi}\int_{\T}(-A_1(\al,\al-\beta),A_2(\al,\al-\beta))\varpi(\al-\beta)d\beta
+\frac{\dpa
z(\al)}{|\da z(\al)|^2}H \varpi(\al)d\al\\
&=J_1+J_2,
\end{split}
\end{align*}
where as before
$$
A_1(\al,\al-\beta)=\frac{V_2(\al,\al-\beta)}{|V(\al,\al-\beta)|^2}
-\frac{1}{|\da z(\al)|^2}\frac{\da
z_2(\al)}{\tan(\frac{\beta}{2})},
$$
and
$$
A_2(\al,\al-\beta)=\frac{V_1(\al,\al-\beta)}{|V(\al,\al-\beta)|^2}
-\frac{1}{|\da z(\al)|^2}\frac{\da
z_1(\al)}{\tan(\frac{\beta}{2})}.
$$
For $J_1$ since $\|A_1\|_{L^{\infty}}\leq \|\F(z)\|_{L^{\infty}}\|z\|_{C^2}$ and $\|A_2\|_{L^\infty}\leq
C\|\F(z)\|_{L^{\infty}}\|z\|^2_{C^2}$ (see appendix) one gets
$\|J_1\|_{L^2}\leq C\|\F(z)\|_{L^{\infty}}\|z\|_{C^2}\|z\|_{L^2}\|\varpi\|_{L^2}$. The inequality $|\da z(\al)|^{-1}\leq \|\F(z)\|^{1/2}_{L^\infty}$ give us $\|J_2\|_{L^2}\leq C\|\F(z)\|^{1/2}_{L^\infty} \|z\|_{L^2}\|\varpi\|_{L^2}$.

 Next it is easy to check that $|\da^2 P_3(\al)|\leq C\|\varpi\|_{L^2}(|\da^2z(\al)|+\|z\|^2_{C^2})$ and to estimate
$\|\da^2 P_3\|_{L^2}$. The kernel in the integral $P_2(\al)$ has order $1$ in $\beta$, and taking two derivatives in $\al$ we get integrals as in $P_3$ and kernels of degree $-1$ which can be estimated as before. Similar terms of lower order are obtained in $\da^2 P_1(\al)$ which are controlled analogously. The most
singular terms are given by
$$
Q_1(\al)=\frac{1}{4\pi}PV\int_{\T}\da^2\varpi(\al-\beta)\frac{V^\bot(\al,\al-\beta)}{|V(\al,\al-\beta)|^2}d\beta,
$$
$$
Q_2(\al)=\frac{1}{8\pi}PV\int_{\T}\varpi(\al-\beta)\frac{\da^2z(\al)-\da^2z(\al-\beta)}{|V(\al,\al-\beta)|^2}d\beta,
$$
$$
Q_3(\al)=-\frac{1}{4\pi}PV\int_{\T}\varpi(\al-\beta)\frac{V^\bot(\al,\al-\beta)}{|V(\al,\al-\beta)|^4}\big(V(\al,\al-\beta)\cdot(\da^2 z(\al)-\da^2 z(\al-\beta))\big)d\beta.
$$
We have
$$
Q_1=\frac{1}{4\pi}PV\!\!\int_{\T}\!\da^2\varpi(\al\!-\!\beta)\big(\frac{V^\bot(\al,\al\!-\!\beta)}{|V(\al,\al\!-\!\beta)|^2}-\frac{\da^{\bot} z(\al)}{|\da z(\al)|^2\tan(\beta/2)}\big)d\beta+\frac{\da^{\bot} z(\al)}{2|\da z(\al)|^2}H(\da^2\varpi)(\al),
$$
giving us
\begin{align}\label{q1nl22dbr}
\begin{split}
|Q_1(\al)|&\leq C\|\F(z)\|^k_{L^\infty}\|z\|^k_{C^2}\|\da^2\varpi\|_{L^2}+\|\F(z)\|^{1/2}_{L^\infty}|H(\da^2\varpi)(\alpha)|\\
&\leq (1+|H(\da^2\varpi)(\alpha)|)\exp C(\|\F(z)\|^2_{L^\infty}(t)+\|z\|^2_{H^3}(t)).
\end{split}
\end{align}
Next we write $Q_2=R_1+R_2+R_3$ where
$$
R_1(\al)=\frac{1}{8\pi}\int_{\T}(\varpi(\al-\beta)-\varpi(\al))\frac{\da^2 z(\al)-\da^2 z(\al-\beta)}{|V(\al,\al-\beta)|^2}d\beta,
$$
$$
R_2(\al)=\frac{\varpi(\al)}{8\pi}\int_{\T}(\da^2z(\al)-\da^2z(\al-\beta))\big(\frac{1}{|V(\al,\al-\beta)|^2}-\frac{4}{|\da z(\al)|^2|\beta|^2}\big)d\beta,
$$
$$
R_3(\al)=\frac{1}{8\pi}\frac{\varpi(\al)}{|\da z(\al)|^2}\int_{\T}(\da^2z(\al)\!-\!\da^2z(\al\!-\!\beta))\big(\frac{4}{|\beta|^2}\!-\!\frac{1}{\sin^2(\beta/2)}\big)d\beta\!+\!\frac{1}{2}\frac{\varpi(\al)}{|\da z(\al)|^2}\la (\da^2z)(\al).
$$
Using that $$|\da^2z(\al)-\da^2z(\al-\beta)|\leq|\beta|^\delta\|z\|_{C^{2,\delta}},$$
we get
$$|R_1(\al)|+|R_2(\al)|\leq \|\varpi\|_{C^1}\|\F(z)\|^{k}\|z\|^k_{C^{2,\delta}}\leq \exp C(\|\F(z)\|^2_{L^\infty}+\|z\|^2_{H^3}).$$
While for $R_3$ we have
\begin{align*}
|R_3(\al)|&\leq C\|\varpi\|_{L^\infty}\|\F(z)\|_{L^\infty}(\|z\|_{C^2}+|\la (\da^2z)(\al)|)\\
&\leq (1+|\la (\da^2z)(\al)|)\exp C(\|\F(z)\|^2_{L^\infty}+\|z\|^2_{H^3}),
\end{align*}
that is
\begin{equation}\label{q2nl22dbr}
|Q_2(\al)|\leq (1+|\la (\da^2z)(\al)|)\exp C(\|\F(z)\|^2_{L^\infty}+\|z\|^2_{H^3}).
\end{equation}
Let us consider  $Q_3=R_4+R_5+R_6+R_7+R_8+R_9$, where
\begin{align*}
R_4=-\frac{1}{4\pi}\int_{\T}(\varpi(\al-\beta)-\varpi(\al))\frac{V^\bot(\al,\al-\beta)}{|V(\al,\al-\beta)|^4}\big(V(\al,\al-\beta)\cdot(\da^2 z(\al)-\da^2 z(\al-\beta))\big)d\beta,
\end{align*}
$$
R_5=-\frac{\varpi(\al)}{4\pi}\int_{\T}\frac{(V(\al,\al-\beta)-\da z(\al)\beta/2)^{\bot}}{|V(\al,\al-\beta)|^4}\big(V(\al,\al-\beta)\cdot(\da^2 z(\al)-\da^2 z(\al-\beta))\big)d\beta,
$$
$$
R_6=-\frac{\varpi(\al)(\da z(\al))^\bot}{8\pi}\int_{\T}\frac{\beta(V(\al,\al-\beta)-\da z(\al)\beta/2)\cdot(\da^2 z(\al)-\da^2 z(\al-\beta))}{|V(\al,\al-\beta)|^4}d\beta,
$$
$$
R_7=-\frac{\varpi(\al)(\da z(\al))^\bot}{16\pi}\da z(\al)\cdot\!\!\int_{\T}\beta^2(\da^2 z(\al)-\da^2 z(\al-\beta))\big(\frac{1}{|V(\al,\al\!-\!\beta)|^4}-\frac{16}{|\da z(\al)|^4|\beta|^4}\big)d\beta,
$$
$$
R_8=-\frac{\varpi(\al)(\da z(\al))^\bot}{4\pi|\da z(\al)|^4}\da z(\al)\cdot\!\!\int_{\T}(\da^2 z(\al)-\da^2 z(\al-\beta))\big(\frac{4}{|\beta|^2}-\frac{1}{\sin^2(\beta/2)}\big)d\beta,
$$
and
$$
R_9=-\frac{\varpi(\al)(\da z(\al))^\bot}{|\da z(\al)|^4}\da z(\al)\cdot \la(\da^2 z(\al)).
$$
Proceeding as before we get
\begin{equation*}
|Q_3(\al)|\leq (1+|\la (\da^2z)(\al)|)\exp
C(\|\F(z)\|^2_{L^\infty}+\|z\|^2_{H^3}),
\end{equation*}
which together with \eqref{q1nl22dbr} and \eqref{q2nl22dbr} gives
us the estimate
$$
|\da^2 P_1(\al)|\leq (1+|\la (\da^2z)(\al)|+|H(\da^2\varpi)(\alpha)|)\exp C(\|\F(z)\|^2_{L^\infty}+\|z\|^2_{H^3}),
$$
and  $\|\da^2 P_1\|_{L^2}\leq \exp
C(\|\F(z)\|^2_{L^\infty}+\|z\|^2_{H^3}).$

Finally we get
\begin{equation}\label{enl22dbr}
\|\da^2 BR(z,\varpi)\|_{L^2}\leq \exp C(\|\F(z)\|^2_{L^\infty}+\|z\|^2_{H^3}).
\end{equation}

%%%%%%%%%%%%%%%%%%%%%%%%%%%%%%%%%%%%%%%%%%%%%%%%%%%%%%%%%%%%%%%%%%%%%%%%%%%%%%%%%%%%%%%%%%%%%%%%%%%%%%%%%%%%%%%%%%%%%%%

%%%%%%%%%%%%%%%%%%%%%%%%%%%%%%%%%%%%%%%%%%%%%%%%%%%%%%%%%%%%%%%%%%%%%%%%%%%%%%%%%%%%%%%%%%%%%%%%%%%%%%%%%%%%%%%%%%%%%%%

\section{Estimates on $z(\al,t)$}

In this section we give the proof of the below lemma for $k=3$. The case $k>3$ is left to the reader.

\begin{lemma}
Let $z(\al,t)$ be a solution of 2DM. Then, the following a priori
estimate holds:
\begin{align}
\begin{split}\label{ntni}
\frac{d}{dt}\|z\|^2_{H^k}(t)&\leq -\frac{\kappa}{2\pi(\mu_1\!+\!\mu_2)}\,\int_\T \frac{\sigma(\al,t)}{|\da z(\al)|^2} \da^k z(\al,t)\cdot \la(\da^k z)(\al,t) d\al\\
&\quad+\exp C(\|\F(z)\|^2_{L^\infty}(t)+\|z\|^2_{H^k}),
\end{split}
\end{align}
for $k\geq 3$.
\end{lemma}

We split the proof in the following four parts.

%%%%%%%%%%%%%%%%%%%%%%%%%%%%%%%%%%%%%%%%%%%%%%%%%%%%%%%%%%%%%%%%%%%%%%%%%%%%%%
%%%%%%%%%%%%%%%%%%%%%%%%%%%%%%%%%%%%%%%%%%%%%%%%%%%%%%%%%%%%%%%%%%%%%%%%%%%%%%

\subsection{Estimates for the $L^2$ norm of the curve}

%%%%%%%%%%%%%%%%%%%%%%%%%%%%%%%%%%%%%%%%%%%%%%%%%%%%%%%%%%%%%%%%%%%%%%%%%%%%%%
%%%%%%%%%%%%%%%%%%%%%%%%%%%%%%%%%%%%%%%%%%%%%%%%%%%%%%%%%%%%%%%%%%%%%%%%%%%%%%

We have

\begin{align*}
\begin{split}
\frac12\frac{d}{dt}\int_{\T}|z(\alpha)|^2
d\alpha=\int_{\T}z(\al)\cdot z_t(\al)d\al
&=\int_{\T}z(\al)\cdot BR(z,\varpi)(\al) d\al+
\int_{\T} c(\al) z(\al)\cdot\da z(\al) d\al\\
&=I_1+I_2.
\end{split}
\end{align*}
Taking $I_1\leq \|z\|_{L^2}\|BR(z,\varpi)\|_{L^2}$ and the inequality \eqref{nsibr} let us estimate $I_1$.

Next we get

$$I_2\leq  A^{1/2}(t) \|c\|_{L^\infty} \int_\T |z(\al)|d\al\leq
2\int_{\T}|\partial_{\al} BR(z,\varpi)(\al)|d\al\int_\T
|z(\al)|d\al
$$
which yields
$$
I_2\leq \exp (C\|F(z)\|^2_{L^{\infty}} + \|z\|^2_{H^3})
$$
if we consider the estimate \eqref{nliibr}.  We conclude that
\begin{equation}\label{nl2}
\D\dt\|z\|^2_{L^2}(t)\leq \exp(C|||z|||^2)
\end{equation}
for an appropriate finite constant $C$, where as before $|||z|||^2
=\|F(z)\|^2_{L^{\infty}} + \|z\|^2_{H^3}$.

%%%%%%%%%%%%%%%%%%%%%%%%%%%%%%%%%%%%%%%%%%%%%%%%%%%%%%%%%%%%%%%%%%%%%%%%%%%%%%%%%%%%%%%%%%%%%%%%%%%%%%%%%
%%%%%%%%%%%%%%%%%%%%%%%%%%%%%%%%%%%%%%%%%%%%%%%%%%%%%%%%%%%%%%%%%%%%%%%%%%%%%%%%%%%%%%%%%%%%%%%%%%%%%%%%%

\subsection{The integrable terms in $\da^3 BR(z,\varpi)$}

%%%%%%%%%%%%%%%%%%%%%%%%%%%%%%%%%%%%%%%%%%%%%%%%%%%%%%%%%%%%%%%%%%%%%%%%%%%%%%%%%%%%%%%%%%%%%%%%%%%%%%%%%
%%%%%%%%%%%%%%%%%%%%%%%%%%%%%%%%%%%%%%%%%%%%%%%%%%%%%%%%%%%%%%%%%%%%%%%%%%%%%%%%%%%%%%%%%%%%%%%%%%%%%%%%%%

Since $z_t (\al)=BR(z,\varpi)(\al) + c(\al)\cdot\da z(\al)$ we have

\begin{align*}
\begin{split}
\int_{\T}\!\partial_{\al}^3 z(\al)\cdot \partial_{\al}^3
z_{t}(\al)d\al&=\int_{\T}\!\partial_{\al}^3 z(\al)\cdot
\partial_{\al}^3 BR(z,\varpi)(\al) d\al+\!\!
\int_{\T}\!\partial_{\al}^3 z(\al)\cdot \da^3(c(\al)\da z(\al)) d\al \\
&=I_1+I_2.
\end{split}
\end{align*}
Here and in 7.3 we study $I_1$. We shall estimate $I_2$ in 7.4.

Let us write $BR(z,\varpi)(\al)=P_1(\al)+P_2(\al)+P_{3}(\al)$ as in \eqref{p1p2p3}.
Then it is easy to check that
$$|\da^3 P_{3}(\al)|\leq C\|\varpi\|_{L^2}(|\da^3 z_2(\al)|+\|z\|^3_{C^2}),$$ giving us a term  controlled by
the energy estimate. The kernel in the integral $P_2(\al)$ has
order $1$ in $\beta$,  therefore taking two derivatives in $\al$
produces regular integrals as in $P_{3}$ and kernels of degree
$-1$ in $\beta$, for which we first exchange $\beta$ by
$\al-\beta$ and then take  one more derivative. We obtain kernels
of grade $-1$ in $\beta$ acting in $\varpi$ or $\varpi_{\al}$
which  can be estimated  as before. For the most singular term
$P_1(\al)$, we have

$$\int_\T \da^3 z(\al)\cdot\da^3 P_1(\al) d\al=I_3+I_4+I_5+I_6,$$
where

$$I_3=\frac{1}{4\pi}\int_\T\int_{\T}\da^3 z(\al)\da^3
\Big(\frac{V^{\bot}(\al,\al-\beta)}{|V(\al,\al-\beta)|^2}\Big)
\varpi(\al-\beta)d\beta,$$

$$I_4=\frac{3}{4\pi}\int_\T\int_{\T}\da^3 z(\al)\da^2
\Big(\frac{V^{\bot}(\al,\al-\beta)}{|V(\al,\al-\beta)|^2}\Big)
\da\varpi(\al-\beta)d\beta,$$

$$I_5=\frac{3}{4\pi}\int_\T\int_{\T}\da^3 z(\al)\da
\Big(\frac{V^{\bot}(\al,\al-\beta)}{|V(\al,\al-\beta)|^2}\Big)
\da^2\varpi(\al-\beta)d\beta,$$

$$I_6=\frac{1}{4\pi}\int_\T\int_{\T}\da^3 z(\al)
\Big(\frac{V^{\bot}(\al,\al-\beta)}{|V(\al,\al-\beta)|^2}\Big)
\da^3\varpi(\al-\beta)d\beta,$$

The most singular terms for $I_3$ are those in which three
derivatives appear and the kernels have degree $-1$. The rest of
the terms have  kernels with degree $k>-1$ and can be estimated as
before. One of the two singular terms of $I_3$ is given by

$$
J_1=\frac{1}{8\pi}\int_{\T}\int_{\T}\partial_{\al}^3
z(\al)\cdot\frac{(\da^3z(\al)-\da^3z(\al-\beta))^\bot}{|V(\al,\al-\beta)|^2}\varpi(\al-\beta)d\beta
d\alpha,$$ which we decompose as follows:
\begin{align*}
\begin{split}
J_1&=\frac{1}{8\pi}\int_{\T}\int_{\T}\partial_{\al}^3
z(\al)\cdot\frac{(\da^3z(\al)-\da^3z(\beta))^\bot}{|V(\al,\beta)|^2}\varpi(\beta)d\beta
d\alpha\\
&=\frac{1}{8\pi}\int_{\T}\int_{\T}\partial_{\al}^3
z(\al)\cdot\frac{(\da^3z(\al)-\da^3z(\beta))^\bot}{|V(\al,\beta)|^2}\frac{\varpi(\beta)+\varpi(\alpha)}{2}d\beta
d\alpha\\
&\quad+\frac{1}{8\pi}\int_{\T}\int_{\T}\partial_{\al}^3
z(\al)\cdot\frac{(\da^3z(\al)-\da^3z(\beta))^\bot}{|V(\al,\beta)|^2}\frac{\varpi(\beta)-\varpi(\alpha)}{2}d\beta
d\alpha\\
&=K_1+K_2.
\end{split}
\end{align*}
That is we have made a kind of integration by parts in $J_1$,
allowing us to show that the most singular term $K_1$ vanishes:
\begin{align*}
\begin{split}
K_1&=-\frac{1}{8\pi}\int_{\T}\int_{\T}\partial_{\al}^3
z(\beta)\cdot\frac{(\da^3z(\al)-\da^3z(\beta))^\bot}{|V(\al,\beta)|^2}\frac{\varpi(\beta)+\varpi(\alpha)}{2}d\beta
d\alpha\\
&=\frac{1}{16\pi}\int_{\T}\int_{\T}(\partial_{\al}^3
z(\al)-\partial_{\al}^3z(\beta))\cdot\frac{(\da^3z(\al)-\da^3z(\beta))^\bot}{|V(\al,\beta)|^2}\frac{\varpi(\beta)+\varpi(\alpha)}{2}d\beta d\alpha\\
&=0,
\end{split}
\end{align*}
whether for $K_2$ we have
\begin{align*}
\begin{split}
K_2&=\frac{1}{16\pi}\int_{\T}\int_{\T}\partial_{\al}^3
z(\al)\cdot(\da^3z(\beta))^\bot \frac{\varpi(\al)-\varpi(\beta)}{|V(\al,\beta)|^2}d\beta d\alpha\\
&=\int_{\T}\partial_{\al}^3 z(\al)\cdot (\da^3 Hz(\al))^\bot
\frac{\da \varpi(\al)}{|\da z(\al)|^2}d\al+
\frac{1}{4\pi}\int_{\T}\int_{\T}\partial_{\al}^3
z(\al)\cdot(\da^3z(\al-\beta))^\bot B_1(\al,\beta)d\beta d\alpha,\\
\end{split}
\end{align*}
where $|B_1(\al,\beta)|\leq
C\|\F(z)\|_{L^\infty}\|z\|_{C^2}\|\varpi\|_{C^{1,\delta}}|\beta|^{\delta-1}.$
The other singular term with three derivatives in $z(\al)$ and
kernel of degree $-1$ inside $I_3$ is given by
\begin{align*}
\begin{split}
J_2&=-\frac{1}{4\pi}\int_{\T}\int_{\T}\partial_{\al}^3
z(\al)\cdot\frac{V^\bot(\al,\al\!-\!\beta)}{|V(\al,\al\!-\!\beta)|^4}
V(\al,\al\!-\!\beta)\cdot(\da^3z(\al)\!-\!\da^3z(\al\!-\!\beta))\varpi(\al\!-\!\beta)d\beta
d\alpha
\end{split}
\end{align*}
Here we take $J_2=K_3+K_4+K_5$ where
\begin{align*}
\begin{split}
K_3&=-\frac{1}{4\pi}\int_{\T}\int_{\T}\partial_{\al}^3
z(\al)\cdot\frac{V^\bot(\al,\beta)}{|V(\al,\beta)|^4}
(V(\al,\beta)-W(\al,\beta))\cdot(\da^3z(\al)\!-\!\da^3z(\beta))\varpi(\beta)d\beta
d\alpha,
\end{split}
\end{align*}
\begin{align*}
\begin{split}
K_4&=-\frac{1}{4\pi}\int_{\T}\int_{\T}\partial_{\al}^3
z(\al)\cdot\frac{V^\bot(\al,\al\!-\!\beta)}{|V(\al,\al\!-\!\beta)|^4}
B_2(\al,\al\!-\!\beta)\cdot(\da^3z(\al)\!-\!\da^3z(\al\!-\!\beta))\varpi(\al\!-\!\beta)d\beta
d\alpha,
\end{split}
\end{align*}
with
$$
B_2(\al,\al-\beta)=W(\al,\al-\beta)-\da z(\al)\beta/2,
$$
and
$$
W(\al,\beta)=((\frac{z_1(\al)-z_1(\beta)}{2})_p,(\frac{z_2(\al)-z_2(\beta)}{2})_p)
$$
is defined in the appendix. Finally we have:
\begin{align*}
\begin{split}
K_5&=-\frac{1}{8\pi}\int_{\T}\int_{\T}\partial_{\al}^3
z(\al)\cdot\frac{V^\bot(\al,\al\!-\!\beta)\beta}{|V(\al,\al\!-\!\beta)|^4}
\da
z(\al)\cdot(\da^3z(\al)\!-\!\da^3z(\al\!-\!\beta))\varpi(\al\!-\!\beta)d\beta
d\alpha.
\end{split}
\end{align*}
The $L^\infty$ norm of
$$\frac{V^\bot(\al,\beta)}{|V(\al,\beta)|^4}
(V(\al,\beta)-W(\al,\beta))$$ is given in the appendix, allowing
us to estimate the term $K_3$ as before.

Next we split $K_4 = L_1 + L_2$, where
\begin{align*}
\begin{split}
L_1&=-\frac{1}{4\pi}\int_{\T}\int_{\T}\partial_{\al}^3
z(\al)\cdot\frac{V^\bot(\al,\al\!-\!\beta)}{|V(\al,\al\!-\!\beta)|^4}
B_2(\al,\al\!-\!\beta)\cdot\da^3z(\al)\varpi(\al\!-\!\beta)d\beta
d\alpha,
\end{split}
\end{align*}
and
\begin{align*}
\begin{split}
L_2&=\frac{1}{4\pi}\int_{\T}\int_{\T}\partial_{\al}^3
z(\al)\cdot\frac{V^\bot(\al,\al\!-\!\beta)}{|V(\al,\al\!-\!\beta)|^4}
B_2(\al,\al\!-\!\beta)\cdot\da^3z(\al\!-\!\beta)\varpi(\al\!-\!\beta)d\beta
d\alpha,
\end{split}
\end{align*}

We have
\begin{align*}
\begin{split}
|L_1|&\leq C\big|\int_{\T}\partial_{\al}^3
z(\al)\cdot\frac{\da^{\bot} z(\al)}{|\da z(\al)|^4}
\da^2 z(\al)\cdot\da^3z(\al)H\varpi(\al)d\beta
d\alpha\big|\\
&\quad +\|\da^3 z\|^2_{L^2}\|\F(z)\|_{L^\infty}\|z\|_{C^{2,\delta}}\|\varpi\|_{L^\infty}
\end{split}
\end{align*}

\begin{align*}
\begin{split}
|L_2|&\leq C\big|\int_{\T}\partial_{\al}^3
z(\al)\cdot\frac{\da^{\bot} z(\al)}{|\da z(\al)|^4}\varpi(\al)
\da^2 z(\al)\cdot H(\da^3z)(\al) d\beta
d\alpha\big|\\
&\quad +\|\da^3 z\|^2_{L^2}\|\F(z)\|_{L^\infty}\|z\|_{C^{2,\delta}}\|\varpi\|_{C^1}
\end{split}
\end{align*}
and the term $K_4$ is  controlled.

 For $K_5$ we split
\begin{align*}
\begin{split}
K_5&=\frac{1}{8\pi}\int_{\T}\int_{\T}\partial_{\al}^3
z(\al)\cdot\frac{V^\bot(\al,\al\!-\!\beta)\beta}{|V(\al,\al\!-\!\beta)|^4}
(\da z(\al)-\da z(\al\!-\!\beta))\cdot\da^3z(\al\!-\!\beta)\varpi(\al\!-\!\beta)d\beta
d\alpha,\\
&\quad -\frac{1}{8\pi}\int_{\T}\int_{\T}\partial_{\al}^3
z(\al)\cdot B_{3}(\al,\al\!-\!\beta)(\da z(\al)\cdot \da^3 z(\al)-\da z(\al\!-\!\beta)\cdot\da^3z(\al\!-\!\beta))d\beta
d\alpha\\
&=L_3+L_4
\end{split}
\end{align*}
where
$$
B_3(\al,\al\!-\!\beta)=\frac{V^\bot(\al,\al\!-\!\beta)\varpi(\al\!-\!\beta)\beta}{|V(\al,\al\!-\!\beta)|^4}.
$$
Then we have
\begin{align*}
\begin{split}
|L_3|&\leq C\big|\int_{\T}\partial_{\al}^3
z(\al)\cdot\frac{\da^{\bot} z(\al)}{|\da z(\al)|^4}
\da^2 z(\al)\cdot H(\da^3z\, \varpi)(\al)d\alpha|\\
&\quad +\|\da^3
z\|^2_{L^2}\|\F(z)\|_{L^\infty}\|z\|_{C^{2,\delta}}\|\varpi\|_{L^{\infty}}.
\end{split}
\end{align*}

For $L_4$ we use an appropriated  integration by part: $$\da
z(\al)\cdot \da^3 z(\al)=-|\da^2 z(\al)|^2,$$ to obtain
\begin{align*}
\begin{split}
L_4&=\frac{1}{8\pi}\int_{\T}\int_{\T}\partial_{\al}^3
z(\al)\cdot B_{3}(\al,\al\!-\!\beta)(|\da^2 z(\al)|^2-|\da^2 z(\al\!-\!\beta)|^2)d\beta
d\alpha.
\end{split}
\end{align*}
Next we write $L_4=M_1+M_2,$ with
\begin{align*}
\begin{split}
M_1&=\frac{1}{8\pi}\int_{\T}\int_{\T}\partial_{\al}^3
z(\al)\cdot C(\al,\al-\beta)(|\da^2 z(\al)|^2-|\da^2 z(\al\!-\!\beta)|^2)d\beta d\alpha
\end{split}
\end{align*}
for
\begin{align*}
C(\al,\al-\beta)&=B_{3}(\al,\al\!-\!\beta)-\frac{2\da^{\bot} z(\al)\varpi(\al)}{|\da^2 z(\al)|^4\sin^2(\beta/2)}\\
&=\frac{V^\bot(\al,\al\!-\!\beta)\varpi(\al\!-\!\beta)\beta}{|V(\al,\al\!-\!\beta)|^4}-\frac{2\da^{\bot} z(\al)\varpi(\al)}{|\da^2 z(\al)|^4\sin^2(\beta/2)},
\end{align*}
and
\begin{align*}
\begin{split}
M_2&=C\int_{\T}\int_{\T}\partial_{\al}^3
z(\al)\cdot\da^{\bot} z(\al) \frac{\varpi(\al)}{|\da^2 z(\al)|^4}\la(|\da^2 z |^2)d\alpha.
\end{split}
\end{align*}
Since
$$
|C(\al,\al-\beta)|\leq \|\F(z)\|_{L^\infty}\|z\|_{C^2}\|\varpi\|_{C^1}\frac{1}{|\beta|}
$$
(see lemma \ref{C} in the appendix for more details)
and
$$
||\da^2 z(\al)|^2-|\da^2 z(\al\!-\!\beta)|^2|\leq 2\|z\|_{C^1}|\beta|\int_0^1|\da^3 z(\al+(s-1)\beta)| ds
$$
we get $|M_1|\leq
\|\F(z)\|_{L^\infty}\|z\|_{C^2}\|\varpi\|_{C^1}\|\da^3
z\|^2_{L^2}$.

For the term $M_2$ we use the estimate $$\|\la (|\da^2 z
|^2)\|_{L^2}=\|\da (|\da^2 z |^2)\|_{L^2}\leq 2\|\da^2 z
\|_{L^\infty}\|\da^3 z \|_{L^2},$$ to obtain $|M_2|\leq
C\|\F(z)\|_{L^\infty}\|\varpi\|_{L^{\infty}}\|\da^2 z
\|_{L^\infty}\|\da^3 z \|^2_{L^2}$.

For $I_4$ the most singular terms are those for which two
derivatives are applied to $z(\al)$. One of those is $J_3$
$$
J_3=C\int_\T\int_{\T}\da^3 z(\al)\cdot(\da^2z(\al)-\da^2z(\al-\beta))\frac{\da\varpi(\al-\beta)}{|V(\al,\al-\beta)|^2}
d\beta.
$$
We split $J_3=K_6+K_7+K_8$
$$
K_6=C\int_\T\int_{\T}\da^3 z(\al)\cdot(\da^2z(\al)-\da^2z(\al-\beta))\frac{\da\varpi(\al-\beta)-\da\varpi(\al)}{|V(\al,\al-\beta)|^2}
d\beta,
$$
$$
K_7=C\int_\T\int_{\T}\da\varpi(\al)\da^3 z(\al)\cdot(\da^2z(\al)-\da^2z(\al-\beta))(\frac{1}{|V(\al,\al-\beta)|^2}-\frac{4}{|\da z(\al)|^2\beta^2})
d\beta,
$$
$$
K_8=C\int_\T\da\varpi(\al)\da^3 z(\al)\cdot \int_{\T}\frac{\da^2z(\al)-\da^2z(\al-\beta)}{\beta^2}
d\beta,
$$
Using that \begin{equation}\label{ntni}\da^2z(\al)-\da^2z(\al-\beta)=\beta \int_0^1 \da^3 z(\al+(s-1)\beta) ds\end{equation} and $|\da\varpi(\al-\beta)-\da\varpi(\al)|\leq \|w\|_{C^{1,\delta}}|\beta|^\delta$ we have
\begin{align*}
|K_6|&\leq C\|\F(z)\|_{L^{\infty}}\|w\|_{C^{1,\delta}}\int_0^1\int_{\T} |\beta|^{\delta-1} \int_{\T}
|\da^3 z(\al)||\da^3 z(\al+(s-1)\beta)| d\alpha d\beta ds \\
&\leq C\|\F(z)\|_{L^{\infty}}\|w\|_{C^{1,\delta}}\int_0^1\int_{\T} |\beta|^{\delta-1} \int_{\T}
(|\da^3 z(\al)|^2+|\da^3 z(\al+(s-1)\beta)|^2) d\alpha d\beta ds\\
&\leq C\|\F(z)\|_{L^{\infty}}\|w\|_{C^{1,\delta}}\|\da^3 z\|^2_{L^3}.
\end{align*}
Due to \eqref{ntni} and the estimates obtained in the appendix we
have
$$
|K_7|\leq C \|\F(z)\|_{L^{\infty}}\|w\|_{C^1}\|z\|_{C^2}\|\da^3 z\|^2_{L^2}.
$$
Then using that $1/\beta-1/2\sin(\beta/2)$ is bounded, we get
$$
K_8\leq C \|\F(z)\|_{L^{\infty}}\|w\|_{C^1}\|\da^3 z\|^2_{L^3}.
$$

Regarding $I_5$, its most singular term is giving by
$$
J_4=C\int_\T\int_{\T}\da^3 z(\al)\cdot(\da z(\al)-\da
z(\al-\beta))\frac{\da^2\varpi(\al-\beta)}{|V(\al,\al-\beta)|^2}
d\beta,
$$
which after being decomposed in the form
\begin{align*}
J_4&=C\int_\T\int_{\T}\da^3 z(\al)\cdot\Big( \frac{\da z(\al)-\da z(\al-\beta)}{|V(\al,\al-\beta)|^2}-\frac{2\da^2z(\al)}{|\da z(\al)|^2\tan(\beta/2)}\Big)\da^2\varpi(\al-\beta)
d\beta\\
&\quad+C\int_\T\da^3 z(\al)\cdot\frac{\da^2 z(\al)}{|\da z(\al)|^2}H(\da^2\varpi)(\al)d\al.
\end{align*}
 can be estimated as before  $|J_4|\leq
C\|\F(z)\|_{L^{\infty}}\|z\|_{C^{2,\delta}}\|\da^2
w\|_{L^2}\|\da^3 z\|_{L^2}$.

%%%%%%%%%%%%%%%%%%%%%%%%%%%%%%%%%%%%%%%%%%%%%%%%%%%%%%%%%%%%%%%%%%%%%%%%%%%%%%%%%%%%%%%%%%%%%%%%%%%%%%%%%%%%%%%%%
%%%%%%%%%%%%%%%%%%%%%%%%%%%%%%%%%%%%%%%%%%%%%%%%%%%%%%%%%%%%%%%%%%%%%%%%%%%%%%%%%%%%%%%%%%%%%%%%%%%%%%%%%%%%%%%%%

\subsection{Looking for $\sigma(\al)$}

%%%%%%%%%%%%%%%%%%%%%%%%%%%%%%%%%%%%%%%%%%%%%%%%%%%%%%%%%%%%%%%%%%%%%%%%%%%%%%%%%%%%%%%%%%%%%%%%%%%%%%%%%%%%%%%%%
%%%%%%%%%%%%%%%%%%%%%%%%%%%%%%%%%%%%%%%%%%%%%%%%%%%%%%%%%%%%%%%%%%%%%%%%%%%%%%%%%%%%%%%%%%%%%%%%%%%%%%%%%%%%%%%%%

The term  $I_6$ will gives us the proper sign (Rayleigh-Taylor
condition) that has to be imposed upon $\sigma(\al)$. Let us recall the formula
$$\sigma(\alpha,t)=\frac{\mu^2-\mu^1}{\kappa}BR(z,\varpi)(\al,t)\cdot\dpa
z(\al,t)+ \mathrm{g}(\rho^2-\rho^1)\da z_1(\al,t)$$
We write $I_6$ in the form $I_6 = J_1 + J_2$ where
$$
J_1=\frac{1}{4\pi}\int_\T\int_{\T}\da^3 z(\al)\cdot
\Big(\frac{V^{\bot}(\al,\al-\beta)}{|V(\al,\al-\beta)|^2}-\frac{(\da z(\al))^\bot}{|\da z(\al)|^2\tan(\beta/2)}\Big)
\da^3\varpi(\al-\beta)d\beta d\alpha,
$$
and
$$
J_2=\frac{1}{4\pi}\int_\T\int_{\T}\da^3 z(\al)\cdot
\frac{(\da z(\al))^\bot}{|\da z(\al)|^2}H(\da^3\varpi)(\al)d\beta d\alpha.
$$
Let us denote the kernel  of $J_1$ by $\Sigma(\al,\al-\beta)$,
which is of degree $0$ in $\beta$. After an integration by parts
we obtain:
\begin{align*}
J_1&=-\frac{1}{4\pi}\int_\T\int_{\T}\da^3 z(\al)\cdot
\Sigma(\al,\al-\beta)\partial_\beta(\da^2\varpi(\al-\beta))d\beta d\alpha\\
&=\frac{1}{4\pi}\int_\T\int_{\T}\da^3 z(\al)\cdot
\partial_\beta\Sigma(\al,\al-\beta)
\da^2\varpi(\al-\beta)d\beta d\alpha.
\end{align*}
Then $\partial_\beta\Sigma(\al,\al-\beta)$  has terms of degree
$0$ which are estimated easily. The term with degree $-1$ is given
by
$$
\frac{(\da z(\al-\beta))^{\bot}}{2|V(\al,\al-\beta)|^2}+\frac{(\da
z(\al))^{\bot}}{2|\da
z(\al)|^2\sin^2(\beta/2)}-\frac{V^{\bot}(\al,\al-\beta)}{|V(\al,\al-\beta)|^4}V(\al,\al-\beta)\cdot
\da z(\al-\beta),
$$
 and we decompose it as a sum of a kernel of degree $0$ (easy to
 estimate)
$$
-\frac{(\da z(\al))^\bot}{|\da z(\al)|^2}\big(\frac{2}{|\beta|^2}-\frac{1}{2\sin^2(\beta/2)}\big),
$$
and  six kernels of degree $-1$, ($P_1$,..., $P_6$) given by
$$
P_1(\al,\al-\beta)=\frac{(\da z(\al-\beta)-\da z(\al))^{\bot}}{2|V(\al,\al-\beta)|^2},
$$
$$
P_2(\al,\al-\beta)=-\frac{(\da z(\al))^{\bot}}{2}\big(\frac{1}{|V(\al,\al-\beta)|^2}-\frac{4}{|\da z(\al)|^2|\beta|^2}\big),
$$
$$
P_3(\al,\al-\beta)=\frac{V^{\bot}(\al,\al-\beta)}{|V(\al,\al-\beta)|^4}V(\al,\al-\beta)\cdot (\da z(\al)-\da z(\al-\beta)),
$$
$$
P_4(\al,\al-\beta)=-\frac{V^{\bot}(\al,\al-\beta)}{|V(\al,\al-\beta)|^4}(V(\al,\al-\beta)-\da z(\al)\beta/2)\cdot \da z(\al),
$$
$$
P_5(\al,\al-\beta)=-\frac{|\da z(\al)|^2\beta}{2}\frac{V^{\bot}(\al,\al-\beta)-\da^{\bot} z(\al)\beta/2}{|V(\al,\al-\beta)|^4},
$$
$$
P_6(\al,\al-\beta)=-\frac{\da^{\bot} z(\al)|\da z(\al)|^2|\beta|^2}{4|V(\al,\al-\beta)|^2}\big(\frac{1}{|V(\al,\al-\beta)|^2}-\frac{4}{|\da z(\al)|^2|\beta|^2}\big).
$$
To control the term with kernel $P_2$ we consider $P_2=Q_1+Q_2$
$$
Q_1(\al,\al-\beta)=P_2(\al,\al-\beta)-\frac{(\da z(\al))^{\bot}}{2}\frac{2\da z(\al)\cdot \da^2 z(\al)}{|\da z(\al)|^4\beta},
$$
$$
Q_2(\al,\al-\beta)=-\da^{\bot} z(\al)\Big(\frac{\da z(\al)\cdot \da^2 z(\al)}{|\da z(\al)|^4}\big(\frac{1}{\beta}-\frac{1}{2\tan(\beta/2)}\big)+
\frac{\da z(\al)\cdot \da^2 z(\al)}{|\da z(\al)|^42\tan(\beta/2)}\Big).
$$
It is shown in the appendix that $\|Q_1\|_{L^{\infty}}\leq
\|\F(z)\|^k_{L^\infty}\|z\|^k_{C^{2,\delta}}|\beta|^{\delta-1}$, (see lemma \ref{Q1})
giving us
$$
\frac{1}{4\pi}\int_\T\int_\T\da^3 z(\al)\cdot Q_1(\al,\al-\beta)
\da^2\varpi(\al-\beta)d\beta d\al \leq C\|\F(z)\|^k_{L^\infty}\|z\|^k_{C^{2,\delta}}\|\da^3 z\|_{L^2}\|\da^2 \varpi\|_{L^2}.
$$
The integral
$$
K_1=\frac{1}{4\pi}\int_\T\int_\T\da^3 z(\al)\cdot Q_2(\al,\al-\beta)
\da^2\varpi(\al-\beta)d\beta d\al
$$
is bounded by
\begin{align*}
|K_1|&\leq C\|\F(z)\|^{3/2}_{L^\infty}\|z\|_{C^2}\int_{\T}|\da^3 z(\al)|(\|\da^2 \varpi\|_{L^2}+|H(\da^2\varpi)(\al)|) d\al\\
&\leq C\|\F(z)\|^{3/2}_{L^\infty}\|z\|_{C^2}\|\da^3 z\|_{L^2}\|\da^2 \varpi\|_{L^2}.
\end{align*}
It is now very clear that the other $P_i$ terms can be estimated
as above or as before i.e. we finally have
$$
J_1\leq \exp C(\|\F(z)\|^2_{L^\infty}+\|z\|^2_{H^3})
$$
We consider now the $J_2$ term which can be written as follows
\begin{align*}
J_2&=\frac{1}{4\pi A(t)}\int_\T\da^3 z(\al)\cdot
\da^{\bot} z(\al)\la(\da^2\varpi)(\al) d\alpha=\frac{1}{4\pi A(t)}\int_\T\la(\da^3 z\cdot
\da^{\bot} z)(\al)\da^2\varpi(\al) d\alpha
\end{align*}
and using the formula \eqref{aux} we separate $J_2$ as a sum of
two parts, $K_2$ and $K_3$, where
\begin{align*}
K_2=-\frac{\kappa \mathrm{g}(\rho^2-\rho^1)}{2\pi (\mu^2+\mu^1)A(t)}\int_\T\la(\da^3 z\cdot
\da^{\bot} z)(\al)\da^3z_2(\al) d\alpha,
\end{align*}
and
\begin{align*}
K_3=-\frac{A_\mu}{4\pi A(t)}\int_\T\la(\da^3 z\cdot
\da^{\bot} z)(\al)\da^2T(\varpi)(\al) d\alpha.
\end{align*}
For $K_2$ we decompose further $K_2=L_1+L_2$, where
\begin{align*}
L_1=\frac{\kappa \mathrm{g}(\rho^2-\rho^1)}{2\pi (\mu^2+\mu^1)A(t)}\int_\T\la(\da^3 z_1
\da z_2)(\al)\da^3z_2(\al) d\alpha,
\end{align*}
and
\begin{align*}
L_2=-\frac{\kappa \mathrm{g}(\rho^2-\rho^1)}{2\pi (\mu^2+\mu^1)A(t)}\int_\T\la(\da^3 z_2
\da z_1)(\al)\da^3z_2(\al) d\alpha.
\end{align*}
Then $L_1$ is written as $L_1=M_1+M_2$ with
$$
M_1=\frac{\kappa \mathrm{g}(\rho^2-\rho^1)}{2\pi (\mu^2+\mu^1)A(t)}\int_\T(\la(\da^3 z_1
\da z_2)(\al)-\la(\da^3 z_1)(\al)\da z_2(\al))\da^3z_2(\al) d\alpha,
$$
and
$$
M_2=\frac{\kappa \mathrm{g}(\rho^2-\rho^1)}{2\pi (\mu^2+\mu^1)A(t)}\int_\T\la(\da^3 z_1)(\al)\da z_2(\al)\da^3z_2(\al) d\alpha.
$$
Using the commutator estimate, we get  $$M_1\leq
C\|\F(z)\|_{L^\infty}\|z\|_{C^{2,\delta}}\|\da^3 z\|^2_{L^2}\leq
\exp C(\|\F(z)\|^2_{L^\infty}+\| z\|^2_{H^3})$$ The identity
$$
\da z_2(\al)\da^3z_2(\al)=-\da z_1(\al)\da^3z_1(\al)-|\da^2 z(\al)|^2,
$$
lets us write $M_2$ as the sum of $N_1$ and $N_2$ where
$$
N_1=-\frac{\kappa \mathrm{g}(\rho^2-\rho^1)}{2\pi (\mu^2+\mu^1)A(t)}\int_\T\la(\da^3 z_1)(\al)|\da^2 z(\al)|^2 d\alpha,
$$
and
$$
N_2=-\frac{\kappa \mathrm{g}(\rho^2-\rho^1)}{2\pi (\mu^2+\mu^1)A(t)}\int_\T\da z_1(\al)\da^3z_1(\al)\la(\da^3 z_1)(\al) d\alpha.
$$
Integration by parts shows that $$N_1\leq
C\|\F(z)\|_{L^\infty}\|z\|_{C^2}\|\da^3 z\|^2_{L^2}\leq \exp
C(\|\F(z)\|^2_{L^\infty}+\| z\|^2_{H^3}).$$ Writing $L_2$ in the
form:
\begin{align*}
L_2=-\frac{\kappa \mathrm{g}(\rho^2-\rho^1)}{2\pi (\mu^2+\mu^1)A(t)}\int_\T\da z_1(\al)\da^3 z_2(\al)
\la(\da^3z_2)(\al) d\alpha,
\end{align*}
we obtain finally
\begin{align*}
K_2\leq \exp C(\|\F(z)\|^2_{L^\infty}+\| z\|^2_{H^3})-\frac{2\kappa \mathrm{g}(\rho^2-\rho^1)}{4\pi (\mu^2+\mu^1)A(t)}\int_\T\da z_1(\al)\da^3 z(\al)\cdot\la(\da^3z)(\al) d\alpha.
\end{align*}
In the  estimate above we can observe how a part of $\sigma(\al)$
 appears in the non-integrable terms. Let us now return to
 $K_3=L_3+L_4+L_5$, where
\begin{align*}
L_3=-\frac{A_\mu}{4\pi A(t)}\int_\T\la(\da^3 z\cdot
\da^{\bot} z)(\al)(2\da^2 BR(z,\varpi)(\al))\cdot\da z(\al) d\alpha,
\end{align*}
\begin{align*}
L_4=-\frac{A_\mu}{2\pi A(t)}\int_\T\la(\da^3 z\cdot
\da^{\bot} z)(\al)\da BR(z,\varpi)(\al))\cdot\da^2 z(\al) d\alpha,
\end{align*}
and
\begin{align*}
L_5=-\frac{A_\mu}{2\pi A(t)}\int_\T\la(\da^3 z\cdot
\da^{\bot} z)(\al) BR(z,\varpi)(\al)\cdot\da^3 z(\al) d\alpha.
\end{align*}
We will control first the terms $L_3$ and $L_4$ and then we will
show how the rest of $\sigma(\al)$ appears in $L_5$. Integrating
by parts in $L_4$ and we obtain
$$
L_4\leq C\|\F(z)\|_{L^\infty}\|H(\da^3 z\cdot
\da^{\bot} z)\|_{L^2}(\|\da^2 BR(z,\varpi)\|_{L^2}\|\da^2 z\|_{L^\infty}+\|\da BR(z,\varpi)\|_{L^\infty}\|\da^3 z\|_{L^2})
$$
and using the estimates for $\|\da^2 BR(z,\varpi)\|_{L^2}$, we get
$L_4\leq \exp C(\|\F(z)\|^2_{L^\infty}+\| z\|^2_{H^3})$ With $L_3$
we also integrate by parts to obtain $L_3=M_3+M_4$ where
\begin{align*}
M_3=-\frac{A_\mu}{4\pi A(t)}\int_\T H(\da^3 z\cdot
\da^{\bot} z)(\al)(2\da^2 BR(z,\varpi)(\al))\cdot\da^2 z(\al) d\alpha,
\end{align*}
and
\begin{align*}
M_4=-\frac{A_\mu}{4\pi A(t)}\int_\T H(\da^3 z\cdot
\da^{\bot} z)(\al)(2\da^3 BR(z,\varpi)(\al))\cdot\da z(\al) d\alpha.
\end{align*}
Easily we have
$$
M_3\leq C\|\F(z)\|_{L^\infty}\|H(\da^3 z\cdot
\da^{\bot} z)\|_{L^2}\|\da^2 BR(z,\varpi)\|_{L^2}\|\da^2 z\|_{L^\infty}\leq \exp C(\|\F(z)\|^2_{L^\infty}\!\!+\| z\|^2_{H^3}).
$$
In $M_4$ the application of Leibniz's rule to $\da^3 BR(z,\varpi)$
produces many terms which can be estimated with the same tools
used before with $I_4$ and $I_5$. For the most singular terms we
have the expressions:
$$
N_3=-\frac{A_\mu}{4\pi A(t)}\int_\T H(\da^3 z\cdot
\da^{\bot} z)(\al)2 \da (BR(z,\da^2\varpi)(\al))\cdot\da z(\al) d\alpha,
$$
$$
N_4=-\frac{A_\mu}{2\pi A(t)}\int_\T H(\da^3 z\cdot
\da^{\bot} z)(\al)\int_{\T}\varpi(\al-\beta)\frac{(\da^3 z(\al)-\da^3z(\al-\beta))^\bot}{|V(\al,\al-\beta)|^2}d\beta\cdot\da z(\al) d\alpha,
$$
$$
N_5=\frac{A_\mu}{2\pi A(t)}\int_\T H(\da^3 z\cdot
\da^{\bot} z)(\al)\int_{\T}\varpi(\al-\beta)B(\al,\al-\beta)\cdot (\da^3z(\al)-\da^3z(\al-\beta)))d\beta d\alpha,
$$
where
$$
B(\al,\al-\beta)=\frac{V^{\bot}(\al,\al-\beta)\cdot\da z(\al)}{|V(\al,\al-\beta)|^4}V(\al,\al-\beta).
$$
Let us consider
\begin{align*}
\da (BR(z,\da^2\varpi)(\al))\cdot\da z(\al)&=
\da(BR(z,\da^2\varpi)(\al))\cdot\da z(\al))-BR(z,\da^2\varpi)(\al)\cdot\da^2 z(\al)\\
&=\frac12\da (T(\da^2\varpi)(\al))-BR(z,\da^2\varpi)(\al)\cdot\da^2 z(\al)
\end{align*}
which yields
$$
N_3\leq C\|\F(z)\|_{L^\infty}\|H(\da^3 z\cdot
\da^{\bot} z)\|_{L^2}(\|T(\da^2\varpi)\|_{H^1}+\| BR(z,\da^2\varpi)\|_{L^2}\|\da^2 z\|_{L^\infty})
$$
and therefore
$$
N_3\leq \exp C(\|\F(z)\|^2_{L^\infty}\!\!+\| z\|^2_{H^3}).
$$
Next  we write $N_4=O_1+O_2+O_3+O_4+O_5,$
$$
O_1=-\frac{A_\mu}{2\pi A(t)}\int_\T H(\da^3 z\cdot
\da^{\bot}z)(\al)(\da^3 z(\al))^\bot\cdot\da z(\al) \int_{\T}\big(\frac{\varpi(\al-\beta)}{|V(\al,\al-\beta)|^2}-\frac{4\varpi(\al)}{|\da z(\al)|^2|\beta|^2}\big)d\beta d\alpha,
$$
$$
O_2=\frac{A_\mu}{2\pi A(t)}\int_\T H(\da^3 z\cdot
\da^{\bot}z)(\al) \int_{\T}(\da^3 z(\al-\beta))^\bot\cdot\da z(\al)\big(\frac{\varpi(\al-\beta)}{|V(\al,\al-\beta)|^2}-\frac{4\varpi(\al)}{|\da z(\al)|^2|\beta|^2}\big)d\beta d\alpha,
$$
$$
O_3=-\frac{A_\mu}{2\pi A^2(t)}\int_\T H(\da^3 z\cdot
\da^{\bot}z)(\al)\varpi(\al)(\da^3 z(\al))^\bot\cdot\da z(\al) \int_{\T}\big(\frac{4}{|\beta|^2}-\frac{1}{\sin^2(\beta/2)}\big)d\beta d\alpha,
$$
$$
O_4=\frac{A_\mu}{2\pi A^2(t)}\int_\T H(\da^3 z\cdot
\da^{\bot}z)(\al)\varpi(\al) \int_{\T}(\da^3 z(\al-\beta))^\bot\cdot\da z(\al)\big(\frac{4}{|\beta|^2}-\frac{1}{\sin^2(\beta/2)}\big)d\beta d\alpha,
$$
$$
O_5=-\frac{A_\mu}{2\pi A^2(t)}\int_\T H(\da^3 z\cdot
\da^{\bot}z)(\al)\varpi(\al)\da z(\al)\cdot\la((\da^3z)^\bot)(\al) d\alpha,
$$
The terms $O_1$, $O_2$, $O_3$ and $O_4$ can be estimated as before. Then we split $O_5=R_1+R_2$ where
$$
R_1=\frac{A_\mu}{2\pi A^2(t)}\int_\T H(\da^3 z\cdot
\da^{\bot}z)(\al)(\la(\varpi\da
z\cdot(\da^3z)^\bot)(\al)-\varpi(\al)\da
z(\al)\cdot\la((\da^3z)^\bot)(\al)) d\alpha,
$$
and
$$
R_2=-\frac{A_\mu}{2\pi A^2(t)}\int_\T H(\da^3 z\cdot
\da^{\bot}z)(\al)\la(\varpi\da z\cdot(\da^3z)^\bot)(\al) d\alpha.
$$
Using the commutator estimate, we obtain
$$
R_1\leq C\|\F(z)\|^2_{L^\infty}\|H(\da^3 z\cdot
\da^{\bot}z)\|_{L^2}\|\varpi\da
z\|_{C^{1,\delta}}\|\da^3z\|_{L^2}\leq\exp
C(\|\F(z)\|^2_{L^\infty}\!\!+\| z\|^2_{H^3}).
$$
The identity $\la (H)=-\da$
gives
\begin{align*}
R_2&=\frac{A_\mu}{2\pi A^2(t)}\int_\T \da(\da^3 z\cdot
\da^{\bot}z)(\al)\varpi(\al)\da z(\al)\cdot(\da^3z(\al))^\bot d\alpha\\
&=-\frac{A_\mu}{2\pi A^2(t)}\int_\T \da(\da^3 z\cdot
\da^{\bot}z)(\al)\varpi(\al)\da^3z(\al)\cdot\da^\bot z(\al) d\alpha
\end{align*}
and integrating by parts we get
$$
R_2\leq C\|\F(z)\|^2_{L^\infty}\|\da^3 z\cdot
\da^{\bot}z\|^2_{L^2}\|\da\varpi\|_{L^{\infty}}\leq \exp
C(\|\F(z)\|^2_{L^\infty}\!\!+\| z\|^2_{H^3}).
$$
Regarding the term $N_5$ we have the expression
$$
B(\al,\al-\beta)=\frac{(V(\al,\al-\beta)-\da z(\al)\beta/2)^{\bot}\cdot\da z(\al)}{|V(\al,\al-\beta)|^4}V(\al,\al-\beta)
$$
which shows that $B(\al,\al-\beta)$ has order $-1$ and, therefore,
the term $N_5$ can be estimated as before.

Finally we have to find $\sigma(\al)$ in $L_5$ to finish the proof
of the lemma. To do that let us  split $L_5 =M_5 +M_6+M_7+M_8$
where
$$
M_5=\frac{A_\mu}{2\pi A(t)}\int_\T\la(\da^3 z_1\da z_2)(\al) BR_1(z,\varpi)(\al)\da^3 z_1(\al) d\alpha,
$$
$$
M_6=\frac{A_\mu}{2\pi A(t)}\int_\T\la(\da^3 z_1\da z_2)(\al) BR_2(z,\varpi)(\al)\da^3 z_2(\al) d\alpha,
$$
$$
M_7=-\frac{A_\mu}{2\pi A(t)}\int_\T\la(\da^3 z_2\da z_1)(\al) BR_1(z,\varpi)(\al)\da^3 z_1(\al) d\alpha,
$$
$$
M_8=-\frac{A_\mu}{2\pi A(t)}\int_\T\la(\da^3 z_2\da z_1)(\al) BR_2(z,\varpi)(\al)\da^3 z_2(\al) d\alpha.
$$
and $BR_j$, $j=1,2$, is the $j$th-component of the Birkhoff-Rott
integral.

Then
\begin{align*}
M_5&=\frac{A_\mu}{2\pi A(t)}\int_\T(\la(\da z_2\da^3 z_1)(\al)-\da z_2(\al)\la(\da^3 z_1)(\al)) BR_1(z,\varpi)(\al)\da^3 z_1(\al) d\alpha\\
&\quad +\frac{A_\mu}{2\pi A(t)}\int_\T\da z_2(\al)BR_1(z,\varpi)(\al)\da^3 z_1(\al)\la(\da^3 z_1)(\al) d\alpha,
\end{align*}
and the commutator estimates yields
\begin{equation}\label{m5}
M_5\leq \exp C(\|\F(z)\|^2_{L^\infty}\!\!+\!\| z\|^2_{H^3})\!+\!\frac{A_\mu}{2\pi A(t)}\!\int_\T\!BR_1(z,\varpi)(\al)\da z_2(\al)\da^3 z_1(\al)\la(\da^3 z_1)(\al) d\alpha.
\end{equation}
In a similar way we have
\begin{equation*}
M_6\leq \exp C(\|\F(z)\|^2_{L^\infty}\!\!+\!\| z\|^2_{H^3})\!+\!\frac{A_\mu}{2\pi A(t)}\!\int_\T\!BR_2(z,\varpi)(\al)\da z_2(\al)\da^3 z_2(\al)\la(\da^3 z_1)(\al) d\alpha.
\end{equation*}
Let us introduce the notation
$$
N_4=\frac{A_\mu}{2\pi A(t)}\!\int_\T\!BR_2(z,\varpi)(\al)\da z_2(\al)\da^3 z_2(\al)\la(\da^3 z_1)(\al) d\alpha,
$$
 the  equality $\da z_2(\al)\da^3 z_2(\al)=-\da z_1(\al)\da^3
z_1(\al)-|\da^2z(\al)|^2$ gives $N_4=O_6+O_7$ where
$$
O_6=-\frac{A_\mu}{2\pi A(t)}\!
\int_\T\!BR_2(z,\varpi)(\al)|\da^2z(\al)|^2\la(\da^3 z_1)(\al) d\alpha,
$$
$$
O_7=-\frac{A_\mu}{2\pi A(t)}\!\int_\T\!BR_2(z,\varpi)(\al)\da z_1(\al)\da^3 z_1(\al)\la(\da^3 z_1)(\al) d\alpha.
$$
Integrating by parts in $O_6$ we get
\begin{align*}
O_6&\leq \!C\|\F(z)\|_{L^\infty}(\|\da BR(z,\varpi)\|_{L^\infty}\|\da^2 z\|^2_{L^\infty}\!\!+\|BR(z,\varpi)\|_{L^\infty}\|\da^2 z\|_{L^\infty}\|\da^3 z\|_{L^2})\|H(\da^3 z_1)\|_{L^2}\\
&\leq \exp C(\|\F(z)\|^2_{L^\infty}\!\!+\!\| z\|^2_{H^3}).
\end{align*}
Finally we get the estimate
$$
M_6\leq\!\exp C(\|\F(z)\|^2_{L^\infty}\!\!+\!\| z\|^2_{H^3})-\frac{A_\mu}{2\pi A(t)}\!\int_\T\!BR_2(z,\varpi)(\al)\da z_1(\al)\da^3 z_1(\al)\la(\da^3 z_1)(\al) d\alpha.
$$
Which together with \eqref{m5} yields
$$
M_5+M_6\leq \!\exp C(\|\F(z)\|^2_{L^\infty}\!\!+\!\| z\|^2_{H^3})
-\frac{A_\mu}{2\pi A(t)}\!\int_\T\!BR(z,\varpi)(\al)\cdot\da^{\bot} z(\al)\,\da^3 z_1(\al)\la(\da^3 z_1)(\al) d\alpha.
$$ With $M_7$ and $M_8$ we use the equality $\da z_1(\al)\da^3 z_1(\al)=-\da z_2(\al)\da^3 z_2(\al)-|\da^2z(\al)|^2$. Then operating similarly as we did
with $M_5$ and $M_6$, we get
$$
M_7+M_8\leq \!\exp C(\|\F(z)\|^2_{L^\infty}\!\!+\!\| z\|^2_{H^3})
-\frac{A_\mu}{2\pi A(t)}\!\int_\T\!BR(z,\varpi)(\al)\cdot\da^{\bot} z(\al)\,\da^3 z_2(\al)\la(\da^3 z_2)(\al) d\alpha.
$$ The addition of both inequalities produces
$$
L_5\leq\!\exp C(\|\F(z)\|^2_{L^\infty}\!\!+\!\| z\|^2_{H^3})
-\frac{A_\mu}{2\pi A(t)}\!\int_\T\!BR(z,\varpi)(\al)\cdot\da^{\bot} z(\al)\,\da^3 z(\al)\cdot\la(\da^3 z)(\al) d\alpha.
$$
and all the previous discussion shows  that $I_5$ satisfies
identical estimates than $L_5$.

%%%%%%%%%%%%%%%%%%%%%%%%%%%%%%%%%%%%%%%%%%%%%%%%%%%%%%%%%%%%%%%%%%%%%%%%%%%%%%%%%%%%%%%%%%%%%%%%%%%%%%%%%%%%%%%%%
%%%%%%%%%%%%%%%%%%%%%%%%%%%%%%%%%%%%%%%%%%%%%%%%%%%%%%%%%%%%%%%%%%%%%%%%%%%%%%%%%%%%%%%%%%%%%%%%%%%%%%%%%%%%%%%%%

\subsection{Estimates on $\da^3(c(\al,t)\da z(\al,t)).$}

%%%%%%%%%%%%%%%%%%%%%%%%%%%%%%%%%%%%%%%%%%%%%%%%%%%%%%%%%%%%%%%%%%%%%%%%%%%%%%%%%%%%%%%%%%%%%%%%%%%%%%%%%%%%%%%%%
%%%%%%%%%%%%%%%%%%%%%%%%%%%%%%%%%%%%%%%%%%%%%%%%%%%%%%%%%%%%%%%%%%%%%%%%%%%%%%%%%%%%%%%%%%%%%%%%%%%%%%%%%%%%%%%%%

In the evolution of the $L^2$ norm of $\da^3 z(\al)$ it
remains to control  the term
$$I_2 = \int_\T \da^3 z(\al)\cdot \da^3(c(\al)\da z(\al))d\alpha.$$
Let us recall the formula
\begin{align}
\begin{split}\label{ladt}
c(\al,t)&=\frac{\al+\pi}{2\pi A(t)}\int_\T \partial_{\beta}
z(\beta,t)\cdot \partial_{\beta} BR(z,\varpi)(\beta,t) d\beta\\
&\quad-\frac{1}{A(t)}\int_{-\pi}^\al \partial_{\beta}
z(\beta,t)\cdot\partial_{\beta} BR(z,\varpi)(\beta,t) d\beta.
\end{split}
\end{align}
We take $I_2 = J_1 + J_2 + J_3 + J_4$, where
$$J_1=\int_\T \da^3 z(\al)\cdot \da^4 z(\al)\, c(\al) d\alpha,\quad J_2=3\int_\T |\da^3 z(\al)|^2 \da c(\al)d\alpha,$$
$$J_3=3\int_\T \da^3 z(\al)\cdot \da^2 z(\al)\, \da^2 c(\al) d\alpha,\quad J_4=\int_\T \da^3 z(\al)\cdot \da z(\al)\, \da^3 c(\al)d\alpha.$$
An integration by parts in $J_1$ yields
$$
J_1=-\frac12\int_\T |\da^3 z(\al)|^2 \da c(\al) d\alpha\leq \|\da c\|_{L^\infty}\|\da^3 z\|^2_{L^2}\leq 2\|\F(z)\|^{1/2}_{L^\infty}\|\da BR(z,\varpi)\|_{L^\infty}\|\da^3 z\|^2_{L^2},
$$
and the estimate for $\|\da BR(z,\varpi)\|_{L^\infty}$ obtained before
gives us
\begin{equation*}
J_1\leq \exp C(\|\F(z)\|^2_{L^\infty}(t)+\|z\|^2_{H^3}(t)).
\end{equation*}
The term $J_2$ satisfies that $J_2=-6 J_1$, therefore
\begin{equation*}
J_2\leq \exp C(\|\F(z)\|^2_{L^\infty}(t)+\|z\|^2_{H^3}(t)).
\end{equation*}
Next we split $J_3=K_1+K_2,$ where
$$
K_1=-3\int_\T \da^3 z(\al)\cdot \da^2 z(\al)\frac{\da^2 z(\al)}{|\da z(\al)|^2}\cdot \da BR(z,\varpi) d\alpha,
$$
$$
K_2=-3\int_\T \da^3 z(\al)\cdot \da^2 z(\al)\frac{\da z(\al)}{|\da z(\al)|^2}\cdot \da^2 BR(z,\varpi) d\alpha.
$$
We have $$K_1\leq 3\|\F(z)\|_{L^\infty}\|z\|^2_{C^2}\|\da
BR(z,\varpi)\|_{L^\infty}\|\da^3 z\|_{L^2}\leq \exp
C(\|\F(z)\|^2_{L^\infty}(t)+\|z\|^2_{H^3}(t)),$$
$$K_2\leq \|\F(z)\|^{1/2}_{L^\infty}\|z\|_{C^2}\|\da^3 z\|_{L^2}\|\da^2 BR(z,\varpi)\|_{L^2},$$
and therefore
\begin{equation*}
I_3\leq \exp C(\|\F(z)\|^2_{L^\infty}(t)+\|z\|^2_{H^3}(t)).
\end{equation*}
The equality $\da^3 z(\al)\cdot \da z(\al)=-|\da^2 z(\al)|^2$
yields
$$
I_4=-\int_{\T}|\da^2 z(\al)|^2\da^3 c(\alpha)d\al=2\int_{\T}\da^3 z(\al)\cdot \da^2 z(\al) \da^2 c(\alpha)d\al=\frac23 I_3,
$$
and finally
\begin{equation*}
I_4\leq \exp C(\|\F(z)\|^2_{L^\infty}(t)+\|z\|^2_{H^3}(t)).
\end{equation*}

%%%%%%%%%%%%%%%%%%%%%%%%%%%%%%%%%%%%%%%%%%%%%%%%%%%%%%%%%%%%%%%%%%%%%%%%%%%%%%
%%%%%%%%%%%%%%%%%%%%%%%%%%%%%%%%%%%%%%%%%%%%%%%%%%%%%%%%%%%%%%%%%%%%%%%%%%%%%%

\section{The arc-chord condition}\label{arc-chord}

%%%%%%%%%%%%%%%%%%%%%%%%%%%%%%%%%%%%%%%%%%%%%%%%%%%%%%%%%%%%%%%%%%%%%%%%%%%%%%
%%%%%%%%%%%%%%%%%%%%%%%%%%%%%%%%%%%%%%%%%%%%%%%%%%%%%%%%%%%%%%%%%%%%%%%%%%%%%%

In this section we analyze the evolution of the quantity
$\|\F(z)\|_{L^\infty}(t)$, which gives the local control of the
arc-chord condition.

\begin{lemma}\label{lemaarcchord} The following estimate holds
\begin{align}
\begin{split}\label{enlif}
\D\dt\|\F(z)\|^2_{L^\infty}(t)&\leq \exp C(\|\F(z)\|^2_{L^\infty}(t)+\|z\|^2_{H^3}(t))
\end{split}
\end{align}
\end{lemma}

Proof: Let take $p>1$. It follows that
\begin{align*}
\D\dt\|\F(z)\|^p_{L^p}(t)&=\dt\int_{\T}\int_\T\big(\frac{|\beta|^2/4}{|V(\al,\al-\beta,t)|^2}\big)^{p}d\beta d\al\\
&=I_1+I_2+I_3,
\end{align*}
where
$$
I_1=-p\int_{\T}\int_\T(|\beta|/2)^{2p}\frac{V(\al,\al-\beta,t)\cdot (z_t(\al,t)-z_t(\al-\beta,t))}{|V(\al,\al-\beta,t)|^{2p+2}}d\beta d\al,
$$
$$
I_2=-p\int_{\T}\int_\T(|\beta|/2)^{2p}\frac{V_1^3(\al,\al-\beta,t)(z_{1\,t}(\al,t)-z_{1\,t}(\al-\beta,t))}{|V(\al,\al-\beta,t)|^{2p+2}}
 d\beta d\al,
$$
and
$$
I_3=p\int_{\T}\int_\T(|\beta|/2)^{2p}\frac{V_2^3(\al,\al-\beta,t)(z_{2\,t}(\al,t)-z_{2\,t}(\al-\beta,t))}{|V(\al,\al-\beta,t)|^{2p+2}} d\beta d\al.
$$
For $I_1$ we have

\begin{align*}
I_1&\leq p\int_{\T}\int_\T\big(\frac{|\beta|/2}{|V(\al,\al-\beta,t)|}\big)^{2p+1}\frac{|z_t(\al,t)-z_t(\al-\beta,t)|}{|\beta|}d\beta d\al.
\end{align*}
Let us consider
\begin{align*}
z_t(\al)-z_t(\al\!-\!\beta)&=(BR(z,\varpi)(\al)-BR(z,\varpi)(\al-\beta))+(c(\al)-c(\al-\beta))\da z(\al)\\
&\quad+c(\al-\beta)(\da z(\al)-\da z(\al-\beta))\\
&=J_1+J_2+J_3.
\end{align*}
Then  for $J_1$ we get $|J_1|\leq \|\da
BR(z,\varpi)\|_{L^{\infty}}|\beta|$, and the estimate
\eqref{nliibr} gives
$$|J_1|\leq \exp C(\|\F(z)\|^2_{L^\infty}(t)+\|z\|^2_{H^3}(t))|\beta|.$$
Using the definition for $c(\al)$ easily we obtain that
$$
|c(\al)-c(\al-\beta)|\leq \frac{\|\da BR(z,\varpi)\|_{L^{\infty}}}{A(t)^{1/2}}|\beta|,
$$
and again using \eqref{nliibr} we get
$$|J_2|\leq \exp C(\|\F(z)\|^2_{L^\infty}(t)+\|z\|^2_{H^3}(t))|\beta|.$$
For $J_3$ we have $|J_3|\leq \|c\|_{L^\infty}\|z\|_{C^2}|\beta|$
that is
$$|J_3|\leq
\exp C(\|\F(z)\|^2_{L^\infty}(t)+\|z\|^2_{H^3}(t))|\beta|.$$ Those
estimates obtained for  $J_i$ allow us to write
\begin{align*}
I_1&\leq p \exp C(\|\F(z)\|^2_{L^\infty}(t)+\|z\|^2_{H^3}(t))\int_{\T}\int_\T\big(\frac{|\beta|/2}{|V(\al,\al-\beta,t)|}\big)^{2p+1}d\beta d\al\\
&\leq p
\|\F(z)\|^{1/2}_{L^\infty}(t)\big(\exp C(\|\F(z)\|^2_{L^\infty}(t)+\|z\|^2_{H^3}(t))\big)\|\F(z)\|^p_{L^p}(t).
\end{align*}
 H\"older inequality implies
\begin{align*}
|I_2|&\leq p C \|z_{1\,t}\|_{L^\infty}\int_\T\int_\T\big(
\frac{|\beta|/2}{|V(\al,\al-\beta,t)|}\big)^{2p-1}d\beta d\al
\leq p C \|z_{1\,t}\|_{L^\infty} (1+\|\F(z)\|^p_{L^p})\\
&\leq p \big(\exp C(\|\F(z)\|^2_{L^\infty}(t)+\|z\|^2_{H^3}(t))\big)\|\F(z)\|^p_{L^p}.
\end{align*}
Since $|V_2(\al,\al-\beta)|\leq 1$ we have
\begin{align*}
|I_3|&\leq 2p \|z_{2\,t}\|_{L^\infty}\int_\T\int_\T\big(
\frac{|\beta|/2}{|V(\al,\al-\beta,t)|}\big)^{2p}d\beta d\al\\
&\leq 2p \big(\exp C(\|\F(z)\|^2_{L^\infty}(t)+\|z\|^2_{H^3}(t))\big)\|\F(z)\|^p_{L^p}.
\end{align*}
Those estimates for $I_1$,$I_2$ and $I_3$ gives us
\begin{align*}
\D\dt\|\F(z)\|^p_{L^p}(t)&\leq p \big(\exp C(\|\F(z)\|^2_{L^\infty}(t)+\|z\|^2_{H^3}(t))\big)\|\F(z)\|^p_{L^p}(t),
\end{align*}
 therefore
\begin{align*}
\D\dt\|\F(z)\|_{L^p}(t)&\leq \big(\exp C(\|\F(z)\|^2_{L^\infty}(t)+\|z\|^2_{H^3}(t))\big)\|\F(z)\|_{L^p}(t).
\end{align*}
After an integration in the time variable $t$ we get
\begin{align*}
\|\F(z)\|_{L^p}(t+h)&\leq \|\F(z)\|_{L^p}(t) \exp\,\big(
\!\int_t^{t+h}e^{C(\|\F(z)\|^2_{L^\infty}(s)+\|z\|^2_{H^3}(s))}ds\big),
\end{align*}
and letting $p\rightarrow \infty$ we obtain
\begin{align*}
\|\F(z)\|_{L^\infty}(t+h)&\leq \|\F(z)\|_{L^\infty}(t) \exp\,\big(
\!\int_t^{t+h}e^{C(\|\F(z)\|^2_{L^\infty}(s)+\|z\|^2_{H^3}(s))}ds\big).
\end{align*}
Therefore
\begin{align*}
\D\dt\|\F(z)\|_{L^\infty}(t)&=\lim_{h\rightarrow
0}(\|\F(z)\|_{L^\infty}(t+h)-\|\F(z)\|_{L^\infty}(t))h^{-1}\\
&\leq \|\F(z)\|_{L^\infty}(t)\lim_{h\rightarrow 0}( \exp\,\big(
\!\int_t^{t+h}e^{C(\|\F(z)\|^2_{L^\infty}(s)+\|z\|^2_{H^3}(s))}ds\big)-1)h^{-1}\\
&\leq \|\F(z)\|_{L^\infty}(t)e^{C(\|\F(z)\|^2_{L^\infty}(t)+\|z\|^2_{H^3}(t))}.
\end{align*}
With this we finish the proof of lemma \ref{lemaarcchord}. q.e.d.

%%%%%%%%%%%%%%%%%%%%%%%%%%%%%%%%%%%%%%%%%%%%%%%%%%%%%%%%%%%%%%%%%%%%%%%%%%%%%%%%%%%%%%%%%%%%%%%%%%%%%%%%%%%%%%%%%%%
%%%%%%%%%%%%%%%%%%%%%%%%%%%%%%%%%%%%%%%%%%%%%%%%%%%%%%%%%%%%%%%%%%%%%%%%%%%%%%%%%%%%%%%%%%%%%%%%%%%%%%%%%%%%%%%%%%%

\section{The evolution of the minimum of $\sigma(\al,t)$}\label{minimumsigma}
%the difference of the gradients of the pressure in the normal direction

%%%%%%%%%%%%%%%%%%%%%%%%%%%%%%%%%%%%%%%%%%%%%%%%%%%%%%%%%%%%%%%%%%%%%%%%%%%%%%%%%%%%%%%%%%%%%%%%%%%%%%%%%%%%%%%%%%%
%%%%%%%%%%%%%%%%%%%%%%%%%%%%%%%%%%%%%%%%%%%%%%%%%%%%%%%%%%%%%%%%%%%%%%%%%%%%%%%%%%%%%%%%%%%%%%%%%%%%%%%%%%%%%%%%%%%

In this section we get an a priori estimate for the evolution of
the minimum of the difference of the gradients of the pressure in
the normal direction to the interface. This quantity is given by
\begin{equation}\label{fsigma}
\sigma(\al,t)=\frac{\mu^2-\mu^1}{\kappa}BR(z,\varpi)(\al,t)\cdot\dpa
z(\al,t)+ \mathrm{g}(\rho^2-\rho^1)\da z_1(\al,t).
\end{equation}
\begin{lemma}
Let $z(\al,t)$ be a solution of the system with $z(\al,t)\in C^{1}([0,T];H^3)$, and $$m(t)=\D\min_{\al\in\T} \sigma(\al,t).$$
Then $$m(t)\geq m(0)-\int_0^t \exp C(\|\F(z)\|^2_{L^\infty}(s)+\|z\|^2_{H^3}(s)) ds.$$
\end{lemma}

Proof: Suppose that $z(\al,t)\in C^1([0,T];H^3)$ is a solution of
the system, then the result obtained in the preceding sections
together with  Sobolev inequalities show that $\sigma(\al,t)\in
C^1([0,T]\times \T).$ Therefore we may consider $\al_t\in\T$ such
that
$$m(t)=\D\min_{\al\in\T} \sigma(\al,t)=\sigma(\al_t,t),$$
which is a Lipschitz function  differentiable almost everywhere.
With an analogous argument to the one used in  \cite{cor2} and
\cite{DY2}, we may calculate the derivative of $m(t)$, to obtain
$$
m'(\al_t,t)=\sigma_t(\al_t,t).
$$
The identity \eqref{fsigma} yields
\begin{align*}
\sigma_t(\al,t)&=\frac{\mu^2-\mu^1}{\kappa}\partial_t (BR(z,\varpi))(\al,t)\cdot\dpa
z(\al,t)\\
&\quad+(\frac{\mu^2-\mu^1}{k}\, BR(z,\varpi)(\al)\cdot\dpa
z_t(\al,t)+\mathrm{g}(\rho^2-\rho^1)\da z_{1\, t}(\al,t))\\
&=I_1+I_2.
\end{align*}
And we have
$$
|I_2|\leq C(\|BR(z,\varpi)\|_{L^\infty}+1)\|\da z_t\|_{L^\infty}.
$$
We can easily estimate $\|BR(z,\varpi)\|_{L^\infty}$, obtaining
$$
|I_2|\leq \exp C(\|\F(z)\|^2_{L^\infty}+\|z\|^2_{H^3})\|\da z_t\|_{L^\infty}.
$$
Next we use equation \eqref{fullpm} to get
\begin{align*}
\|\da z_t\|_{L^\infty}&\leq (\|\da BR(z,\varpi)\|_{L^\infty}+\|\da c\|_{L^\infty}\|\da z\|_{L^\infty}+\|c\|_{L^\infty}\|\da^2 z\|_{L^\infty})\\
&\leq C\|\da BR(z,\varpi)\|_{L^\infty}(1+\|\F(z)\|^{1/2}_{L^\infty}\|z\|_{C^2}),
\end{align*} and with the bound obtained before for $\|\da BR(z,\varpi)\|_{L^\infty}$  \eqref{nliibr}, we
have
$$
|I_2|\leq \exp C(\|\F(z)\|^2_{L^\infty}+\|z\|^2_{H^3}).
$$
Let us write
$BR(z,\varpi)(\al,t)=P_1(\al,t)+P_2(\al,t)+P_3(\al,t)$, where the
$P_j$ were defined in \eqref{p1p2p3}.  We have
$$
|\dpt P_2(\al)|+|\dpt P_3(\al)|\leq C(\|\varpi_t\|_{L^2}+\|\varpi\|_{L^2}\|z_t\|_{L^\infty}).
$$
The norm $\|z_t\|_{L^\infty}$ is bounded by \eqref{fullpm}, and
the adequate estimates for $\|\varpi_t\|_{L^2}$ which will be
introduced later. In $\dpt P_1$ there are terms of lower order
which can be estimated as $\dpt P_2$ and $\dpt P_3$, but the most
singular ones are given by
$$
J_1=\frac{1}{4\pi}PV\int_{\T}\varpi_t(\al-\beta)\frac{V^\bot(\al,\al-\beta)}{|V(\al,\al-\beta)|^2}d\beta,
$$
$$
J_2=\frac{1}{8\pi}PV\int_{\T}\varpi(\al-\beta)\frac{z_t(\al)-z_t(\al-\beta)}{|V(\al,\al-\beta)|^2}d\beta,
$$
$$
J_3=-\frac{1}{4\pi}PV\int_{\T}\varpi(\al-\beta)\frac{V^\bot(\al,\al-\beta)}{|V(\al,\al-\beta)|^4}\big(V(\al,\al-\beta)\cdot(z_t(\al)-z_t(\al-\beta))\big)d\beta,
$$
Let us now  split $J_1$ in a similar way as we did before,  to
obtain
$$
J_1=\frac{1}{4\pi}PV\!\!\int_{\T}\!\varpi_t(\al\!-\!\beta)\big(\frac{V^\bot(\al,\al\!-\!\beta)}{|V(\al,\al\!-\!\beta)|^2}-\frac{\da^{\bot} z(\al)}{|\da z(\al)|^2\tan(\beta/2)}\big)d\beta+\frac{\da^{\bot} z(\al)}{2|\da z(\al)|^2}H(\varpi_t)(\al),
$$
and
$$
|J_1|\leq C\|\F(z)\|^k_{L^\infty}\|z\|^k_{C^2}\|\varpi_t\|_{L^2}+\|\F(z)\|^{1/2}_{L^\infty}\|\varpi_t\|_{C^\delta}.
$$
Next we divide $J_2=K_1+K_2+K_3$ where
$$
K_1=\frac{1}{8\pi}\int_{\T}(\varpi(\al-\beta)-\varpi(\al))\frac{z_t(\al)-z_t(\al-\beta)}{|V(\al,\al-\beta)|^2}d\beta,
$$
$$
K_2=\frac{\varpi(\al)}{8\pi}\int_{\T}(z_t(\al)-z_t(\al-\beta))\big(\frac{1}{|V(\al,\al-\beta)|^2}-\frac{4}{|\da z(\al)|^2|\beta|^2}\big)d\beta,
$$
$$
K_3=\frac{1}{8\pi}\frac{\varpi(\al)}{|\da z(\al)|^2}\int_{\T}(z_t(\al)-z_t(\al-\beta))\big(\frac{4}{|\beta|^2}-\frac{1}{\sin^2(\beta/2)}\big)d\beta+\frac{1}{2}\frac{\varpi(\al)}{|\da z(\al)|^2}\la (z_t)(\al).
$$
The identity $$z_t(\al)-z_t(\al-\beta)=\beta\int_0^1\da
z_t(\al+(s-1)\beta) ds,$$ gives us
$$|K_1|+|K_2|\leq \|\varpi\|_{C^1}\|\F(z)\|^{k}\|z\|^k_{C^2}\|\da z_t\|_{L^\infty}.$$
And for $K_3$ we have
$$|K_3|\leq C\|\varpi\|_{L^\infty}\|\F(z)\|_{L^\infty}\|z_t\|_{C^{1,\delta}}.$$

In order to control $\|\varpi_t\|_{C^\delta}$ we will use the inequality
$$
\|f\|_{C^\delta}\leq C(\|f\|_{L^2}+\|f\|_{\overline{C}^{\delta}}),
$$
with
$$
\|f\|_{\overline{C}^{\delta}}=\sup_{\al\neq \beta}\frac{|f(\al)-f(\beta)|}{|\al-\beta|^{\delta}}.
$$
Let us now take the time derivative of the identity
\eqref{formulavarpi}, we get
$$
\varpi_t(\al)+A_\mu T(\varpi_t)(\al)=-\frac{A_\mu}{2\pi}
R(\al)-2\kappa g\frac{\rho^2-\rho^1}{\mu^2+\mu^1}\da z_{2\,t}(\al)
$$
which yields
$$
\|\varpi_t\|_{H^{\frac12}}\leq C\|(I+A_\mu
T)^{-1}\|_{H^{\frac12}}(\|R\|_{H^{\frac12}}+\|\da z_{t}\|_{H^{\frac12}})
$$
and since we control $\|(I+A_\mu T)^{-1}\|_{H^{\frac12}}$ it remains
to estimate $\|R\|_{H^{\frac12}}$.

Instead we will estimate $\|R\|_{H^{1}},$ and to do that we consider
the splitting $R= S_1 + S_2 + S_3$ where

$$
S_1(\al)=\int_{\T}\varpi(\al-\beta)\dpt \Big(\frac{(V(\al,\al-\beta)-\da z(\al)\beta/2)^{\bot}\cdot \da z(\al)}{|V(\al,\al-\beta)|^2}\Big) d\beta,
$$
$$
S_2(\al)=-\int_{\T}\varpi(\al-\beta)\dpt\Big(V^2
_2(\al,\al-\beta)\frac{V^{\bot}(\al,\al-\beta)\cdot \da z(\al)}{|V(\al,\al-\beta)|^2}\Big) d\beta,
$$
$$
S_3(\al)=-\int_{\T}\varpi(\al-\beta)\dpt\Big(V_2(\al,\al-\beta)\da z_1(\al)\Big)d\beta.
$$
The terms $S_2(\al)$ and $S_3(\al)$ are controlled as follows:
$$
|S_2(\al)|+|S_3(\al)|\leq C\|z_t\|_{C^1}\|\varpi\|_{L^2}.
$$
For $S_1$ we split $S_1(\al)=Q_1(\al)+Q_2(\al)+Q_3(\al)+Q_4(\al)+Q_5(\al)$, where
$$
Q_1(\al)=\int_{\T}\varpi(\al-\beta)\frac{(V(\al,\al-\beta)-\da z(\al)\beta/2)^{\bot}\cdot \da z_t(\al)}{|V(\al,\al-\beta)|^2} d\beta,
$$
$$
Q_2(\al)=\frac{1}{2}\int_{\T}\varpi(\al-\beta)\frac{(z_t(\al)-z_t(\al-\beta)-\da z_t(\al)\beta)^{\bot}\cdot \da z(\al)}{|V(\al,\al-\beta)|^2} d\beta,
$$
$$
Q_3(\al)=\frac{1}{2}\int_{\T}\varpi(\al-\beta)\frac{(V(\al,\al\!-\!\beta)\!-\!\da z(\al)\beta/2)^{\bot}\cdot \da z(\al)}{|V(\al,\al-\beta)|^4}B(\al,\al-\beta) d\beta,
$$
with
\begin{equation}\label{fb}
B(\al,\al-\beta)= V(\al,\al\!-\!\beta)\cdot (z_t(\al)-z_t(\al\!-\!\beta)),
\end{equation}
$$
Q_4(\al)=\frac{1}{2}\int_{\T}\varpi(\al-\beta)\frac{C(\al,\al-\beta)}{|V(\al,\al-\beta)|^2} d\beta,
$$
\begin{align*}
C(\al,\al-\beta)&=V_1^2(\al,\al-\beta)(z_{1\, t}(\al)- z_{1\, t}(\al-\beta))\da z_2(\al)\\
&\quad+V_2^2(\al,\al-\beta)( z_{2\, t}(\al)-z_{2\, t}(\al-\beta))\da z_1(\al)
\end{align*}
$$
Q_5(\al)=-\int_{\T}\varpi(\al-\beta)\frac{(V(\al,\al\!-\!\beta)\!-\!\da z(\al)\beta/2)^{\bot}\cdot \da z(\al)}{|V(\al,\al-\beta)|^4}D(\al,\al-\beta) d\beta,
$$
$$
D(\al,\al-\beta)=V_1^3(\al,\al-\beta)(z_{1\, t}(\al)- z_{1\, t}(\al-\beta))-V_2^3(\al,\al-\beta)(z_{2\,t}(\al)-z_{2\,t}(\al-\beta)).
$$
We have $|Q_4(\al)|\leq
C\|z\|_{C^1}\|\varpi\|_{L^2}\|z_t\|_{L^\infty}$. In a similar way
this estimates follows for $Q_5$: $$|Q_5(\al)|\leq
C(\|z\|_{C^1}+\|\F(z)\|^{1/2}_{L^\infty}\|z\|^2_{C^1})\|\varpi\|_{L^2}\|z_t\|_{L^\infty}.$$
For $Q_1$ we proceed as before to obtain
$$|Q_1(\al)|\leq C \|\F(z)\|_{L^\infty}\|z\|_{C^2}\|\varpi\|_{L^2}\|z_t\|_{C^1}.$$
The inequality $$|z_t(\al)- z_t(\al-\beta)-\da z_t(\al)\beta|\leq
\|z_t\|_{C^{1,\delta}}|\beta|^{1+\delta}$$ gives
$$|Q_2(\al)|\leq C \|\F(z)\|_{L^\infty}\|z\|_{C^1}\|\varpi\|_{L^{\infty}}\|z_t\|_{C^{1,\delta}}.$$
And
$$|z_t(\al)-z_t(\al-\beta)|\leq \|z_t\|_{C^1}|\beta|,$$ yields
$$|Q_3(\al)|\leq \|\F(z)\|^{3/2}_{L^\infty}\|z\|^2_{C^2}\|\varpi\|_{L^2}\|z_t\|_{C^1}.$$
Finally we have
$$
|R_1(\al)|\leq \|\F(z)\|^{3/2}_{L^\infty}\|z\|^2_{C^2}\|\varpi\|_{H^1}\|z_t\|_{C^{1,\delta}}.
$$
Using \eqref{formulavarpi} we obtain
$$
\|\varpi_t\|_{\overline{C}^\delta}\leq C(\|T(\varpi_t)\|_{\overline{C}^\delta}+\|R\|_{\overline{C}^\delta}+\|\da z_t\|_{\overline{C}^\delta}).
$$
For $\delta\leq 1/2$ we have
$$
\|T(\varpi_t)\|_{\overline{C}^\delta}\leq\|T(\varpi_t)\|_{H^1} \leq 2\|\da T(\varpi_t)\|_{L^2}\leq \|\F(z)\|^2_{L^\infty}\|z\|^4_{C^{2,\delta}}\|w_t\|_{L^2}.
$$
Now to estimate $\|R\|_{\overline{C}^\delta}\leq \|R\|_{H^1}$  we
consider $\|\da R\|_{L^2}$. The most singular terms for this
quantity are those with two derivatives in $\alpha$ and one in
time, or with one derivative in $\alpha$, one in time and a
principal value. Let us write:

$$
Q_6(\al)=\int_{\T}\varpi(\al-\beta)\frac{(V(\al,\al-\beta)-\da z(\al)\beta/2)^{\bot}\cdot \da^2 z_t(\al)}{|V(\al,\al-\beta)|^2} d\beta,
$$
$$
Q_7(\al)=\frac{1}{2}\int_{\T}\varpi(\al-\beta)\frac{(\da z_t(\al)-\da z_t(\al-\beta)-\da^2 z_t(\al)\beta)^{\bot}\cdot \da z(\al)}{|V(\al,\al-\beta)|^2} d\beta,
$$
$$
Q_8(\al)=\frac{1}{2}\int_{\T}\varpi(\al-\beta)\frac{(V(\al,\al\!-\!\beta)\!-\!\da z(\al)\beta/2)^{\bot}\cdot \da z(\al)}{|V(\al,\al-\beta)|^4}D(\al,\al-\beta) d\beta,
$$
with
\begin{equation}\label{fb}
D(\al,\al-\beta)= V(\al,\al\!-\!\beta)\cdot (\da z_t(\al)-\da z_t(\al\!-\!\beta)).
\end{equation}
We have
$$
|Q_6(\al)|\leq \|\F(z)\|^k_{L^\infty}\|z\|^k_{C^2}\|w\|_{L^2}|\da^2 z_t(\al)|,
$$
and
$$
|Q_8(\al)|\leq \|\F(z)\|^k_{L^\infty}\|z\|^k_{C^2}\|w\|_{L^\infty}\| z_t(\al)\|_{C^{1,\delta}}.
$$
Let us split $Q_7(\al)=J_4+J_5$ where
$$
J_4=\frac{1}{2}PV\int_{\T}\varpi(\al-\beta)\frac{(\da z_t(\al)-\da z_t(\al-\beta))^{\bot}\cdot \da z(\al)}{|V(\al,\al-\beta)|^2} d\beta,
$$
and
$$
J_5=-\frac{1}{2}(\da^2 z_t(\al))^{\bot}\cdot \da z(\al)PV\int_{\T}\varpi(\al-\beta)\frac{\beta}{|V(\al,\al-\beta)|^2} d\beta.
$$
For $J_5$ we have  $|J_5|\leq
\|\F(z)\|^k_{L^\infty}\|z\|^k_{C^2}\|w\|_{H^1}|\da^2 z_t(\al)|$.
Next we divide $J_4=K_4+K_5+K_6+K_7$ where
$$
K_4=\frac{1}{2}\int_{\T}(\varpi(\al-\beta)-\varpi(\al))\frac{(\da z_t(\al)-\da z_t(\al-\beta))^{\bot}\cdot \da z(\al)}{|V(\al,\al-\beta)|^2} d\beta,
$$
$$
K_5=\frac{\varpi(\al)}{2}\int_{\T}(\da z_t(\al)-\da z_t(\al-\beta))^{\bot}\cdot \da z(\al)\big(\frac{1}{|V(\al,\al-\beta)|^2}-\frac{4}{|\da z(\al)|^2|\beta|^2}\big) d\beta,
$$
$$
K_6=\frac{2\varpi(\al)}{|\da z(\al)|^2}\int_{\T}(\da z_t(\al)-\da z_t(\al-\beta))^{\bot}\cdot \da z(\al)\big(\frac{1}{|\beta|^2}-\frac{1}{4\sin^2(\beta/2)}\big) d\beta,
$$
$$
K_7=2\pi\frac{\varpi(\al)}{|\da z(\al)|^2}(\la(\da z_t)(\al))^\bot\cdot \da z(\al).
$$
We have $|K_4|+|K_5|\leq C
\|\F(z)\|^k_{L^\infty}\|z\|^k_{C^2}\|w\|_{H^1}\|
z_t\|_{C^{1,\delta}}$, $|K_6|\leq
\|\F(z)\|_{L^\infty}\|z\|_{C^2}\|w\|_{H^1}\|z_t\|_{C^1}$ and
$|K_7|\leq \|\F(z)\|_{L^\infty}\|z\|_{C^2}\|w\|_{H^1}|\la(\da
z_t)(\al)|$.

Finally let us observe $\|z_t\|_{C^{1,\delta}}\leq \|z_t\|_{H^2}$,
which provide us the control of $\|\da^2 z_t\|_{L^2}$. We consider
now the terms of $\da^2 z_t(\al)$ given by
$$I_3=\da^2 BR(z,\varpi)(\al),\quad I_4=\da^2(c(\al)\da z(\al)).$$
Easily we get
$$
|I_4|\leq \|\F(z)\|^k_{L^\infty}\|z\|^k_{C^2}(1+|\da^3 z(\al)|+|\da^2 BR(z,\varpi)(\al)|),
$$
which yields
$$
\|I_4\|_{L^2}\leq (\|\F(z)\|^k_{L^\infty}\|z\|^k_{C^2}(1+\| z\|_{H^3}+\|\da^2 BR(z,\varpi)\|_{L^2}),
$$
so that we can control $\|\da^2 BR(z,\varpi)\|_{L^2}$ as in
\eqref{enl22dbr}, and finish the  estimate of $I_3$.

The upper bound $$|\sigma_t(\al,t)|\leq \exp
C(\|\F(z)\|^2_{L^\infty}(t)+\|z\|^2_{H^3}(t)),$$ gives us
$$
m'(t)\geq -\exp C(\|\F(z)\|^2_{L^\infty}(t)+\|z\|^2_{H^3}(t)),
$$
for almost every $t$. And a further integration yields lemma 8.1
$$
m(t)\geq m(0) - \int_0^t \exp(C|||z|||^2)ds.
$$

%%%%%%%%%%%%%%%%%%%%%%%%%%%%%%%%%%%%%%%%%%%%%%%%%%%%%%%%%%%%%%%%%%%%%%%%%%%%%
%%%%%%%%%%%%%%%%%%%%%%%%%%%%%%%%%%%%%%%%%%%%%%%%%%%%%%%%%%%%%%%%%%%%%%%%%%%%%

\section{Regularization and approximation}

Our next step is to use the a priori estimates to get
local-existence. For that purpose we introduce a regularized
 evolution equation having local-existence independently of the sign condition on $\sigma(\al,t)$ at $t=0$.
 But for
$\sigma(\al,0)>0$, we find a time of existence for the Muskat
problem uniformly in the regularization, allowing us to take the
limit.

Let $z^\ep (\al,t)$ be a solution of the following system:

\begin{align*}
\begin{split}
z^{\ep,\delta}_t(\al,t)&=BR^{\delta}(z^{\ep,\delta},\varpi^{\ep,\delta})(\al,t)+c^{\ep,\delta}(\al,t)\da z^{\ep,\delta}(\al,t),\\
z^{\ep,\delta}(\al,0)&=z_0(\al),
\end{split}
\end{align*}
where
\begin{align*}
\begin{split}
BR^\delta(z,\varpi)(\al,t)=\big(&-\!\frac{1}{4\pi}
\int_{\T}\varpi(\beta,t)\frac{\tanh
(\frac{z_2(\al,t)-z_2(\beta,t)}{2})
(1+\tan^2(\frac{z_1(\al,t)-z_1(\beta,t)}{2}))}{\tan^2(\frac{z_1(\al,t)-z_1(\beta,t)}{2})+\tanh^2(\frac{z_2(\al,t)-z_2(\beta,t)}{2})+\delta}d\beta,\\
&\frac{1}{4\pi}
\int_{\T}\varpi(\beta,t)\frac{\tan(\frac{z_1(\al,t)-z_1(\beta,t)}{2})(1-\tanh^2(\frac{z_2(\al,t)-z_2(\beta,t)}{2}))}{\tan^2(\frac{z_1(\al,t)-z_1(\beta,t)}{2})+\tanh^2(\frac{z_2(\al,t)-z_2(\beta,t)}{2})+\delta}d\beta\big),
\end{split}
\end{align*}
\begin{align*}
\begin{split}
\varpi^{\ep,\delta}(\al,t)=-A_\mu\phi_\ep\ast \phi_\ep\ast (2BR(z^{\ep,\delta},\varpi^{\ep,\delta})\cdot \da z^{\ep,\delta})(\al)-2\kappa\mathrm{g}\frac{\rho^2-\rho^1}{\mu^2+\mu^1}\phi_\ep\ast \phi_\ep\ast(\da
z^{\ep,\delta}_2)(\al),
\end{split}
\end{align*}
%\begin{align*}
%\begin{split}
%BR^\delta(z{\ep,\delta},\varpi^{\ep,\delta})(\al,t)=\big(&-\!\frac{1}{4\pi}
%\int_{\T}\varpi^{\ep,\delta}(\beta,t)\frac{\tanh
%(\frac{z^{\ep,\delta}_2(\al,t)-z^{\ep,\delta}_2(\beta,t)}{2})
%(1+\tan^2(\frac{z^{\ep,\delta}_1(\al,t)-z^{\ep,\delta}_1(\beta,t)}{2}))}{\tan^2(\frac{z^{\ep,\delta}_1(\al,t)-z^{\ep,\delta}_1(\beta,t)}{2})+\tanh^2(\frac{z^{\ep,\delta}_2(\al,t)-z^{\ep,\delta}_2(\beta,t)}{2})+\delta}d\beta,\\
%&\frac{1}{4\pi}
%\int_{\T}\varpi^{\ep,\delta}(\beta,t)\frac{\tan(\frac{z^{\ep,\delta}_1(\al,t)5-z^{\ep,\delta}_1(\beta,t)}{2})(1-\tanh^2(\frac{z^{\ep,\delta}_2(\al,t)-z^{\ep,\delta}_2(\beta,t)}{2}))}{\tan^2(\frac{z^{\ep,\delta}_1(\al,t)-z^{\ep,\delta}_1(\beta,t)}{2})+\tanh^2(\frac{z^{\ep,\delta}_2(\al,t)-z^{\ep,\delta}_2(\beta,t)}{2})+\delta}d\beta\big),
%\end{split}
%\end{align*}

\begin{align*}
\begin{split}
c^{\ep,\delta}(\al,t)&=\frac{\al+\pi}{2\pi}\int_\T\frac{\da z^{\ep,\delta}(\al,t)}{|\da
z^{\ep,\delta}(\al,t)|^2}\cdot \da BR^{\delta}(z^{\ep,\delta},\varpi^{\ep,\delta})(\al,t) d\al\\
&\quad-\int_{-\pi}^\al
\frac{\da z^{\ep,\delta}(\beta,t)}{|\da z^{\ep,\delta}(\al,t)|^2}\cdot \partial_{\beta}
BR^{\delta}(z^{\ep,\delta},\varpi^{\ep,\delta})(\beta,t) d\beta,
\end{split}
\end{align*}
$$
\phi\in C^\infty_c(\R),\quad \phi(\al)\geq 0,\quad \phi(-\al)=\phi(\al),\quad \int_\R \phi(\al)d\al=1,\quad \phi_\ep(\al)=\phi(\al/\ep)/\ep,
$$
for $\ep>0$ and $\delta>0$.

Then the operator $I+A_\mu\phi_\ep\ast\phi_\ep\ast T$ has a
bounded inverse in $H^{\frac12}$, for $\ep$ small enough,  with a
norm bounded independently of $\ep>0$. For this system there is
local-existence for initial data with $\F(z_0)(\al,\beta)< \infty$
even if $\sigma(\al,0)$ does not have the proper sign (see
\cite{Y}). So that there exists a time $T^{\ep,\delta}$ and a
solution of the system $z^{\ep,\delta}\in
C^1([0,T^{\ep,\delta}],H^k)$ for $k\leq 3$, and as long as the
solution  exists, we have $|\da z^{\ep,\delta}
(\al,t)|^2=A^{\ep,\delta}(t)$. Taking advantage of this property,
and using that $\varpi^{\ep,\delta}$ is regular, we obtain
estimates which are independent of $\delta$. Letting now
$\delta\rightarrow 0$ we get local-existence for the following
system:
\begin{align*}
\begin{split}
z^{\ep}_t(\al,t)&=BR(z^{\ep},\varpi^{\ep})(\al,t)+c^{\ep}(\al,t)\da z^{\ep}(\al,t),\\
z^{\ep}(\al,0)&=z_0(\al),
\end{split}
\end{align*}
where
\begin{align*}
\begin{split}
c^{\ep}(\al,t)&=\frac{\al+\pi}{2\pi}\int_\T\frac{\da z^{\ep}(\al,t)}{|\da
z^{\ep}(\al,t)|^2}\cdot \da BR(z^{\ep},\varpi^{\ep})(\al,t) d\al\\
&\quad-\int_{-\pi}^\al
\frac{\da z^{\ep}(\beta,t)}{|\da z^{\ep}(\al,t)|^2}\cdot \partial_{\beta}
BR(z^\ep,\varpi^\ep)(\beta,t) d\beta,
\end{split}
\end{align*}
\begin{align*}
\begin{split}
\varpi^{\ep}(\al,t)=-A_\mu\phi_\ep\ast \phi_\ep\ast (2BR(z^{\ep},\varpi^{\ep})\cdot \da z^{\ep})(\al)-2\kappa\mathrm{g}\frac{\rho^2-\rho^1}{\mu^2+\mu^1}\phi_\ep\ast \phi_\ep\ast(\da
z^{\ep}_2)(\al).
\end{split}
\end{align*}
%$$
%\phi\in C^\infty_c(\R),\quad \phi(\al)\geq 0,\quad \phi(-\al)=\phi(\al),\quad \int_\R \phi(\al)d\al=1,\quad \phi_\ep(\al)=\phi(\al/\ep)/\ep,
%$$
%for $\ep>0$.
Next we will show that for this system we have
\begin{align}
\begin{split}\label{ntniepsilon}
\frac{d}{dt}\|z^\ep\|^2_{H^k}(t)&\leq -\frac{\kappa}{2\pi(\mu_1\!+\!\mu_2)}\!\int_\T \frac{\sigma^\ep(\al,t)}{|\da z^\ep(\al,t)|^2}  \phi_\ep\ast (\da^k z^\ep)(\al,t)\cdot \la(\phi_\ep\ast(\da^k z^\ep))(\al,t) d\al\\
&\quad+\exp C(\|\F(z^\ep)\|^2_{L^\infty}(t)+\|z^\ep\|^2_{H^k}(t)).
\end{split}
\end{align}
where
$$
\sigma^\ep(\al,t)=\frac{\mu^2-\mu^1}{\kappa}BR(z^\ep,\varpi^\ep)(\al,t)\cdot\dpa
z^\ep(\al,t)+ \mathrm{g}(\rho^2-\rho^1)\da z^\ep_1(\al,t).
$$
To do the task we have to repeat the arguments in our previous
sections, with the exception of 7.3. (looking for
$\sigma^\ep(\al)$) where we   proceed  differently using the
following well-known estimate for the commutator of the
convolution:
\begin{equation}\label{conconvo}
\|\phi_\ep\ast(gf)-g\phi_\ep\ast(f)\|_{H^1}\leq C\|g\|_{C^1}\|f\|_{L^2}
\end{equation}
regarding, where the constant $C$ is independent of $\ep$.

In the following we will present the details of the evolution of
the $L^2$ norm of the third derivatives, being the case  of the
kth-derivative ($k>3$) completely analogous. Furthermore, with
regards of the different decompositions introduced  in the
previous sections, in the following we shall select only the more
singular terms, showing for them the corresponding uniform
estimates and leaving to the reader the remainder easy cases.

If we consider the term corresponding to $K_2$ in section 7.3 we have
$$
K^\ep_2=-\frac{\kappa\mathrm{g}(\rho^2\!-\!\rho^1)}{2\pi(\mu^2\!+\!\mu^1)A^\ep(t)}\int_\T\la\big(\da^3z^\ep\cdot\da^\bot z^\ep\big)(\al)\phi_\ep\ast\phi_\ep\ast (\da^3z^\ep_2)(\al)d\al.
$$
which, we write in the following manner
$$
K^\ep_2=-\frac{\kappa\mathrm{g}(\rho^2\!-\!\rho^1)}{2\pi(\mu^2\!+\!\mu^1)A^\ep(t)}\int_\T\la\big(\phi_\ep\ast(\da^3z^\ep\cdot\da^\bot z^\ep)\big)(\al)\phi_\ep\ast (\da^3z^\ep_2)(\al)d\al.
$$
Then we have $K^\ep_2=L^\ep_1+L^\ep_2$, where
\begin{align*}
L^\ep_1=\frac{\kappa\mathrm{g}(\rho^2-\rho^1)}{2\pi (\mu^2+\mu^1)A^\ep(t)}\int_\T\la\big(\phi_\ep\ast(\da^3 z^\ep_1
\da z^\ep_2)\big)(\al)\phi_\ep\ast(\da^3z^\ep_2)(\al) d\alpha,
\end{align*}
\begin{align*}
L^\ep_2=-\frac{\kappa\mathrm{g}(\rho^2-\rho^1)}{2\pi
(\mu^2+\mu^1)A^\ep(t)}\int_\T\la\big(\phi_\ep\ast(\da^3 z^\ep_2 \da
z^\ep_1)\big)(\al)\phi_\ep\ast (\da^3z^\ep_2)(\al) d\alpha.
\end{align*}
Next we write $L^\ep_1=M^\ep_1+M^\ep_2+M^\ep_3+M^\ep_4$ where
$$
M^\ep_1=\frac{\kappa\mathrm{g}(\rho^2\!-\!\rho^1)}{2\pi (\mu^2\!+\!\mu^1)A^\ep(t)}\int_\T\la\big(\phi_\ep\ast(\da^3 z^\ep_1
\da z^\ep_2)-\phi_\ep\ast(\da^3 z^\ep_1) \da z^\ep_2\big)(\al)\phi_\ep\ast(\da^3z^\ep_2)(\al) d\alpha,
$$
$$
M^\ep_2=\frac{\kappa\mathrm{g}(\rho^2\!-\!\rho^1)}{2\pi (\mu^2\!+\!\mu^1)A^\ep(t)}\int_\T
[\la\big(\phi_\ep\ast(\da^3 z^\ep_1)\da z^\ep_2\big)(\al)-\la\big(\phi_\ep\ast(\da^3 z^\ep_1)\big)(\al)\da z^\ep_2(\al)]\phi_\ep\ast(\da^3z^\ep_2)(\al) d\alpha,
$$
$$
M^\ep_3=\frac{\kappa\mathrm{g}(\rho^2\!-\!\rho^1)}{2\pi (\mu^2\!+\!\mu^1)A^\ep(t)}\int_\T
\la\big(\phi_\ep\ast(\da^3 z^\ep_1)\big)(\al)[\da z^\ep_2(\al)\phi_\ep\ast(\da^3z^\ep_2)(\al)-\phi_\ep\ast(\da z^\ep_2\da^3z^\ep_2)(\al)] d\alpha,
$$
$$
M^\ep_4=\frac{\kappa\mathrm{g}(\rho^2\!-\!\rho^1)}{2\pi (\mu^2\!+\!\mu^1)A^\ep(t)}\int_\T
\la\big(\phi_\ep\ast(\da^3 z^\ep_1)\big)(\al)\phi_\ep\ast(\da z^\ep_2\da^3z^\ep_2)(\al) d\alpha.
$$
Using \eqref{conconvo}, we get
\begin{align*}
M^\ep_1&\leq C\|\F(z^\ep)\|_{L^\infty}\|\la\big(\phi_\ep\ast(\da^3 z^\ep_1
\da z^\ep_2)-\phi_\ep\ast(\da^3 z^\ep_1) \da z^\ep_2\big)\|_{L^2}\|\phi_\ep\ast(\da^3z^\ep_2)\|^2_{L^2}\\
&\leq C\|\F(z^\ep)\|_{L^\infty}\|\phi_\ep\ast(\da^3 z^\ep_1
\da z^\ep_2)-\phi_\ep\ast(\da^3 z^\ep_1) \da z^\ep_2\|_{H^1}\|\da^3z^\ep_2\|^2_{L^2}\\
&\leq C\|\F(z^\ep)\|_{L^\infty}\|\da^3 z^\ep_1\|_{L^2}\|\da z^\ep_2\|_{C^1}\|\da^3z^\ep_2\|^2_{L^2},
\end{align*}
and therefore
$$
M^\ep_1\leq \exp C(\|\F(z^\ep)\|^2_{L^\infty}+\| z^\ep\|^2_{H^3}).
$$
For $M^\ep_2$ we use the commutator estimate for the operator
$\la$ to obtain
$$
M^\ep_2\leq C\|\F(z^\ep)\|_{L^\infty}\|\phi_\ep\ast(\da^3 z^\ep_1)\|_{L^2}\|\da z^\ep_2\|_{C^{1,\delta}}\|\phi_\ep\ast(\da^3z^\ep_2)\|_{L^2}\leq \exp C(\|\F(z^\ep)\|^2_{L^\infty}+\| z^\ep\|^2_{H^3}).
$$
Regarding $M_3^\ep$ we have
$$
M^\ep_3=\frac{\kappa\mathrm{g}(\rho^2\!-\!\rho^1)}{2\pi (\mu^2\!+\!\mu^1)A^\ep(t)}\int_\T
\phi_\ep\ast(\da^3 z^\ep_1)(\al)\la\big(\da z^\ep_2\phi_\ep\ast(\da^3z^\ep_2)-\phi_\ep\ast(\da z^\ep_2\da^3z^\ep_2)\big)(\al) d\alpha,
$$
showing that it can be estimated as $M_1^\ep.$

 The identity
$$
\da z^\ep_2(\al)\da^3z^\ep_2(\al)=-\da z^\ep_1(\al)\da^3z^\ep_1(\al)-|\da^2 z^\ep(\al)|^2,
$$
allow us to write $M^\ep_4$ as the sum of $N^\ep_1$ and $N^\ep_2$
where
$$
N^\ep_1=-\frac{\kappa\mathrm{g}(\rho^2\!-\!\rho^1)}{2\pi (\mu^2\!+\!\mu^1)A^\ep(t)}\int_\T\la\big(\phi_\ep\ast\da^3 z^\ep_1\big)(\al)\phi_\ep\ast(|\da^2 z^\ep|^2)(\al) d\alpha,
$$
and
$$
N^\ep_2=-\frac{\kappa\mathrm{g}(\rho^2\!-\!\rho^1)}{2\pi (\mu^2\!+\!\mu^1)A^\ep(t)}\int_\T\phi_\ep\ast(\da z^\ep_1\da^3z^\ep_1)(\al)\la\big(\phi_\ep\ast\da^3 z^\ep_1\big)(\al) d\alpha.
$$
Then an integration by parts shows that $$N^\ep_1\leq
C\|\F(z^\ep)\|_{L^\infty}\|z^\ep\|_{C^2}\|\da^3
z^\ep\|^2_{L^2}\leq \exp C(\|\F(z^\ep)\|^2_{L^\infty}+\|
z^\ep\|^2_{H^3}).$$ Using again the identity \eqref{conconvo} in
$N^\ep_2$, we obtain finally
\begin{align*}
L^\ep_1&\leq -\frac{\kappa\mathrm{g}(\rho^2\!-\!\rho^1)}{2\pi (\mu^2\!+\!\mu^1)A^\ep(t)}\int_\T\da z^\ep_1(\al)\phi_\ep\ast(\da^3 z^\ep_1)(\al)
\la\big(\phi_\ep\ast(\da^3z^\ep_1)\big)(\al) d\alpha\\
&\quad +\exp C(\|\F(z^\ep)\|^2_{L^\infty}+\|z^\ep\|^2_{H^3}).
\end{align*}
In a similar way we get for $L^\ep_2$
\begin{align*}
L^\ep_2&\leq -\frac{\kappa\mathrm{g}(\rho^2\!-\!\rho^1)}{2\pi (\mu^2\!+\!\mu^1)A^\ep(t)}\int_\T\da z^\ep_1(\al)\phi_\ep\ast(\da^3 z^\ep_2)(\al)
\la\big(\phi_\ep\ast(\da^3z^\ep_2)\big)(\al) d\alpha\\
&\quad +\exp C(\|\F(z^\ep)\|^2_{L^\infty}+\|z^\ep\|^2_{H^3}),
\end{align*}
 giving us
\begin{align*}
K^\ep_2&\leq -\frac{\kappa\mathrm{g}(\rho^2\!-\!\rho^1)}{2\pi (\mu^2\!+\!\mu^1)A(t)}\int_\T\da z^\ep_1(\al)\phi_\ep\ast(\da^3 z)(\al)\cdot\la\big(\phi_\ep\ast(\da^3z)\big)(\al) d\alpha\\
&\quad +\exp C(\|\F(z^\ep)\|^2_{L^\infty}+\| z^\ep\|^2_{H^3}).
\end{align*}
The formula for $\sigma^\ep(\al,t)$ begins to appear in the
non-integrable terms.  Using a similar method for the rest of the
non-integrable terms we obtain the inequality \eqref{ntniepsilon}
for $k=3$.

The next step is to integrate the system during a time $T$
independent of $\ep$. First let us observe that if $z_0(\al)\in
H^k$, then we have the solution $z^\ep\in C^1([0,T^\ep];H^k)$. And
if initially $\sigma(\al,0)>0$, there is a time depending on
$\ep$, denoted by $T^\ep$ again, in which $\sigma^\ep(\al,t)>0$.
Now, for $t\leq T^\ep$ we have \eqref{ntniepsilon}, and then we
use the following pointwise inequality (see \cite{CC}):
$$
f(\al)\la f(\al)-\frac{1}{2}\la(f^2)(\al)\geq 0,
$$
to obtain
\begin{align*}
\frac{d}{dt}\|z^\ep\|^2_{H^k}(t)&\leq I+\exp C(\|\F(z^\ep)\|^2_{L^\infty}(t)+\|z^\ep\|^2_{H^k}(t)),
\end{align*}
where
$$
I=-\frac{\kappa}{2\pi(\mu_1\!+\!\mu_2)A^\ep(t)}\!\int_\T
\sigma^\ep(\al,t)\frac{1}{2}\la(|\phi_\ep\ast (\da^k
z^\ep)|^2)(\al,t) d\al.
$$
We have
$$
\|\la(\sigma^\ep)\|_{L^\infty}(t)\leq C\|\sigma^\ep\|_{H^2}(t)\leq C(\|BR(z^\ep,\varpi^\ep)\|_{L^2}(t)+\|\da^2BR(z^\ep,\varpi^\ep)\|_{L^2}(t)+1)\|z\|_{H^3}(t),
$$
and writing
$$
I=-\frac{\kappa}{2\pi(\mu_1\!+\!\mu_2)A^\ep(t)}\!\int_\T
\la(\sigma^\ep)(\al,t)\frac{1}{2}|\phi_\ep\ast (\da^k
z^\ep)|^2(\al,t) d\al,
$$
we obtain
$$
I\leq C\|\F(z^\ep)\|_{L^\infty}\|\la(\sigma^\ep)\|_{L^\infty}\|\da^k z^\ep\|^2_{L^2}\leq \exp C(\|\F(z^\ep)\|^2_{L^\infty}+\| z^\ep\|^2_{H^k}).
$$
Finally, for $t\leq T^\ep$ we have
\begin{align}\label{porfin}
\frac{d}{dt}\|z^\ep\|^2_{H^k}(t)&\leq C\exp C(\|\F(z^\ep)\|^2_{L^\infty}(t)+\|z^\ep\|^2_{H^k}(t)).
\end{align}
We have also (see section \ref{arc-chord}):
\begin{align*}
\frac{d}{dt}\|\F(z^\ep)\|^2_{L^\infty}(t)&\leq C\exp C(\|\F(z^\ep)\|^2_{L^\infty}(t)+\|z^\ep\|^2_{H^3}(t)),
\end{align*}
and from \eqref{porfin} it follows that
\begin{align*}
\frac{d}{dt}(\|z^\ep\|^2_{H^k}(t)+\|\F(z^\ep)\|^2_{L^\infty}(t))&\leq C\exp C(\|z^\ep\|^2_{H^k}(t)+\|\F(z^\ep)\|^2_{L^\infty}(t))
\end{align*}
for $t\leq T^\ep$. Integrating
\begin{align}\label{porfin2}
\|z^\ep\|^2_{H^k}(t)+\|\F(z^\ep)\|^2_{L^\infty}(t)&\leq -\frac{1}{C}
\ln \big(-t+\exp(-C(\|z_0\|^2_{H^k}+\|\F(z_0)\|^2_{L^\infty}))\big),
\end{align}
$t\leq T^\ep$. Let us mention that at this point of the proof we
can not assume local-existence, because we have the above estimate
for $t\leq T^\ep$, and if we let $\ep\rightarrow 0$, it could be
possible that $T^\ep\rightarrow 0$ i.e. we cannot assume that if
the initial data satisfy $\sigma(\al,0)>0$,  there must be  a time
$T$, independent of $\ep$, in which \eqref{porfin2} is satisfied.
In other words, at this stage of the proof  we do not have
local-existence when $\ep\rightarrow 0$. But since in the
evolution equation everything depends upon the sign of
$\sigma^\ep(\al,t)$,  the  following argument will allow us to
continue the proof. First let us observe that as in section
\ref{minimumsigma} we have
\begin{equation}\label{queseacaba}
m^\ep(t)\geq m(0)-\int_0^t \exp C(\|\F(z^\ep)\|^2_{L^\infty}(s)+\|z^\ep\|^2_{H^3}(s)) ds,
\end{equation}
where
$$
m^\ep(t)=\min_{\al\in\T}\sigma^\ep(\al,t),
$$
and $t\leq T^\ep$. Using \eqref{porfin2} in \eqref{queseacaba} we get
\begin{equation}\label{unpocomas}
m^\ep(t)\geq m(0)+C(\|z_0\|^2_{H^k}+\|\F(z_0)\|^2_{L^\infty})+\ln\big(-t+\exp(-C(\|z_0\|^2_{H^k}+\|\F(z_0)\|^2_{L^\infty}))\big),
\end{equation}
for $t\leq T^\ep$. Using \eqref{unpocomas} and \eqref{porfin2},
now we find that if $\ep\rightarrow 0$, then $T^\ep\nrightarrow
0$, because if we take $T=\min(T_1,T_2)$ where $T_1$ satisfies
$$
m(0)+C(\|z_0\|^2_{H^k}+\|\F(z_0)\|^2_{L^\infty})+\ln\big(-T_1+\exp(-C(\|z_0\|^2_{H^k}+\|\F(z_0)\|^2_{L^\infty}))\big)>0,
$$
and $T_2$
$$
-\frac{1}{C}
\ln \big(-T_2+\exp(-C(\|z_0\|^2_{H^k}+\|\F(z_0)\|^2_{L^\infty}))\big)<\infty.
$$
For $t\leq T$ we have $m^\ep(t)>0$ and
$$
\|z^\ep\|^2_{H^k}(t)+\|\F(z^\ep)\|^2_{L^\infty}(t)\leq -\frac{1}{C}
\ln \big(-T+\exp(-C(\|z_0\|^2_{H^k}+\|\F(z_0)\|^2_{L^\infty}))\big)<\infty,
$$
and $T$ only depends on the initial data $z_0$. Now we let  $\ep$
tends to $0$, to conclude the existence result.

%%%%%%%%%%%%%%%%%%%%%%%%%%%%%%%%%%%%%%%%%%%%%%%%%%%%%%%%%%%%%%%%%%%%%%%%%%%%%%%%%%%%%%%%%%%%%%%%%%%%%%%%%%%%%%%%%%%%%%%%%%%
%%%%%%%%%%%%%%%%%%%%%%%%%%%%%%%%%%%%%%%%%%%%%%%%%%%%%%%%%%%%%%%%%%%%%%%%%%%%%%%%%%%%%%%%%%%%%%%%%%%%%%%%%%%%%%%%%%%%%%%%%%%

\section{Appendix}

%%%%%%%%%%%%%%%%%%%%%%%%%%%%%%%%%%%%%%%%%%%%%%%%%%%%%%%%%%%%%%%%%%%%%%%%%%%%%%%%%%%%%%%%%%%%%%%%%%%%%%%%%%%%%%%%%%%%%%%%%%%
%%%%%%%%%%%%%%%%%%%%%%%%%%%%%%%%%%%%%%%%%%%%%%%%%%%%%%%%%%%%%%%%%%%%%%%%%%%%%%%%%%%%%%%%%%%%%%%%%%%%%%%%%%%%%%%%%%%%%%%%%%%

Let us denote
\begin{align*}
V(\al,\beta)&=(V_1(\al,\beta), V_2(\al,\beta))=(\tan(\frac{z_1(\al)-z_1(\beta)}{2}),\tanh(\frac{z_2(\al)-z_2(\beta)}{2})),
\end{align*}
and
\begin{align*}
W(\al,\beta)&=(W_1(\al,\beta),W_2(\al,\beta))=((\frac{z_1(\al)-z_1(\beta)}{2})_p,(\frac{z_2(\al)-z_2(\beta)}{2})_p),
\end{align*}
where $(\al)_p$ is the periodic extension of the function $\al$ in $\T$. We give the following equalities for the hyperbolic tangent function:
\begin{equation}\label{ac1}
(\tanh(\al)-(\al)_p)/\tanh^2(\al)=(\al)_pf(\al)\quad\mbox{ with
}\quad f\in L^{\infty}(\R),
\end{equation}

\begin{equation}\label{ac2}
(\tanh(\al)-(\al)_p)/\tanh^3(\al)=g(\al)\quad\mbox{ with
}\quad g\in L^\infty(\R).
\end{equation}
For the tangent function it holds
\begin{equation}\label{ac3}
(\tan(\al/2)-(\al/2)_p)/\tan(\al/2)=(\al/2)_p h(\al)\quad\mbox{ with
}\quad h\in L^\infty(\R),
\end{equation}

\begin{equation}\label{ac4}
(\tan(\al/2)-(\al/2)_p)/\tan^2(\al/2)=(\al/2)_p j(\al)\quad\mbox{ with
}\quad j\in L^\infty(\R),
\end{equation}

\begin{equation}\label{ac5}
(\tan(\al/2)-(\al/2)_p)/|\tan^3(\al/2)|=k(\al)\quad\mbox{ with }\quad
k\in L^{\infty}(\R).
\end{equation}
Also we shall use that the below functions are bounded on $[-\pi,\pi]$
\begin{equation}\label{ac6}
2/\al-1/\tan(\al/2),\, 4/\al^2-1/\sin^2(\al/2)\in L^{\infty}(\T),
\end{equation}
and the following estimates:
\begin{equation}\label{nh1}
|W(\al,\al-\beta)-\da z(\al)\beta/2|\leq
\frac{1}{2}\|z\|_{C^2}|\beta|^2,
\end{equation}
\begin{equation}\label{nh2}
|W(\al,\al-\beta)-\da z(\al)\beta/2-\da^2 z(\al)\beta^2/4|\leq
\frac12 \|z\|_{C^{2,\delta}}|\beta|^{2+\delta}.
\end{equation}

\begin{lemma}\label{A12}
Given
$$
A_1(\al,\al-\beta)=\frac{V_2(\al,\al-\beta)}{|V(\al,\al-\beta)|^2}
-\frac{1}{|\da z(\al)|^2}\frac{\da
z_2(\al)}{\tan(\frac{\beta}{2})}
$$

$$
A_2(\al,\al-\beta)=\frac{V_1(\al,\al-\beta)}{|V(\al,\al-\beta)|^2}
-\frac{1}{|\da z(\al)|^2}\frac{\da
z_1(\al)}{\tan(\frac{\beta}{2})},
$$
we have
$$\|A_1(\al,\al-\beta)\|_{L^{\infty}}\leq \|\F(z)\|_{L^{\infty}}\|z\|_{C^2},$$
and
$$\|A_2(\al,\al-\beta)\|_{L^{\infty}}\leq \|\F(z)\|_{L^{\infty}}\|z\|^2_{C^2},$$
\end{lemma}
Proof: We introduce the splitting $A_1(\al,\al-\beta)=I_1+I_2+I_3+I_4$ where

$$I_1=\frac{\tanh(\frac{z_2(\al)-z_2(\al-\beta)}{2})-(\frac{z_2(\al)-z_2(\al-\beta)}{2})_p}{V(\al,\al-\beta)}, $$

$$I_2=\F(z)(\al,\beta)\frac{((z_2(\al)-z_2(\al-\beta))/2)_p-\da z_2(\al)\beta/2}{\beta^2/4},$$

$$I_3=\frac{\da z_2(\al)}{\beta/2}(\F(z)(\al,\beta)-\frac{1}{|\da z(\al)|^2}),$$

$$I_4=\frac{\da z_2(\al)}{|\da z(\al)|^2}(\frac2\beta-\frac{1}{\tan(\frac\beta2)}),$$
and $\F(z)(\al,\beta)$ was defined in \eqref{df}.

Since
\begin{align*}
\begin{split}
I_1&=\frac{1}{1+\frac{V^2_1(\al,\al-\beta)}{V^2_2(\al,\al-\beta)}}f(\frac{z_2(\al)-z_2(\al-\beta)}{2})
(\frac{z_2(\al)-z_2(\al-\beta)}{2})_p
\end{split}
\end{align*}
by \eqref{ac1}, we get $I_1\leq C$. Also $I_2\leq \|\F(z)\|_{L^{\infty}}\|z\|_{C^2}$ using \eqref{nh1}, and $I_4\leq C\|\F(z)\|_{L^{\infty}}^{1/2}.$ We rewrite
\begin{align*}
\begin{split}
I_3&=\frac{\da z_2(\al)}{\beta/2}\frac{(\da
z(\al)\beta/2+V(\al,\al-\beta))\cdot
(\da z(\al)\beta/2-V(\al,\al-\beta))}{|\da z(\al)|^2|V(\al,\al-\beta)|^2},\\
\end{split}
\end{align*}
and split further
$$
I_3=J_1+J_2,
$$
where
$$
J_1=\frac{\da z_2(\al)}{\beta/2}\frac{(\da
z_1(\al)\beta/2+V_1(\al,\al-\beta)) (\da
z_1(\al)\beta/2-V_1(\al,\al-\beta))}{|\da
z(\al)|^2|V(\al,\al-\beta)|^2},
$$
$$
J_2=\frac{\da z_2(\al)}{\beta/2}\frac{(\da
z_2(\al)\beta/2+V_2(\al,\al-\beta)) (\da
z_2(\al)\beta/2-V_2(\al,\al-\beta))}{|\da
z(\al)|^2|V(\al,\al-\beta)|^2}.
$$
We continue as follows
$$
J_1=K_1+K_2,
$$
for
$$
K_1=\frac{\da z_2(\al)\da z_1(\al) (\da
z_1(\al)\beta/2-V_1(\al,\al-\beta))}{|\da
z(\al)|^2|V(\al,\al-\beta)|^2},
$$
$$
K_2=\frac{\da z_2(\al)V_1(\al,\al-\beta)(\da
z_1(\al)\beta/2-V_1(\al,\al-\beta))}{|\da
z(\al)|^2|V(\al,\al-\beta)|^2\beta/2},
$$
to take
$K_1=L_1+L_2$,
$$
L_1=\frac{\da z_2(\al)\da z_1(\al)
((\frac{z_1(\al)-z_1(\al-\beta)}{2})_p-V_1(\al,\al-\beta))}{|\da
z(\al)|^2|V(\al,\al-\beta)|^2},
$$
$$
L_2=\frac{\da z_2(\al)\da z_1(\al) (\da
z_1(\al)\beta/2-(\frac{z_1(\al)-z_1(\al-\beta)}{2})_p)}{|\da
z(\al)|^2|V(\al,\al-\beta)|^2}.
$$
We find
$$
L_1=\frac{\da z_2(\al)\da z_1(\al)}{|\da
z(\al)|^2}\frac{1}{1+\frac{V^2_2(\al,\al-\beta)}{V^2_1(\al,\al-\beta)}}
\frac{(\frac{z_1(\al)-z_1(\al-\beta)}{2})_p-V_1(\al,\al-\beta)}{V^2_1(\al,\al-\beta)}
$$
and  using \eqref{ac4}, we obtain $L_1\leq C$.

Since
$$
L_2=\frac{\da z_2(\al)\da z_1(\al)}{|\da
z(\al)|^2}\F(z)(\al,\beta) \frac{(\da
z_1(\al)\beta/2-(\frac{z_1(\al)-z_1(\al-\beta)}{2})_p)}{\beta^2/4}
$$
we have $L_2\leq C\|\F(z)\|_{L^{\infty}}\|z\|_{C^2}$. Next let us
write $K_2=L_3+L_4$, for
$$
L_3=\frac{\da
z_2(\al)V_1(\al,\al-\beta)((\frac{z_1(\al)-z_1(\al-\beta)}{2})_p-V_1(\al,\al-\beta))}{|\da
z(\al)|^2|V(\al,\al-\beta)|^2\beta/2},
$$
$$
L_4=\frac{\da z_2(\al)V_1(\al,\al-\beta)(\da
z_1(\al)\beta/2-(\frac{z_1(\al)-z_1(\al-\beta)}{2})_p)}{|\da
z(\al)|^2|V(\al,\al-\beta)|^2\beta/2}.
$$
In a similar way we find that
$$
L_3=\frac{\da z_2(\al)}{|\da z(\al)|^2}
\frac{1}{1+\frac{V^2_2(\al,\al-\beta)}{V^2_1(\al,\al-\beta)}}\frac{(\frac{z_1(\al)-z_1(\al-\beta)}{2})_p-
V_1(\al,\al-\beta)}{V_1(\al,\al-\beta)}\frac{2}{\beta}.
$$ By \eqref{ac3} one gets
$$
L_3\leq
C\frac{|\da z_2(\al)|}{|\da z(\al)|^2}\frac{|(\frac{z_1(\al)-z_1(\al-\beta)}{2})_p|}{|\beta|/2}\leq C.
$$
As before we conclude that
$$
L_4\leq C\|\F(z)\|_{L^{\infty}}\frac{|\da
z_1(\al)\beta/2-(\frac{z_1(\al)-z_1(\al-\beta)}{2})_p|}{|\beta|^2}\leq
C\|\F(z)\|_{L^{\infty}}\|z\|_{C^2}.
$$

We consider now $J_2=K_3+K_4$, where
$$
K_3=\frac{|\da z_2(\al)|^2}{|\da z(\al)|^2} \frac{\da
z_2(\al)\beta/2-V_2(\al,\al-\beta)}{|V(\al,\al-\beta)|^2}
$$
and
$$
K_4=\frac{\da z_2(\al)}{|\da
z(\al)|^2}\frac{V_2(\al,\al-\beta)}{|V(\al,\al-\beta)|^2}
\frac{\da z_2(\al)\beta/2-V_2(\al,\al-\beta)}{\beta/2}.
$$
Using \eqref{ac1}, we find

\begin{align*}
\begin{split}
K_3 &\leq C+\frac{|\da z_2(\al)|^2}{|\da z(\al)|^2}\frac{|\da z_2(\al)\beta/2-(\frac{z_2(\al)-z_2(\al-\beta)}{2})_p|}{|V(\al,\al-\beta)|^2}\leq C\|\F(z)\|_{L^{\infty}}\|z\|_{C^2}.
\end{split}
\end{align*}
and
\begin{align*}
\begin{split}
K_4 &\leq \frac{|\da z_2(\al)|}{|\da
z(\al)|^2}\Big(\frac{|V_2(\al,\al-\beta)-W_2(\al,\al-\beta)|}{(1\!+\!
\frac{(V_1(\al,\al\!-\!\beta))^2}{(V_2(\al,\al\!-\!\beta))^2})|V_2(\al,\al\!-\!\beta)||\beta/2|}+\frac{|\da z_2(\al)\beta/2\!-\!W_2(\al,\al\!-\!\beta)|}{|V(\al,\al\!-\!\beta)||\beta/2|}\Big)\\
& \leq C\|\F(z)\|^{1/2}_{L^{\infty}}\frac{|(\frac{z_2(\al)-z_2(\al-\beta)}{2})_p|}{|\beta/2|}+\|\F(z)\|_{L^{\infty}}\frac{|\da z_2(\al)\beta/2-(\frac{z_2(\al)-z_2(\al-\beta)}{2})_p|}{(\beta/2)^2}\\
&\leq C\|\F(z)\|_{L^{\infty}}\|z\|_{C^2},
\end{split}
\end{align*}
that is $K_3+K_4\leq C\|\F(z)\|_{L^{\infty}}\|z\|_{C^2}$.

Putting all the previous estimates together we get
$|A_1(\al,\al-\beta)|\leq C\|\F(z)\|_{L^{\infty}}\|z\|_{C^2}$.

Regarding
$$
A_2(\al,\al-\beta)=\frac{V_1(\al,\al-\beta)}{|V(\al,\al-\beta)|^2}
-\frac{1}{|\da z(\al)|^2}\frac{\da z_1(\al)}{\tan(\beta/2)}
$$
 we have the splitting $A_2=I_5+I_6+I_7+I_8$, where
$$
I_5=\frac{V_1(\al,\al-\beta)-(\frac{z_1(\al)-z_1(\al-\beta)}{2})_p}{|V(\al,\al-\beta)|^2},
$$

$$
I_6=\F(z)(\al,\beta)\frac{(\frac{z_1(\al)-z_1(\al-\beta)}{2})_p-\da
z_1(\al)\beta/2}{\beta^2/4},
$$

$$
I_7=\frac{\da z_1(\al)}{\beta/2}(\F(z)(\al,\beta)-\frac{1}{|\da
z(\al)|^2}),
$$

$$
I_8=\frac{\da z_1(\al)}{|\da
z(\al)|^2}(\frac{2}{\beta}-\frac{1}{\tan(\beta/2)}).
$$
Then the same arguments used above allows us to obtain $|A_2|\leq
C\|\F(z)\|_{L^{\infty}}\|z\|_{C^2}$.

\begin{lemma}\label{B}
Let $B(\al,\beta)$ be defined by
$$B(\al,\al-\beta)=V_1(\al,\al-\beta)\frac{V(\al,\al-\beta)^{\bot}\cdot\da z(\al)}{|V(\al,\al-\beta)|^4}-\da z_1(\al)\frac{(\da^2 z(\al))^{\bot}\cdot\da z(\al)}{|\da z(\al)|^4\tan(\beta/2)}$$
Then it satisfies the inequality
$$|B(\al,\al-\beta)|\leq C\|\F(z)\|^2_{L^{\infty}}\|z\|^3_{C^{2,\delta}}|\beta|^{\delta-1}.$$
\end{lemma}

Proof: Let us  decompose  $B(\al,\beta)=I_1+I_2$ where
$$
I_1=(V_1(\al,\al-\beta)-(\frac{z_1(\al)-z_1(\al-\beta)}{2})_p)
\frac{V(\al,\al-\beta)^{\bot}\cdot\da
z(\al)}{|V(\al,\al-\beta)|^4},
$$
and
$$
I_2=(\frac{z_1(\al)-z_1(\al-\beta)}{2})_p
\frac{V(\al,\al-\beta)^{\bot}\cdot\da
z(\al)}{|V(\al,\al-\beta)|^4}-\da z_1(\al)\frac{(\da^2
z(\al))^{\bot}\cdot\da z(\al)}{|\da z(\al)|^4\tan(\beta/2)}.
$$
Using the identity \eqref{ac5}, we can rewrite $I_1$ as follows:
$$
I_1=k(z_1(\al)-z_1(\al-\beta))
\frac{1}{\big(1+\frac{V^2_2(\al,\beta)}{V^2_1(\al,\beta)}\big)^{3/2}}
\frac{V(\al,\al-\beta)^{\bot}\cdot\da z(\al)}{|V(\al,\al-\beta)|}
$$
to get $|I_1|\leq C\|z\|_{C^1}$.

Next we consider $I_2=J_1+J_2$, where
$$
J_1=W_1(\al,\al-\beta)\frac{(V(\al,\al-\beta)-W(\al,\al-\beta))^{\bot}\cdot\da
z(\al)}{|V(\al,\al-\beta)|^4},
$$
and
$$
J_2=W_1(\al,\al-\beta)\frac{W(\al,\al-\beta)^{\bot}\cdot\da
z(\al)}{|V(\al,\al-\beta)|^4}-\da z_1(\al)\frac{(\da^2
z(\al))^{\bot}\cdot\da z(\al)}{|\da z(\al)|^4\tan(\beta/2)}.
$$
Using \eqref{ac2}, \eqref{ac5}, and the fact that
$(\beta/2)_p/\tan(\beta/2)$ is bounded, we obtain $|J_1|\leq
C\|z\|_{C^1}$. To continue we can rewrite $J_2$ as follows:
$$
J_2=W_1(\al,\al-\beta) \frac{(W(\al,\al-\beta)-\da
z(\al)\beta/2)^{\bot}\cdot\da z(\al)}{|V(\al,\al-\beta)|^4} -\da
z_1(\al)\frac{(\da^2 z(\al))^{\bot}\cdot\da z(\al)}{|\da
z(\al)|^4\tan(\beta/2)},
$$
and $J_2=K_1+K_2+K_3+K_4$, where
$$
K_1=(W_1(\al,\al-\beta)-\da z_1(\al)\beta/2)
\frac{(W(\al,\al-\beta)-\da z(\al)\beta/2)^{\bot}\cdot\da
z(\al)}{|V(\al,\al-\beta)|^4},
$$
$$
K_2=\da z_1(\al)\beta/2 \frac{(W(\al,\al-\beta)-\da
z(\al)\beta/2-\da^2 z(\al)\beta^2/4 )^{\bot}\cdot\da
z(\al)}{|V(\al,\al-\beta)|^4},
$$

$$
K_3=2\da z_1(\al)(\da^2 z(\al)^{\bot}\cdot\da
z(\al))(\F(z)(\al,\beta)^2-\frac{1}{|\da z(\al)|^4})/\beta,
$$

$$
K_4=\frac{\da z_1(\al)(\da^2 z(\al))^{\bot}\cdot\da z(\al)}{|\da
z(\al)|^4}(\frac{2}{\beta}-\frac{1}{\tan(\beta/2)}).
$$
Clearly we have $|K_4|\leq C\|\F(z)\|_{L^{\infty}}\|z\|_{C^2}$,
and using \eqref{nh1} we obtain $|K_1|\leq C\|\F(z)\|^2_{L^{\infty}}\|z\|^3_{C^2}$.
Furthermore the estimate \eqref{nh2} allows us to obtain  $|K_2|\leq
C\|\F(z)\|^2_{L^{\infty}}\|z\|^3_{C^{2,\delta}}|\beta|^{\delta-1}$.
Next we consider in $K_3$ the factor $L(\al,\beta)$ given by
$$
L(\al,\beta)=(\F(z)(\al,\beta)^2-\frac{1}{|\da z(\al)|^4})/\beta.
$$
We can write $L(\al,\beta)$ as follows:
\begin{equation}\label{nscll}
\frac{(|\da
z(\al)|^2\beta^2/4\!+\!|V(\al,\al\!-\!\beta)|^2)}{|\da
z(\al)|^4|V(\al,\al\!-\!\beta)|^2}\frac{(\da
z(\al)\beta/2\!+\!V(\al,\al\!-\!\beta))\!\cdot\! (\da
z(\al)\beta/2\!-\!V(\al,\al\!-\!\beta))}{|V(\al,\al\!-\!\beta)|^2\beta}.
\end{equation}
Then proceeding as in the previous lemma we get  $|K_3|\leq
C\|\F(z)\|^2_{L^{\infty}}\|z\|^3_{C^2}$ and this ends the proof.
q.e.d.

\begin{lemma}\label{C}
Given $C(\al,\beta)$ by the following equality
$$
C(\al,\al-\beta)=\frac{V^\bot(\al,\al-\beta)\varpi(\al-\beta)\beta}{|V(\al,\al-\beta)|^4}-\frac{2\da^{\bot} z(\al)\varpi(\al)}{|\da^2 z(\al)|^4
\sin^2(\beta/2)},
$$
we obtain
$$
|C(\al,\al-\beta)|\leq C\|\F(z)\|^2_{L^\infty}\|z\|_{C^2}\|\varpi\|_{C^1}\frac{1}{|\beta|}.
$$
\end{lemma}
Proof: We decompose $C(\al,\al-\beta)=I_1+I_2+I_3+I_4+I_5$ where
$$
I_1=\frac{(V(\al,\al-\beta)-W(\al,\al-\beta))^{\bot}\varpi(\al-\beta)\beta}{|V(\al,\al-\beta)|^4},
$$
$$
I_2=\frac{(W(\al,\al-\beta)-\da z(\al)\beta/2)^{\bot}\varpi(\al-\beta)\beta}{|V(\al,\al-\beta)|^4},
$$
$$
I_3=\frac{\dpa z(\al)\beta^2(\varpi(\al-\beta)-\varpi(\al))}{2|V(\al,\al-\beta)|^4},
$$
$$
I_4=8\dpa z(\al)\varpi(\al)(\F(z)(\al,\beta)^2-\frac{1}{|\da z(\al)|^4})/\beta^2,
$$
$$
I_5=2\frac{\dpa z(\al)\varpi(\al)}{|\da z(\al)|^4}(4/\beta^2-\sin^2(\beta/2)).
$$
Using \eqref{ac2} and \eqref{ac5} we get $|I_1|\leq C \|\F(z)\|^{1/2}_{L^\infty}\|\varpi\|_{L^\infty}$. Using \eqref{nh1}
clearly we obtain $|I_2|\leq C\|\F(z)\|^2_{L^\infty}\|z\|_{C^2}\|\varpi\|_{L^\infty}/|\beta|$. For the next term it holds
$$|I_3|\leq C\|\F(z)\|^2_{L^\infty}\|z\|_{C^1}\|\varpi\|_{C^1}/|\beta|.$$
The reference \eqref{ac6} gives $|I_5|\leq C\|\F(z)\|^{3/2}_{L^\infty}\|\varpi\|_{L^\infty}$. Finally, the estimate given in the previous lemma for the term
$$
(\F(z)(\al,\beta)^2-\frac{1}{|\da z(\al)|^4})/\beta,
$$
written in \eqref{nscll} allows us to conclude $|I_4|\leq\|\F(z)\|^2_{L^\infty}\|z\|_{C^1}\|\varpi\|_{L^\infty}/|\beta|$.

\begin{lemma}\label{Q1}
Let $Q_1(\al,\beta)$ be given by
$$
Q_1(\al,\al-\beta)=-\frac{(\da z(\al))^{\bot}}{2}\big(\frac{1}{|V(\al,\al-\beta)|^2}-\frac{4}{|\da z(\al)|^2|\beta|^2}
+\frac{2\da z(\al)\cdot \da^2 z(\al)}{|\da z(\al)|^4\beta}\big).
$$
Then it satisfies the estimate
$\|Q_1\|_{L^{\infty}}\leq
\|\F(z)\|^k_{L^\infty}\|z\|^k_{C^{2,\delta}}|\beta|^{\delta-1}$.

\end{lemma}

Proof: To simplify we will consider
$$
C(\al,\al-\beta)=\frac{1}{|V(\al,\al-\beta)|^2}-\frac{4}{|\da z(\al)|^2|\beta|^2}+\frac{4\da z(\al)\cdot \da^2 z(\al)}{|\da z(\al)|^4\beta},
$$
and we will show that $\|C\|_{L^{\infty}}\leq \|\F(z)\|^k_{L^\infty}\|z\|^k_{C^{2,\delta}}|\beta|^{\delta-1}$. We can rewrite
$$
C(\al,\al-\beta)=\frac{(\da z(\al)\beta+2V(\al,\al\!-\!\beta))\cdot(\da z(\al)\beta-2V(\al,\al\!-\!\beta))}{|V(\al,\al-\beta)|^2|\da z(\al)|^2|\beta|^2}+\frac{4\da z(\al)\cdot \da^2 z(\al)}{|\da z(\al)|^4\beta},
$$
and then take $C(\al,\al-\beta)=I_1+I_2+I_3$ where

$$
I_1=-\frac{|2V(\al,\al-\beta)-\da z(\al)\beta|^2}{|V(\al,\al-\beta)|^2|\da z(\al)|^2|\beta|^2},
$$
$$
I_2=-\frac{2 \da z(\al)\beta\cdot(2V(\al,\al-\beta)-\da z(\al)\beta-\da^2 z(\al)\beta^2/2)}{|V(\al,\al-\beta)|^2|\da z(\al)|^2|\beta|^2},
$$
$$
I_3=-\frac{\da z(\al)\cdot\da^2 z(\al)\beta}{|\da z(\al)|^2}(\frac{1}{|V(\al,\al-\beta)|^2}-\frac{4}{|\da z(\al)|^2|\beta|^2}).
$$
Since
$$
|I_1|\leq \frac{|2V(\al,\al-\beta)-2W(\al,\al-\beta)|^2}{|V(\al,\al-\beta)|^2|\da z(\al)|^2|\beta|^2}+\frac{|2W(\al,\al-\beta)-\da z(\al)\beta|^2}{|V(\al,\al-\beta)|^2|\da z(\al)|^2|\beta|^2},
$$
using \eqref{ac1}, \eqref{ac3} and the inequality \eqref{nh1} we control the term $I_1$. For $I_2$ it holds
$$
|I_2|\leq\frac{4|V(\al,\al-\beta)-W(\al,\al-\beta)|}{|V(\al,\al-\beta)|^2|\da z(\al)||\beta|}+
\frac{4|W(\al,\al-\beta)-\da z(\al)\beta-\da^2 z(\al)\beta^2/2|}{|V(\al,\al-\beta)|^2|\da z(\al)||\beta|},
$$
and using \eqref{ac1}, \eqref{ac4}, and \eqref{nh2} we get the appropriate inequality. For $I_3$ we write
$$
I_3=-\frac{4\da z(\al)\cdot\da^2 z(\al)}{|\da z(\al)|^2}(\F(z)(\al,\beta)-\frac{4}{|\da z(\al)|^2|\beta|^2})/\beta,
$$
and proceed as before.

\subsection*{{\bf Acknowledgements}}

\smallskip
The first author was partially supported by the grant {\sc MTM2005-04730} of the MEC (Spain).
The other two authors were partially supported by the grant {\sc MTM2005-05980} of the MEC (Spain).

%%%%%%%%%%%%%%%%%%%%%%%%%%%%%%%%%%%%%%%%%%%%%%%%%%%%%%%%%%%%%%%%%%%%%%%%%%%%%%%%%%%

\begin{quote}
\begin{tabular}{l}
\textbf{Antonio C\'ordoba} \\
{\small Departamento de Matem\'aticas}\\{\small Facultad de
Ciencias} \\ {\small Universidad Aut\'onoma de Madrid}
\\ {\small Crta. Colmenar Viejo km.~15,  28049 Madrid,
Spain} \\ {\small Email: antonio.cordoba@uam.es}
\end{tabular}
\end{quote}
\begin{quote}
\begin{tabular}{ll}
\textbf{Diego C\'ordoba} &  \textbf{Francisco Gancedo}\\
{\small Instituto de Ciencias Matem\'aticas} & {\small Department of Mathematics}\\
{\small Consejo Superior de Investigaciones Cient\'ificas} & {\small University of Chicago}\\
{\small Serrano 123, 28006 Madrid, Spain} & {\small 5734 University Avenue, Chicago, IL 60637}\\
{\small Email: dcg@imaff.cfmac.csic.es} & {\small Email: fgancedo@math.uchicago.edu}
\end{tabular}
\end{quote}


\begin{thebibliography}{99}

\bibitem{Ambrose} D. Ambrose. Well-posedness of Two-phase
Hele-Shaw Flow without Surface Tension. \emph{Euro. Jnl. of
Applied Mathematics} 15 597-607, 2004.

%\bibitem{AM} D. Ambrose and N. Masmoudi. The zero surface tension limit of two-dimensional water waves.
%\emph{Comm. Pure Appl. Math.} 58 1287-1315, 2005.

\bibitem{bear} J. Bear, Dynamics of Fluids in Porous Media, \emph{American
Elsevier}, New York, 1972.

\bibitem{BMO} G. Baker, D. Meiron and S. Orszag. Generalized vortex methods for free-surface flow problems.
\emph{J. Fluid Mech.} 123 477-501, 1982.

%\bibitem{bertozzi-Constantin} A.~L. Bertozzi and P. Constantin. Global regularity for vortex patches.
%\emph{Comm. Math. Phys.} 152 (1): 19--28, 1993.

%\bibitem{bertozzi-Majda} A.~L. Bertozzi and A.~J. Majda. \newblock Vorticity
%and the Mathematical Theory of Incompresible Fluid Flow. \newblock \emph{Cambridge Press}, 2002.

%\bibitem{Birkhoff} G. Birkhoff. Helmholtz and Taylor instability.
%\emph{Hydrodynamics Instability, Proc. Symp. Appl. Math. XII
%A.M.S.}, 55--76, 1962.

%\bibitem{CO} R. Caflisch and O. Orellana. Singular solutions and ill-posedness for the evolution of vortex
%sheets. \emph{SIAM J. Math. Anal.} 20 (2): 293--307, 1989.

%\bibitem{Ch} J.Y. Chemin. Persistence of geometric structures in
%two-dimensional incompressible fluids. \emph{Ann. Sci. Ecole. Norm. Sup.} 26 (4): 517--542, 1993.

\bibitem{CaC}  L.A. Caffarelli and A. C\'{o}rdoba, Phase transitions: Uniform regularity of the intermediate layers.,
 \textit{J. Reine Angew. Math., } 593, (2006), 209-235.

\bibitem{Peter} P. Constantin and M. Pugh. Global solutions for small data to the
Hele-Shaw problem. \emph{Nonlinearity}, 6 (1993), 393 - 415.

%\bibitem{CDGKS} P. Constantin, T.F. Dupont, R.E. Goldstein, L.P. Kadanoff, M.J. Shelley and S.M. Zhou.
%Droplet breakup in a model of the Hele-Shaw cell. \newblock
%\emph{Physical Review E}, 47, 4169-4181, 1993.

%\bibitem{CMT} P.~Constantin, A.~J. Majda, and E.~Tabak. \newblock Formation
%of strong fronts in the 2-{D} quasigeostrophic thermal active scalar.
%\newblock \emph{Nonlinearity}, 7:1495--1533, 1994.

\bibitem{CC}  A. C\'{o}rdoba and D. C\'{o}rdoba, A pointwise estimate for
fractionary derivatives with applications to P.D.E., \textit{Proc.
Natl. Acad. Sci., } 100, 26, (2003), 15316-15317.

\bibitem{cor2} A. C\'{o}rdoba and D. C\'{o}rdoba. A maximum principle applied to Quasi-geostrophic
equations. \emph{Comm. Math. Phys.} 249 (2004), no. 3, 511--528.

\bibitem{DY3} A. C\'ordoba, D. C\'ordoba and F. Gancedo. Interface evolution: the full water wave problem in 2D.
Preprint 2008.

\bibitem{DY} D. C\'ordoba and F. Gancedo. Contour dynamics of incompressible 3-D fluids
in a porous medium with different densities. \emph{Comm. Math.
Phys.} 273, 2, 445-471 (2007).

\bibitem{DY2} D. C\'ordoba and F. Gancedo. A maximum principle for the Muskat problem for fluids with different densities. To appear in \emph{Comm. Math.
Phys.}, Preprint 2007.

\bibitem{Dahlberg} Dahlberg, Bjorn E. J. On the Poisson integral for Lipschitz and
$C\sp{1}$-domains. Studia Math. 66 (1979), no. 1, 13--24.

%\bibitem{ES} J. Escher and G. Simonett. Classical solutions for Hele-Shaw models
%with surface tension. \emph{Adv. Differential Equations},
%2:619-642, 1997.

\bibitem{FKP} R. A. Fefferman, C. E. Kenig and J. Pipher. The theory of weights and the Dirichlet problem for elliptic equations.  \emph{Ann. of Math.} (2)  134  (1991),  no. 1, 65--124.

\bibitem{Y} F. Gancedo. Existence for the $\alpha$-patch model and the QG sharp front in Sobolev spaces.
\emph{Adv. Math.}, Vol 217/6: 2569-2598, 2008.

\bibitem{H-S} Hele-Shaw.
\newblock \emph{Nature} 58, 34, 1898.

\bibitem{Hou} T.Y. Hou, J.S. Lowengrub and M.J. Shelley. Removing the Stiffness
from Interfacial Flows with Surface Tension. \emph{J. Comput.
Phys.}, 114: 312-338, 1994.


\bibitem{Muskat} M. Muskat. \newblock The flow of homogeneous fluids through porous media.
\newblock \emph{New York}, 1937.

%\bibitem{Nirenberg} L. Nirenberg. An abstract form of the nonlinear Cauchy-Kowalewski theorem.
%\emph{J. Differential Geometry}, 6 561-576, 1972.

%\bibitem{Nishida} T. Nishida. A note on a theorem of Nirenberg. \emph{J. Differential Geometry}, 12 629-633, 1977.

\bibitem{S-T} P.G. Saffman and Taylor.
\newblock The penetration of a fluid into a porous medium or Hele-Shaw cell containing a more viscous liquid.
\newblock \emph{Proc. R. Soc. London, Ser. A} 245, 312-329, 1958.

\bibitem{SCH} M. Siegel, R. Caflisch and S. Howison. Global
Existence, Singular Solutions, and Ill-Posedness for the Muskat
Problem. \emph{Comm. Pure and Appl. Math.}, 57: 1374-1411, 2004.

%\bibitem{Rodrigo} J.L. Rodrigo. On the Evolution of Sharp Fronts for the
%Quasi-Geostrophic Equation. \emph{Comm. Pure and Appl. Math.}, 58: 0821-0866, 2005.

%\bibitem{St2} E. Stein and G. Weiss. \newblock Introduction to Fourier Analysis
%on Euclidean spaces. \newblock \emph{Princeton University Press.}
%Princeton, NJ, 1971.

\bibitem{St3} E.~Stein. \newblock Harmonic Analysis. \newblock \emph{
Princeton University Press.} Princeton, NJ, 1993.

%\bibitem{Taylor} G. Taylor. The instability of liquid surfaces when accelerated in a direction perpendicular
%to their planes. I. \emph{Proc. Roy. Soc. London. Ser. A.} 201,
%192-196 1950.

%\bibitem{Wu} S. Wu. \newblock Well-posedness in Sobolev spaces of the full water wave problem in 2-D.
%\emph{Invent. math.} 130, 39-72 1997.

%\bibitem{Wu2} S. Wu. \newblock Well-posedness in Sobolev spaces of the full water wave problem in 3-D.
%\emph{J. Amer. Math. Soc.} 12, 445-495 1999.

\end{thebibliography}
\end{document}